\documentclass{article}

\usepackage[T1]{fontenc}
\usepackage[utf8]{inputenc}
\usepackage{csquotes}
\usepackage[british]{babel}

\usepackage{mathtools}

\usepackage{amsmath}
\usepackage{amssymb}
\usepackage{amsfonts}
\usepackage{dsfont} 
\usepackage{amsthm}
\usepackage{thmtools}
\usepackage{tabularx}

\usepackage{tabularray}

\usepackage[font=small]{caption}
\usepackage{subcaption}

\usepackage[left=3.4cm,right=3.4cm]{geometry} 

\usepackage{float}

\usepackage{xcolor} 

\usepackage{multirow} 
\usepackage{adjustbox} 
\usepackage{placeins} 

\usepackage[style=alphabetic,backend=bibtex, url=false, doi=false]{biblatex}
\usepackage{enumitem} 
\usepackage{array} 

\usepackage{hyperref}
\usepackage[capitalise]{cleveref}
\numberwithin{equation}{section}

\usepackage{authblk}

\title{Central limit theorems describing isolation by distance under varying population size}
\date{January 2025}
\author[1]{Rapha\"el Forien}
\author[2]{Bastian Wiederhold}
\affil[1]{BioSP, INRAE, 84914, Avignon, France}
\affil[2]{Department of Statistics, University of Oxford, OX1 3LB, Oxford, UK}
\DeclareMathOperator{\blangle}{\big\langle}
\DeclareMathOperator{\brangle}{\big\rangle}
\DeclareMathOperator{\Blangle}{\Big\langle}
\DeclareMathOperator{\Brangle}{\Big\rangle}

\DeclareMathOperator{\dimension}{d}

\theoremstyle{definition}
\newtheorem{definition}{Definition}
\numberwithin{definition}{section}

\newtheorem{assumption}[definition]{Assumption}

\theoremstyle{plain}
\newtheorem{theorem}[definition]{Theorem}
\newtheorem{proposition}[definition]{Proposition}
\newtheorem{corollary}[definition]{Corollary}
\newtheorem{conjecture}[definition]{Conjecture}
\newtheorem{lemma}[definition]{Lemma}

\newtheorem{remark}[definition]{Remark}

\begin{document}


\maketitle

\begin{abstract}
We derive a central limit theorem for a spatial $\Lambda$-Fleming-Viot model with fluctuating population size. At each reproduction, a proportion of the population dies and is replaced by a not necessarily equal mass of new individuals. The mass depends on the local population size and a function thereof. Additionally, as new individuals have a single parental type, with growing population size, events become more frequent and of smaller impact, modelling the successful reproduction of a higher number of individuals. From the central limit theorem we derive a Wright-Mal\'ecot formula quantifying the asymptotic probability of identity by descent and thus isolation by distance. The formula reflects that ancestral lineages are attracted by centres of population mass and coalesce with a rate inversely proportional to the population size. Notably, we obtain this information despite the varying population size rendering the dual process intractable. 
\end{abstract}


\section{Introduction} \label{ch:4-introduction}
\paragraph{Isolation by distance}
The genetic patterns of a scattered population are spatially correlated if individuals move across space and successfully reproduce at positions different to their origin. In many species, individual movements during their lifetime will not exhaust the complete range of the species due to environmental, behavioural or energetic constraints. New genetic information created by mutation and recombination is therefore diffusing over space at a finite rate. The genetic similarity of individuals separated in space will decrease as a function of their separation. This phenomenon is called isolation by distance and the patterns of genetic isolation will reflect the spatial gene flow in the population.

\paragraph{Wright-Mal\'ecot formula}
To quantify patterns of genetic isolation we need to find an appropriate scale of genetic similarity. The simplest measure for the genetic relatedness of two individuals is to ask whether they are identical or not. While in reality the genetic information of an individual from a non-clonal species rarely matches any other individual, considering only a part of the genetic information such as a single locus renders the measure applicable. We could then measure isolation by distance by the probability of identity.  

Wright and Mal\'ecot (\parencite{Wright}, \parencite{Malecot}) observed that the decay of the probability of identity with spatial separation is roughly exponential. This approximation holds for many models and is today named after them. The universality of the approximation reflects the fact that ancestral lineages, the sequences of spatial positions of the ancestors, are as jump processes well-approximated by Brownian motions over sufficiently large scales. To see the same type in two sampled individuals, tracing backwards-in-time ancestral lineages need to have coalesced before either acquiring mutations. In turn, the Brownian motions approximating their motion are described by Gaussian transition kernels leading to the observed exponential decay in relatedness.
Later works noticed that some of their assumptions were inconsistent (\cite{felsenstein}) and reproduced the results more convincingly in space (\parencite{steppingstonemodel}, \parencite{neutralevolution2002}). 

\paragraph{Natural forces}
An exact exponential decay in real populations is unlikely, due to influences such as the local structure, rare events or inhomogenous dispersal dynamics.
Ideally, we would like to be able to describe isolation by distance under the multifarious forces acting on natural populations to provide a tool to infer information about the population from its genetic patterns. Simulation-based approaches have incorporated different influences on isolation by distance with recent examples found in (\parencite{Ringbauer2018}, \parencite{IBDSimulation1}).

A first step towards an analytical description of new influences was undertaken with the introduction of the class of spatial $\Lambda$-Fleming-Viot processes (\parencite{anewmodel}, \parencite{anewmodelsurvey}). In this framework, it was shown how long-range reproduction events might alter the genealogy. Several ancestral lineages of sampled individuals can merge at the same time, as the contained ancestors were offspring to the same parent of a long-range event. This is accompanied by a slower than exponential decay of the probability of identity.
Corresponding long-range Wright-Mal\'ecot formulae were discovered in \parencite{weissman}, \parencite{forien} and \parencite{forienwiederhold}. Here, ancestral lineages behave according to stable processes, which can be thought of as a generalization of Brownian motion with larger jumps. Moreover, even if ancestral lineages are still separated in space, there exists a positive rate of coalescence due to large spatial events.

\paragraph{Aim and motivation of this work}
In this paper, we derive a Wright-Mal\'ecot formula to quantify isolation by distance under the influence of varying population size. To this end, we propose a spatial $\Lambda$-Fleming-Viot type model with a reproduction mechanism which might change the local population density. In reproduction events, a proportion of the population in a region dies and is replaced by a not necessarily equal mass governed by a so-called growth function that depends on the local population size. If the carrying capacity is larger, more individuals are expected to have reproductive success. We model this by letting the frequency of events increase with higher population size while reducing the relative fraction replaced by offspring in the population. Our motivation is twofold:
\begin{enumerate}
\item
The resulting formula is more realistic than previous results, as it explains for example the increased genetic isolation if the population is separated through barriers to gene flow such as by unfavourable terrain (see \Cref{fig:1driftcomparison}). This is an important step towards bridging the gap between the aforementioned analytical work and more applied papers incorporating precisely these effects such as \parencite{Ringbauer2018}. Additionally, although in this paper the bias in the gene flow results from differences in population size, we expect a similar formula to hold if there are general biases in the displacement of individuals. For example, octopus dispersal patterns are more likely to follow the circular drift of the Antarctic circumpolar current (see \cite{antarctic23}).
\item 
Mathematically, the varying population size makes it more difficult to find duality relationships with link the forwards-in-time process with a corresponding backwards-in-time structure. Ancestral lineages might be non-Markovian and correlated in a complicated way, as the paths of ancestral lineages might contain information about the population profile. In such a setting, only few rigorous mathematical results have been obtained so far (\parencite{BirknerRandomWalks}, \cite{EKLRT24}, \cite{birkner2024quenched}). 
Notably, in \cite{EKLRT24} the motion of ancestral lineages is extracted from a complex spatial model using lookdown constructions. We will arrive at a comparable motion of ancestral lineages from a completely different angle using a strategy based on a central limit theorem (\cite{forien}). 

In previous papers (\cite{forien}, \cite{forienwiederhold}), the considered process was of population size one allowing one to derive the probability of identity and associated behaviour of ancestral lineages in a rather straight-forward manner from a corresponding central limit theorem. In the present paper, we will consider a central limit theorem originating from a process which is of varying population size. This renders the derivation of the formula for the probability of identity more complex, but demonstrates the feasibility of the technique in a much wider context.
\end{enumerate}

\paragraph{Forwards-in-time central limit theorem}
The main tool of this work will be a connection between the probability of identity and central limit theorems.  
Individuals will be equipped with neutral genetic types, tracing related individuals in the population. Mutation then corresponds to assuming uniform new types from the real interval $[0,1]$ with a neutral mutation rate $\mu$. As two uniform samples on the real interval $[0,1]$ are never equal, in this way, two individuals are only of the same type if neither has experienced mutation since their most recent common ancestor.

The process describing the limiting fluctuations of a spatial model is given as the solution of an SPDE driven by a Gaussian random field. A spatial differential operator such as the Laplacian disperses the `spikes' of the Gaussian random field across space. 
The genetic fluctuations of the SPDE in two distinct spatial positions can only be correlated if individuals from those positions had a common ancestor back in the past. This explains the emergence of the coalescence behaviour of the ancestral lineages of two samples in the central limit theorem. 
The differential operator carries the fluctuations across space and is therefore linked to the behaviour of ancestral lineages. A coalescence will cause a significant deviation from the average behaviour, and hence the noise and its spatial correlations reflect the coalescence behaviour. 

This work illustrates that the correspondence between central limit theorems and the probability of identity and the behaviour ancestral lineages is a general phenomenon, while providing more realistic analytical expressions to approximate isolation by distance.

\paragraph{Outline} The next section introduces and motivates our model. In \Cref{sec:results} we state and interpret our main results. In \Cref{sec:convergenceofpopulationsize} we deduce the convergence of the population size to the solution of a reaction-diffusion equation. \Cref{sec:martingaleproblems} starts analysing the evolution of genetic types under the model by stating and scaling the corresponding martingale problems. The proof of the central limit theorem is contained in \Cref{sec:proofclt}, and in \Cref{sec:proofofwmf} we deduce the corresponding Wright-Mal\'ecot formula. In \Cref{appendix:simulation} we will provide details on \Cref{fig:1driftcomparison}, which illustrates the accuracy of the predicted probability of identity by comparing to simulations. Further appendices deal with technical details.

\paragraph{Acknowledgements}
Rapha\"el Forien was supported in part by the Chaire Modélisation Mathématique et Biodiversité of Veolia-\'Ecole Polytechnique-Muséum National d'Histoire Naturelle-Fondation X.
Bastian Wiederhold was supported by the Engineering and Physical Sciences Research Council Grant [EP/V520202/1].

\section{Definition of the model}
\paragraph{Class of models} Our analysis will use a spatial $\Lambda$-Fleming-Viot (SLFV) model, a class of measure-valued stochastic processes first developed in \parencite{Eth08} and \parencite{anewmodel}. Up to null-sets, at every space-time position the process assigns a distribution on a type space, representing the composition of the population at this point. The population is driven by a Poisson point process that determines the time and location of events, the region affected and the `impact' $u$ of the event. Inside the area affected by an event, first a parental type is sampled from the population, and then a proportion $u$ of the population is replaced by this type. 

Varying population size SLFV models were first defined rigorously in \parencite{EK19}. There a general offspring mechanism was considered based on the assumption that the `mass' of offspring is random and has an expectation given by the affected area multiplied by the event impact. This idea is also present in our more specific choice of reproduction events. Additionally, events fall at a rate proportional to the local population size and the replaced proportion of the population decreases if the population size is larger.
Another modification of our model compared to the classical SLFV model is a mutation mechanism, which was first studied in \parencite{VW15}. Although the population is already an infinite-population limit, it can be thought of as each individual mutating with a certain rate $\mu$ to a new type represented by a uniform choice from the real interval $[0,1]$. To guarantee existence of the model, we will work on the $d$-dimensional torus $\mathbb{T}^d$.

\paragraph{Choice of events} Throughout this work, we will denote the population size by $n = n_t (z)$ (defined almost everywhere) and the local average of the population size over a ball $B(x,R)$ with radius $R > 0$ as 
\begin{equation*}
\overline{n} (x, R) := \frac{1}{V_R} \int_{B(x,R)} n_t(z) dz.
\end{equation*}
Here, $V_R$ denotes the volume of the ball $B(x,R)$. Heuristically, during a time period of length one, we wish to replace a certain proportion $u$ of the population reflecting the turnover through death and reproduction. However, if the local population size $\overline{n}$ is larger, we expect more individuals have reproductive success and to be parental to some part of the proportion $u$ of new individuals. There are two different ways of implementing this in the SLFV:
\begin{enumerate}
    \item 
    Either we have events replacing a proportion $u \in (0,1)$, but allow a number of parents proportional to the local population size $\overline{n}$,
    \item
    or we could have the rate of events proportional to the local population size $\overline{n}$, each with a single parent and replacing only a proportion $u/\overline{n}$ of the population.
\end{enumerate}
We opted for the second possibility. Events will cover a ball $B(x,R)$ of radius $R> 0 $ around the centre $x$ of the event.
To be exact, for each spatial position $y \in B(x,R)$ the population size changes according to
    \begin{equation} \label{eq:choiceofevents}
    n (y) \mapsto \Big(1 - \frac{u}{\overline{n} (x, R) + 1} \Big) n (y) + \frac{u}{\overline{n} (x,R) + 1} \Big( 1 + r_x( \overline{n} (x, R)) \Big) \overline{n} \text{ at rate } \overline{n} + 1,
    \end{equation}
where $r_x : (0, \infty) \rightarrow [-1, \infty)$ is a function capturing the regulation of the population size depending on the spatial location $x \in \mathbb{R}^d$. 
Typically, the growth function $r_x (\overline{n})$ will be positive if the local population size $\overline{n}$ is smaller than the carrying capacity of the location and negative if the local population size is larger. 

One might note that in \eqref{eq:choiceofevents} we accelerate by $\overline{n} + 1 $ instead of $\overline{n}$ and replace a proportion of $u / (\overline{n} +1)$ instead of $u/ \overline{n}$. Adding one to the rate of reproduction events captures that there is always reproductive activity even if the local population density is low. The same replacement in the denominator also prevents the population size from becoming negative, as otherwise for small $\overline{n}$ the term $(1 - u / \overline{n})$ could be negative. 

To provide a precise definition, we will now first clarify the state space $ \Xi $ of the process and then describe its evolution (similar to \parencite[Section 1]{forien}).

\paragraph{State space} We will choose the growth function in such a way that the population size does not pass a threshold $n_{\max} > 0$. The type space will be the real interval $[0,1]$. At any location, the population can be thought of assuming values in $\mathcal{M}_{n_{\max}} ([0,1])$, the set of all measures with mass smaller or equal to $n_{\max}$ on $[0,1]$. The full state space $\Xi = (\Xi', d)$ will be
\begin{equation*}
\Xi' = \Big\{ \rho : \mathbb{T}^d \rightarrow \mathcal{M}_{n_{\max}} ( [0,1] ) \text{ is a measurable function} \Big\} \Big\vert_\sim
\end{equation*}
under the equivalence relationship
\begin{equation*}
\rho_1 \sim \rho_2 \Leftrightarrow \{ z \in \mathbb{T}^d : \rho_1 (z) \neq \rho_2 (z) \} \text{ is a Lebesgue null set}.
\end{equation*}
We equip $\Xi'$ with the topology of weak convergence. To this end, let $(\phi_n)_{n \geq 1}$ be a sequence of test functions $\phi_n : \mathbb{T}^d \times [0,1] \rightarrow \mathbb{R}$ such that
\begin{equation*}
    \exists k : \langle \rho_1, \phi_k \rangle \neq \langle \rho_2 , \phi_k \rangle \Leftrightarrow \rho_1 \neq \rho_2 \text{ in } \Xi'.
\end{equation*}
The angular bracket $\langle \rho, \phi \rangle$ denotes the integration $\int \phi(z,k) \rho (z, dk) dz$ of $\phi$ with respect to $\rho$ and will be used frequently throughout this work.
To simplify some of the later calculations, in the same way as \parencite{forien}, we introduce $p$-norms as
\begin{equation} \label{eq:pnormdefinition}
\Vert \phi \Vert_p := \Bigg(\int_{\mathbb{T}^{\dimension}} \sup_{k \in [0,1]} \vert \phi (z, k) \vert^p dz \Bigg)^{1/p}.
\end{equation}
and additionally assume a uniform bound of the partial derivatives for indices $\kappa \in \mathbb{N}^d$
\begin{equation} \label{eq:assumption_on_test_sequence}
 \sup_{n \in \mathbb{N}} \Big\{ \Vert \partial_\kappa \phi_n \Vert_p : \vert \kappa \vert \leq 2 , p \in \{ 1, \infty \} \Big\} < \infty.
\end{equation}
The weak metric $d$ inducing the topology is given by
\begin{equation} \label{eq:vaguemetric}
	d (\rho_1, \rho_2) := \sum_{n = 1}^\infty \frac{1}{2^n} \big\vert \langle \rho_1, \phi_n \rangle - \langle \rho_2 , \phi_n \rangle \big\vert.
\end{equation}
In the same way as for the local population size $\overline{n}$, we will write
\begin{equation*}
\overline{\phi} (x,k,R) := \frac{1}{V_R} \int_{B(x,R)} \phi (y,k) dy 
\end{equation*}
for the average of a function $\phi : \mathbb{T}^d \times [0,1] \rightarrow \mathbb{R} $ over a ball of radius $R$ with volume $V_R$. Suppose a stochastic process evolves in $\Xi$. For a space-time position $(z,t)$, the value $\rho_t(z)$ can be interpreted as assigning masses to genetic types of the interval $[0,1]$ (up to a Lebesgue null-set). The total mass \[n_t(z) = \int_0^1 \rho_t(z,dk) := \mathcal{H} \rho_t(z)\] encodes the population size.
We will use the symbol $\mathcal{H}$ to denote the mass of a measure $M \in \mathcal{M} ([0,1])$ or a function $\phi : [0,1] \rightarrow \mathbb{R}$ if interpreted as a density, i.e.
\begin{equation} \label{eq:massoperator}
\mathcal{H} M := \int_0^1 M(dk), \hspace{1cm} \mathcal{H} \phi := \int_0^1 \phi (k) dk.
\end{equation}
\paragraph{Restrictions on the growth function}
We wish to make three additional assumption on our model, which we can guarantee by assuming certain properties of the growth function. Recall from \eqref{eq:choiceofevents} that we add a mass of size
\[
\frac{u}{\overline{n} (x,R) + 1} \Big( 1 + r_x( \overline{n} (x, R)) \Big) \overline{n}
\]
of the sampled parental type in reproduction events. 
\begin{enumerate}
\item 
It is natural to assume that the new-born mass in events is positive. Strictly speaking, the above term $1 + r_x(\overline{n})$ requires the growth function to be bounded from below by minus one. As we will scale the growth function by a term converging to zero \eqref{eq:scaling1}, a lower bound will suffice to ensure this property for the scaled processes.
\item 
We wish to ensure upper and lower bounds on the population size $0 < n_{\min} \leq n_{\max}$. If the growth function $r_x (\overline{n})$ is positive, the population grows locally. To ensure a lower bound it is therefore enough to assume that as the local population size $\overline{n}$ becomes small, the growth function stays larger than a positive value for all locations. An upper bound can be guaranteed by assuming that $(1 + r_x (\overline{n} ))\overline{n}$ stays bounded. 

These assumptions significantly reduce the technicalities required for the proofs in this paper. For example, if the population size had no upper bound, one would need to establish a control on large values in order to ensure convergence of the sequence of rescaled population size processes. 
\item
In addition, we will make a regularity assumption on the growth function. This allows us to follow the arguments of \parencite{forien_central_2017} for the convergence of the population size and ensures that the limit, which will be a solution to a reaction-diffusion equation, exists. 
\end{enumerate}
Each of these ideas can be formalized as follows.
\begin{assumption} \label{assumption}
We assume the following properties of the growth function $r: \mathbb{T}^d \times \mathbb{R}_+ \rightarrow \mathbb{R}_+$:
\begin{enumerate}
    \item (Positive birth) There exists a lower bound $C \in \mathbb{R}$ on the growth function such that $\forall x \in \mathbb{T}^d$ and $\forall n > 0$
    \[ r_x(n) \geq C. \]
    If we consider the model without rescaling, we assume $C = -1$.
    \item (Upper and lower bound on the population size)
  There exists $C' > 0$ such that $\forall x \in \mathbb{T}^d$ we have
    \begin{equation*}
        \lim_{n \to 0+} r_x(n) > C'
    \end{equation*}
    and there exists $n_{\max} > 0$ such that all population sizes $n > 0$ and all locations $x \in \mathbb{T}^d$
    \begin{equation*}
     (1 + r_x(n)) n < n_{\max}.
    \end{equation*}
    \item (Regularity of the growth function) Assume that for all population sizes $ n > 0$ the function $r_\cdot (n) \in C(\mathbb{T}^d)$ and that
    \begin{equation*}
       \sup_{x \in \mathbb{T}^d} \sup_{n \geq 0} \sup_{0 \leq \vert \kappa \vert \leq 4} \Big\vert \frac{\partial^\kappa}{\partial n^\kappa} r_x (n) \Big\vert < \infty.
    \end{equation*}
\end{enumerate}
\end{assumption}

\begin{definition} \label{def:mSLFV}
Let $\mu, R$ be positive constants, $u \in (0,1)$, $r : \mathbb{T}^d \times \mathbb{R}_+ \rightarrow \mathbb{R}_+$ a growth function satisfying \Cref{assumption} and $ \Pi $ be a Poisson point process on $ \mathbb{R}_+ \times \mathbb{T}^d \times [0,1]$ with intensity $ (n_{\max} +1) ( dt \otimes dx \otimes dw) $.

	We define the modified spatial $ \Lambda $-Fleming-Viot process (mSLFV) as the stochastic process $ \rho = (\rho_t, t \geq 0) $ with population size $n_t (z) = \mathcal{H} \rho_t(z)$, which evolves in $\Xi$ according to: 
\begin{enumerate}
\item Reproduction: For any point $(t,x,w) \in \Pi$ we perform the following steps. In case $w \geq (\overline{n}_t (x) + 1) / (n_{\max} +1)$, we ignore the event. Otherwise, we choose a parental type $ k_0 $ according to the normalized mass distribution in the ball $B(x,R)$, i.e. with respect to
\begin{equation*}
    \frac{1}{\int_{B(x, R)} n_{t-} (z) dz} \int_{B(x,R)} \rho_{t-} (y, dk) dy.
\end{equation*}
  Then we remove a proportion 
  \[ u' := \frac{u}{\overline{n} (x, R) +1}\]
  of the population inside the affected area $ B(x,R) $ and introduce a mass 
  \[
  u' \Big(1 + r (\overline{n} (x, R)) \Big) \overline{n} (x, R)
  \] 
  of offspring of type $k_0$, i.e. for every $z \in B(x,R)$ we set
\begin{equation*}
        \rho_t (z, dk) = (1 - u') \rho_{t-} (z, dk) + u' \Big(1 + r_x (\overline{n} (x, R)) \Big) \overline{n} (z, R) \delta_{k_0} (dk).
    \end{equation*}
    \item
    Mutation: If $t_1$ and $t_2$ are two consecutive times at which a spatial location $ z\in \mathbb{T}^d$ is affected by events, we set for times $s$ with $ t_1 \leq s < t_2 $,
    \begin{equation*} 
        \rho_s (z, dk) = e^{- \mu (s -t_1)} \rho_{t_1} (z, dk) + n_{t_1} (z) \left(1 - e^{- \mu (s- t_1 )}\right) dk.
    \end{equation*} 
\end{enumerate}
\end{definition}
An interpretation of the mutation mechanism was given in \Cref{ch:4-introduction}
If there is no ambiguity, we will abbreviate $\overline{n}(z, R) $ by $\overline{n}$. 
As the core $D^{core}(\mathcal{G})$ of the domain $D(\mathcal{G})$ of the generator $\mathcal{G}$ of the mSLFV process, we use the set of functions
\begin{equation*}
D^{core} (\mathcal{G} ) := \Big\{ F_{f,\phi} : \Xi \rightarrow \mathbb{R}, \rho \mapsto f (\langle \rho, \phi \rangle ) \enspace : \enspace \phi \in C_c(\mathbb{T}^d), f \in C^1 (\mathbb{R}) \Big\}.
\end{equation*}
The generator of the process acting on functions in $D^{core}(\mathcal{G})$ is
\begin{equation} \label{eq:generator}
    \begin{aligned}
    &\mathcal{G} F_{f, \phi} (\rho)\\ 
    &= \mu \Big\langle \rho, \int_{[0,1]} ( \phi (\cdot,k') - \phi (\cdot,k)) dk' \Big\rangle f'(\langle \rho, \phi \rangle) \\
& \hspace{1cm} + \int_{\mathbb{T}^d}  \frac{1}{\int_{B(x,R)\times [0,1]} \rho(y,dk) dy} \int_{B(x,R) \times [0,1]}  \\
& \hspace{2cm} (\overline{n} + 1) \Big( f \Big( \blangle \rho + \frac{u}{\overline{n} +1} \mathds{1}_{B(x,R)} (\delta_{k_0} [1 + r(\overline{n}) ] \overline{n} - \rho ) , \phi \brangle \Big)  - f(\langle \rho, \phi \rangle) \Big) \\
& \hspace{11cm} \rho (y, dk_0) dy dx.
    \end{aligned}
\end{equation}
We are working on the torus $\mathbb{T}^d$ and our assumptions ensure a finite total rate of events.
\begin{proposition}
Under \Cref{assumption}, the martingale problem associated to $(\mathcal{G}, D(\mathcal{G}))$ is well-posed, i.e. there exists a unique solution.
\end{proposition}

\begin{remark} \label{rem:otherassumptions} 
The model could be generalized; this remark presents some of the possibilities. 

As we will see in the next section, ancestral lineages follow the gradient of the population size. This vector field is irrotational. There are situations such as the aforementioned octopus population subject to the circular drift around Antarctica \cite{antarctic23}, where one can imagine rotational biases of the ancestral lineages. In general, an additional vector field might be used to describe biases in the dispersal patterns. One could implement this in the model, for example, by incorporating a bias in the sampling location of the parental type.

As already mentioned, weakening the existence of $n_{\max}$ would require more work to control how rarely large population sizes are attained to prove e.g. tightness of the sequence of processes under scaling. Removing the assumed lower bound $n_{\min}$ on the population size, would complicate any bound used on the inverse of the population.

Currently, as we fix the impact $u$ and the radius $R > 0$, given the population profile $n$ and event centre $z \in \mathbb{T}^d$ the mass added in an event is deterministic. We could incorporate a random factor $\xi$ with $\mathbb{E} (\xi) = 1$ in the mass added in reproduction.
\end{remark}

\begin{remark}
Working on a torus $\mathbb{T}^d$ ensures uniqueness of the process, as events then occur at a finite rate. Existence of the process can be guaranteed even on $\mathbb{R}^d$ using a sequence of approximating processes. Uniqueness of a varying population size SLFV on $\mathbb{R}^d$ has yet to be shown. To our knowledge uniqueness for SLFV processes has been obtained so far using either duality with a simpler process or lookdown constructions, where the process is constructed as the De-Finetti measure of an exchangeable sequence. The varying population size renders the dual process intractable in our case. A potential lookdown construction of our model could be based on \parencite[Section 4.2]{EK19}. Unfortunately, the varying population size leads to a particle representation in which individuals change level in every event (so at an infinite rate), which prevents a uniqueness proof. 
\end{remark}

\section{Main results} \label{sec:results}
The centrepiece of this work is a central limit theorem for the process. We will see that it allows us to calculate the probability of identity by descent and reflects the behaviour of ancestral lineages. These results generalize the findings of the fixed radius case of \parencite{forien} to varying population size.
\subsection{Law of large numbers for the population size}
We will apply the following scaling to the process. As $N \to \infty$, we will let the impact $u \to 0$. This corresponds to letting the strength of the genetic drift become smaller and smaller, which in turn emulates a larger population size (even though the model is already an infinite-population limit). As `individuals' participate in events at a rate proportional to the impact, we will accelerate time by a scaling inverse to the scaling of the impact. In addition, we rescale space and time in the usual way to obtain a Brownian motion from a jump process. The mutation rate requires to be scaled inversely to time, so that mutation remains of an order visible in the limit. Equally, the growth function needs to be scaled down to behave reasonably as we accelerate time.

Under this scaling the population size will converge to the solution $n_t$ of a reaction-diffusion equation, while the typed process $\rho_t^N$ converges to Lebesgue measure $\lambda$ on $\mathbb{T}^d \times [0,1]$ multiplied by the limiting population size $n_t$. Lebesgue measure appears as the mutation mechanism dominates in the limit. Equivalently, we could view $\lambda$ as the element of $\Xi$ associating to almost every position a uniform distribution on $[0,1]$.

More formally, let $\delta_N$ be an arbitrary sequence of positive numbers converging to zero. We define the rescaled process $\rho_t^N(z,dk)$ in two steps. 
\begin{enumerate}
\item 
    First, we let $\mathfrak{p}^{N} = (\mathfrak{p}_t^{N}, t \geq 0)$ be the mSLFV process with rescaled parameters and growth function
\begin{equation} \label{eq:scaling1}
u_N := \frac{u_0}{N}, \hspace{1cm} \mu_N := \delta_N^2 \frac{\mu}{N}, \hspace{1cm} r_{z}^N := \delta_N^2 r_{\delta_N z}
\end{equation}
for $z \in \delta_N^{-1} \mathbb{T}^d$.
\item
Then, for $z \in \mathbb{T}^d$, we apply the spatial and temporal scaling
\begin{equation} \label{eq:scaling2}
\rho_t^N (z,dk) := \mathfrak{p}_{Nt/\delta_N^2}^{N} \Big( \frac{z}{\delta_N} , dk \Big).
\end{equation}
\end{enumerate}

As we see in \Cref{lem:mppop} the evolution of the rescaled population size can be characterised by a martingale problem. The parental choice and subsequent dispersal of offspring results in the double average operator
\begin{equation} \label{eq:jumpgenerator}
\begin{aligned}
\mathcal{L}^N \psi (z) &:=  \frac{u V_R \delta_N^{-2}}{V_{\delta_N R}^2} \int_{B(z,\delta_N R)} \int_{B(x,\delta_N R)} \Big[ \psi (y) - \psi (z) \Big] dy dx 
= u V_R \delta_N^{-2} \big( \overline{\overline{\psi}} (z, \delta_N R) - \psi (z) \big),
\end{aligned}
\end{equation}
whereas the growth of the population size will be described by
\begin{equation} \label{eq:growthgenerator}
\begin{aligned}
\mathcal{R}_n^N \psi (z) &:= u V_R \frac{1}{V_{\delta_N R}^2} \int_{B(z,\delta_N R)}  r_x(\overline{n} (x, \delta_N R)) \int_{B(x,\delta_N R)} \psi (y) dy dx 
= u V_R \overline{r(\overline{n}) \overline{\psi}} (z, \delta_N R).
\end{aligned}
\end{equation}

\begin{lemma} \label{lem:mppop}
Let $\phi \in L^1 (\mathbb{T}^d)$. Then, 
\begin{equation*}
Y_t^N (\phi) := \blangle n^N_t , \phi \brangle - \blangle n_0^N, \phi \brangle - u V_R \int_0^t \blangle n_s^N, (\mathcal{L}^N + \mathcal{R}_{n_s^N}^N ) \phi \brangle ds
\end{equation*}
is a square-integrable martingale with predictable quadratic variation
\begin{equation} \label{eq:quadratic_variation_pop}
\begin{aligned}
    &\blangle Y^N (\phi) \brangle_t\\
    &= \int_0^t \frac{u^2}{\delta_N^{d+2}} \int_{\mathbb{T}^d} \frac{1}{(\overline{n_s^N} (x, \delta_N R) +1)N} \\
&\hspace{1.5cm} \Bigg( \int_{\mathbb{T}^d} \mathds{1}_{B(x,\delta_N R)} (z) \Big([1 + \delta_N^2 r_x(\overline{n_s^N} (x, \delta_N ) ]\overline{n_s^N} (x, \delta_N R) -n_s^N (z) \Big) \phi (z) dz \Bigg)^2  dx ds.
\end{aligned}
\end{equation}
\end{lemma}
The proof of \Cref{lem:mppop} can be found at the beginn of \Cref{sec:convergenceofpopulationsize}. In the limit as $N \to \infty$, the first operator $\mathcal{L}^N$ will result in a Laplacian, and the growth operator $\mathcal{R}_n^N$ in a reaction term so that the population size converges to the solution of a reaction-diffusion equation.

\begin{restatable}{proposition}{theollnpop} \label{prop:llnpop}
Fix $T > 0$, and let $n_t^N$ denote the population size of the rescaled process $\rho_t^N$. Then,
\begin{equation}
\lim_{N \to \infty} \mathbb{E} \Bigg[ \sup_{t \in [0, T]} d \Big( n_t^N, n_t \Big) \Bigg] = 0,
\end{equation}
where $n_t$ solves the reaction-diffusion equation
\begin{equation} \label{eq:reactiondiffusion}
\frac{d n_t}{dt} = u V_R \Big( \frac{R^2}{d + 2} \Delta n_t + r (n_t) n_t \Big),
\end{equation}
and $d$ is the metric defined in \eqref{eq:vaguemetric}.
\end{restatable}
The proof of this theorem is contained in \Cref{sec:convergenceofpopulationsize}. It follows from \parencite[Section 3.2-4.3]{forien_central_2017}, but we included it as the preliminary result \Cref{lem:boundonFsquared} plays a crucial role in the later proof of the central limit theorem.
\subsection{Central limit theorem}
We will consider the fluctuations of $\rho_t^N$ around Lebesgue measure on $\mathbb{T}^d \times [0,1]$ multiplied by a deterministic approximation $m_t^N$ to the population size $n_t^N$
\begin{equation} \label{eq:fluctuations}
Z_t^N := \sqrt{N \delta_N^{2 - d}} \Big( \rho_t^N - \lambda m_t^N \Big).
\end{equation}
Underlying this scaling is the idea that the process $\rho_t^N$ converges to $\lambda n_t$ - although we will not explicitly formulate this as a result. Naturally, we assume that the sequence $\delta_N$ is chosen in such a way that $N \delta_N^{2-d} \to \infty$ as $N \to \infty$. In fact, for the proof of \Cref{lem:boundflucpop} and \Cref{lem:compactcontainment}, we will assume even
\begin{equation} \label{eq:assumptionondelta}
N^{\frac{1}{2}} \delta_N^{1 + \frac{d}{2} } \xrightarrow[N \to \infty]{} \infty.
\end{equation}
The deterministic approximation to the population size that we take is defined by
\begin{equation} \label{eq:middleterm}
\blangle m_t^N, \phi \brangle := \blangle m_0^N , \phi \brangle + \int_0^t \blangle m_s^N , (\mathcal{L}^N + \mathcal{R}_{m_s^N}^N ) \phi \brangle ds.
\end{equation}
Comparing to \Cref{lem:mppop}, this formula differs from a representation of $n_t^N$ by a missing martingale term which would contain the randomness.
The scaling \eqref{eq:fluctuations} is the same as in specific settings of \parencite{forien} and \parencite{forienwiederhold}. 
The term $N \delta_N^{2-d}$, the square of $\sqrt{N \delta_N^{2-d}}$, will compensate a corresponding term in the predictable quadratic variation term \eqref{eq:quadratic_variation_pop} of \Cref{lem:mppop}, which results form a change of variables of \eqref{eq:quadratic_variation_pop}.

We denote by $\mathcal{D} (\mathbb{T}^d \times [0,1])$ functions $\phi: \mathbb{T}^d \times [0,1] \rightarrow \mathbb{R}$ which are smooth in the spatial component. Let $\mathcal{D}'(\mathbb{T}^d \times [0,1])$ be the dual space of generalized functions and $\mathbb{D}(\mathbb{R}_+, \mathcal{D}' (\mathbb{T}^d \times [0,1]))$ the space of c\`adl\`ag paths in the dual space. Recall from \eqref{eq:massoperator} that $\mathcal{H} Z_t (\cdot) := \int_0^1 Z_t (\cdot, dk)$.
\begin{theorem}\label{theo:clt}
The sequence of fluctuations $( Z_t^N, t \geq 0 ) , N = 1,2,\ldots$ converges in distribution in $D (\mathbb{R}_+, \mathcal{D}' ( \mathbb{T}^d \times [0,1] ))$ to the mild solution $(Z_t, t \geq 0)$ of the stochastic partial differential equation
\begin{equation} \label{eq:spde}
\left\{
\begin{aligned}
dZ_t &= \Bigg[ u V_R \Bigg( \frac{R^2}{d+2} \Delta + r(n_t) \Bigg) Z_t - \mu Z_t + \big(\mu + u V_R n_t r'(n_t)\big) \mathcal{H} Z_t  \Bigg] dt + d W_t, \\
Z_0 &= 0.
\end{aligned}
\right.
\end{equation}
The Gaussian random field $W$ is uncorrelated in time and its covariation measure $\mathcal{Q}_t$ acts on functions $\phi \in \mathcal{D} (\mathbb{T}^d \times [0,1])$ through
\begin{equation} \label{eq:limitingcovariation}
\blangle \mathcal{Q}_t , \phi \otimes \phi \brangle := u^2 V_R^2 \int_{\mathbb{T}^d} \frac{n_t (z)^2}{n_t(z) + 1} \int_{[0,1]} \Bigg( \phi (z,k ) - \int_{[0,1]} \phi (z,k') dk' \Bigg)^2 dk dz.
\end{equation}
\end{theorem}
We will prove the theorem in \Cref{sec:proofclt}. The arguments synthesize and generalize ideas of \cite{forien} and \cite{forien_central_2017}. 
\begin{remark}
The idea behind the classical central limit theorem is to recentre by the mean of the random variables involved. In fact, there is more freedom regarding the recentring, as long the sequence used to recentre converges to the `limiting mean'. The same freedom is present in the spatial central limit theorems considered here. In \Cref{eq:fluctuations}, recentring by $\lambda n_t^N$ would equally have been possible and would result in the above SPDE without the term involving $\mathcal{H} Z_t$. We chose to recentre by the deterministic approximation $\lambda m_t^N$ as the calculations will illustrate that the Wright-Mal\'ecot formula is the same with or without the terms involving $\mathcal{H} Z_t$, so invariant under the recentring. 
\end{remark}
\subsection{Wright-Mal\'ecot formula} \label{subsec:wmf}
Next, we explain how the fluctuations $Z_t^N$ can be related to the probability of identity. As the processes assume values in $\Xi$ and are only unique up to a spatial nullset, the probability of identity $P_t^N (\psi_1, \psi_2)$ will depend on two probability densities $\psi_1, \psi_2$ describing the sampling location. The probability of sampling the same type from two probability distributions $d_1, d_2$ is given by evaluating the distributions $d_1,d_2$ with respect to the indicator of the diagonal $ \mathds{1}_{\diagdown} (k_1, k_2) := \mathds{1}_{k_1 = k_2}$. Accounting for the varying population size, we hence define the probability of identity as
\begin{equation} \label{eq:individualbasedsampling}
\begin{aligned}
&P_t^N (\psi_1 , \psi_2) \\
&= \mathbb{E} \Bigg[ \frac{1}{\langle n_t^N , \psi_1 \rangle \langle n_t^N, \psi_2 \rangle} \int_{(\mathbb{T}^d)^2} \mathds{1}_{\Delta} (k_1, k_2)  \psi_1 (z_1) \psi_2 (z_2) \rho_t^N (z_1, dk_1) \rho_t^N (z_2, dk_2) dz_1 dz_2 \Bigg] \\
&= \mathbb{E} \left[ \frac{\blangle \rho_t^N \otimes \rho_t^N , (\psi_1 \otimes \psi_2) \mathds{1}_{\diagdown} \brangle}{\langle n_t^N , \psi_1 \rangle \langle n_t^N, \psi_2 \rangle} \right].
\end{aligned}
\end{equation}
We think of this sampling scheme as giving weight $\psi_1 (z)$ to all individuals in location $ z \in \mathbb{T}^d$, and then sampling an individual proportionally to these weights. The procedure is repeated with $\psi_2$ for the second sampled individual.

The question still remains how the probability of identity is related to the central limit theorem \Cref{theo:clt}. The connection is straightforward: We can utilize the relationship $ Z_t^N := \sqrt{N \delta_N^{2 - d}} \big( \rho_t^N - \lambda m_t^N \big) $ to rewrite \eqref{eq:individualbasedsampling} as
\begin{equation*} 
P_t^N (\psi_1, \psi_2) = \frac{1}{N \delta_N^{2-d}} \mathbb{E} \left[ \frac{\blangle Z_t^N \otimes Z_t^N , (\psi_1 \otimes \psi_2) \mathds{1}_{\diagdown} \brangle}{\langle n_t^N , \psi_1 \rangle \langle n_t^N, \psi_2 \rangle} \right].
\end{equation*}
This step uses that Lebesgue measure does not charge the diagonal and was first noted in \cite{forien}. Rearranging,
\begin{equation} \label{eq:samplingwithclt}
N \delta_N^{2-d} P_t^N (\psi_1, \psi_2) = \mathbb{E} \left[ \frac{\blangle Z_t^N \otimes Z_t^N , (\psi_1 \otimes \psi_2) \mathds{1}_{\diagdown} \brangle}{\langle n_t^N , \psi_1 \rangle \langle n_t^N, \psi_2 \rangle} \right].
\end{equation}
The above equation shows that the scaling factor $N \delta_N^{2-d}$ of the central limit theorem is precisely the scaling factor required for $P_t^N (\psi_1, \psi_2)$ to converge to a non-trivial limit. 

To use the central limit theorem, we would like to take the limit $N \to \infty$ on both sides of \eqref{eq:samplingwithclt} and exchange limit and expectation on the right-hand side to obtain
\begin{equation*}
\lim_{N \to \infty} N \delta_N^{2-d} P_t^N (\phi , \psi) = \mathbb{E} \left[ \frac{\blangle Z_t \otimes Z_t , (\psi_1 \otimes \psi_2) \mathds{1}_{\diagdown} \brangle}{\langle n_t , \psi_1 \rangle \langle n_t, \psi_2 \rangle} \right].
\end{equation*}
When the population size is constant, this can be achieved via an application of the dominated convergence theorem (see \cite{forien}).
As we will elaborate in \Cref{remark:technicality}, in the general case there remains a technicality resulting from the intrinsically varying population size, which prevented us from applying dominated convergence. Replacing the dependency of the growth function on the population size $n^N$ in \Cref{def:mSLFV} by a dependency on the deterministic approximation $m^N$ would eliminate this technicality. In this case, bounding the population size would become more involved as 
\[ \big( 1+ r(\overline{m_t^N}) \big) \overline{n_t^N}\]
could potentially grow to infinity under our current assumption. As this would not add to the biological interpretation of this work, we refrain from dealing with this additional technicality.
\begin{conjecture} \label{conj:swapexpectationandlimit}
For any $t \geq 0$, any $\psi_1, \psi_2 \in \mathcal{D} (\mathbb{T}^d)$, we have
    \begin{equation*}
    \lim_{N \to \infty} N \eta_N P_t^N (\psi_1, \psi_2) = \mathbb{E} \Bigg[ \frac{\langle Z_t \otimes Z_t , (\psi_1 \otimes \psi_2) \mathds{1}_{\diagdown} \rangle}{\langle n_t, \psi_1 \rangle \langle n_t, \psi_2 \rangle} \Bigg].
    \end{equation*}
\end{conjecture}
Our aim is now to compute the right-hand side. To avoid dealing with the fraction $\langle n_t , \psi_i \rangle$, we choose 
\[ 
\vartheta_i := \frac{n_t \psi_i}{\langle n_t, \psi_i \rangle},
\]
which transforms the quantity to
\[
\mathbb{E} \Bigg[ \frac{\langle Z_t \otimes Z_t , \mathds{1}_{\diagdown} \psi_1 \otimes \psi_2 \rangle}{\langle n_t, \psi_1 \rangle \langle n_t, \psi_2 \rangle} \Bigg] = \mathbb{E} \Big[ \Big\langle Z_t \otimes Z_t , \mathds{1}_{\diagdown} \Big( \frac{\vartheta_1}{n_t} \otimes \frac{\vartheta_2}{n_t} \Big) \Big\rangle \Big].
\]
In \Cref{sec:proofofwmf}, we will use this expression to prove the following theorem.
\begin{restatable}{theorem}{theoremWMF} 
\label{theo:varying_size_wmf}
Let $G_{s,t}$ be the fundamental solution associated to
\begin{equation*}
- \frac{d}{ds} G_{s,t} \psi = A_s^* G_{s,t} \psi,
\end{equation*}
where $A_s^*$ is the adjoint of 
\begin{equation} \label{eq:generatorwmflineage}
A_s \phi = \frac{u V_R R^2}{d + 2} \Big( \Delta \phi + 2\frac{\nabla n_s}{n_s} \cdot \nabla \phi \Big)
\end{equation}
and let $\vartheta_i$ be as above. Then the Wright-Mal\'ecot formula for two probability densities $\psi_1, \psi_2$ is given by 
\begin{equation} \label{eq:wmfformula}
\begin{aligned}
\mathbb{E} \Bigg[ \frac{\langle Z_t \otimes Z_t , \mathds{1}_{\diagdown} \psi_1 \otimes \psi_2 \rangle}{\langle n_t, \psi_1 \rangle \langle n_t, \psi_2 \rangle} \Bigg] 
&= u^2 V_R^2 \int_0^t \int_{\mathbb{T}^d} \frac{e^{- 2 \mu (t-s)}}{n_s(z) +1} G_{s,t} \vartheta_1 (z) G_{s,t} \vartheta_2 (z) dz ds.
\end{aligned}
\end{equation}
\end{restatable}
The function $\vartheta_i = \frac{n_t \psi_i}{\langle n_t, \psi_i \rangle}$ can be seen as being the limiting density according to which the sampling of positions occurs.

\begin{remark}
There is at least one more way in which one could define the probability of identity. Informally, we could instead sample two positions $z_1, z_2$ according to two the densities $\psi_1, \psi_2$. Then, the probability that two individuals sampled according to the probability measures $(n_t^N(z_1))^{-1} \rho_t^N(z_1)$ and $ (n_t^N(z_2))^{-1} \rho_t^N (z_2)$ are of the same type is again given by the integral with respect to $ \mathds{1}_{\diagdown} (k_1, k_2) = \mathds{1}_{k_1 = k_2} $. Formally, this would correspond to the integral
\begin{align*}
P_t^N (\psi_1 , \psi_2) 
&= \mathbb{E} \Bigg[ \int_{(\mathbb{T}^d)^2} \mathds{1}_{\diagdown} (k_1, k_2)  \psi_1 (z_1) \psi_2 (z_2) \frac{\rho_t^N (z_1, dk_1)}{n_t^N(z_1)} \frac{\rho_t^N (z_2, dk_2)}{n_t^N(z_2)} dz_1 dz_2 \Bigg] \\
&= \mathbb{E} \left[ \Bigg\langle \frac{\rho_t^N}{n_t^N} \otimes \frac{\rho_t^N}{n_t^N} , (\psi_1 \otimes \psi_2) \mathds{1}_{\diagdown} \Bigg\rangle \right].
\end{align*}
Our choice of sampling circumvents some problems with this approach arising from the regularity of $\rho_t^N/n_t^N$, whose fractional character was leading to difficulties. 
\end{remark}

\begin{remark} \label{remark:technicality}
Here, we wish to elaborate on the technicalities which prevented us from proving \Cref{conj:swapexpectationandlimit}. 
We will represent the fluctuations $\langle Z^N , \phi \rangle$ as a stochastic integral plus error terms, which result from the varying population size and vanish in the limit (see \eqref{eq:fluctuationsrepresentation}). The martingale problem associated to $Z^N \otimes Z^N$ would contain very similar error terms. Basically each error term would involve an additional copy of $Z^N$.
The additional $Z^N$ becomes problematic as it is not spatially averaged, preventing an application of \Cref{lem:boundonZsquaredwithtypes} or an appropriate generalization. One might consider the test function $\Phi' = \mathds{1}_{\diagdown} \overline{\phi_1} \otimes \overline{\phi_2}$ to gain the additional regularity in form of a spatial average. In the probability of identity, this would correspond to sampling a position and then adding a uniformly distributed vector from the unit ball. Unfortunately, most of the representations of the fluctuations used in this work involve choosing solutions to differential equations as test function to eliminate parts of the generator. If the growth term is among the eliminated ones, so is featured in the differential equation, the solution of the differential equation with averaged initial condition is not the same as the average of the original solution. This prevents the average from being transferred from the test function onto $Z^N$.
\end{remark}

\subsection{Interpretation and simulation}
Next, we wish elaborate on how \Cref{theo:varying_size_wmf} reflects the behaviour of the ancestral lineages in the process and provide visualizations of the corresponding Wright-Mal\'ecot formulae.
\subsubsection{Behaviour of two ancestral lineages}
Remember that in \Cref{theo:varying_size_wmf} we describe a probability of identity: in order for two individuals sampled from two different locations at time $t$ to have the same type, we need to track both their ancestral lineages to a common ancestor at some time $s < t$ and require that neither lineage acquires a mutation up to this point. The formula \eqref{eq:wmfformula} can then be interpreted as follows:
\begin{enumerate}[leftmargin=4\parindent]
\item[Mutation]
From an individual-based perspective each individual mutates at rate $\mu$. Hence, backwards in time starting from time $t$ the probability that neither ancestral lineage has seen a mutation at time $s$ in the past is given by $e^{- 2 \mu (t-s)}$.
\item[Motion]
We can write
\[ G_{s,t} \vartheta (z) = \int_{\mathbb{T}^d} p_{s,t} (z,x) \vartheta (x) dx,\]
where $z \mapsto p_{t-u, t} (z,x)$ is the transition density of $(X_u)_{u \in [0,t]}$, started from $X_0 = x$ and with the generator $u \mapsto A_{t-u}$ comprising the Laplacian and the drift according to twice the gradient of the population profile from \eqref{eq:generatorwmflineage}. The quantity $G_{s,t} \vartheta_1 (z) G_{s,t} \vartheta_2 (z)$ can be seen as the `probability' for both lineages to be at the same position $z$ a time $t-s$ in the past. This illustrates that the lineages behave according to Brownian motions with a drift according to the gradient given by the population profile.
\item[Coalescence]
The term $(n_s (z) +1)^{-1}$ captures the fact that once both lineages are at the same position $z$, coalescence occurs at a rate inversely proportional to the population size $n_s(z)$. This fits the intuition that with a larger local population size, the time to the most recent common ancestor while staying in the same deme should be longer.
\end{enumerate}
Notably, the motion of ancestral lineages is very much comparable to \cite[Theorem 2.24, Corollary 2.28 and Examples 3.2]{EKLRT24}. For example, in the case of profile $w$ of the population behaving according to a travelling wave solution to the Fisher-KPP equation
the authors obtain for the generator of ancestral lineages
\[
\mathcal{L} f = \partial_{xx} f + 2 \frac{\partial_x w}{w} \partial_x f + 2 \partial_x f.
\]
The $2 \partial_x f$ originates from the moving frame of the travelling wave with the remainder precisely equalling \eqref{eq:generatorwmflineage}. Their work is based on the lookdown framework and we arrive at a similar expression from a central limit theorem perspective. We can even go further and describe the behaviour of two ancestral lineages with coalescence happening inversely proportional to the population size. A hint that this is heuristically similar in their framework can be found in the discussion of Terence Tsui Ho Lung's Phd thesis \cite[Section 8.1.1]{phdthesis}. 

\subsubsection{Simulation and analytical prediction}
\begin{figure}[!htb]
    \centering
    \includegraphics[width=\linewidth, clip=true, trim = 10mm 10mm 10mm 10mm]{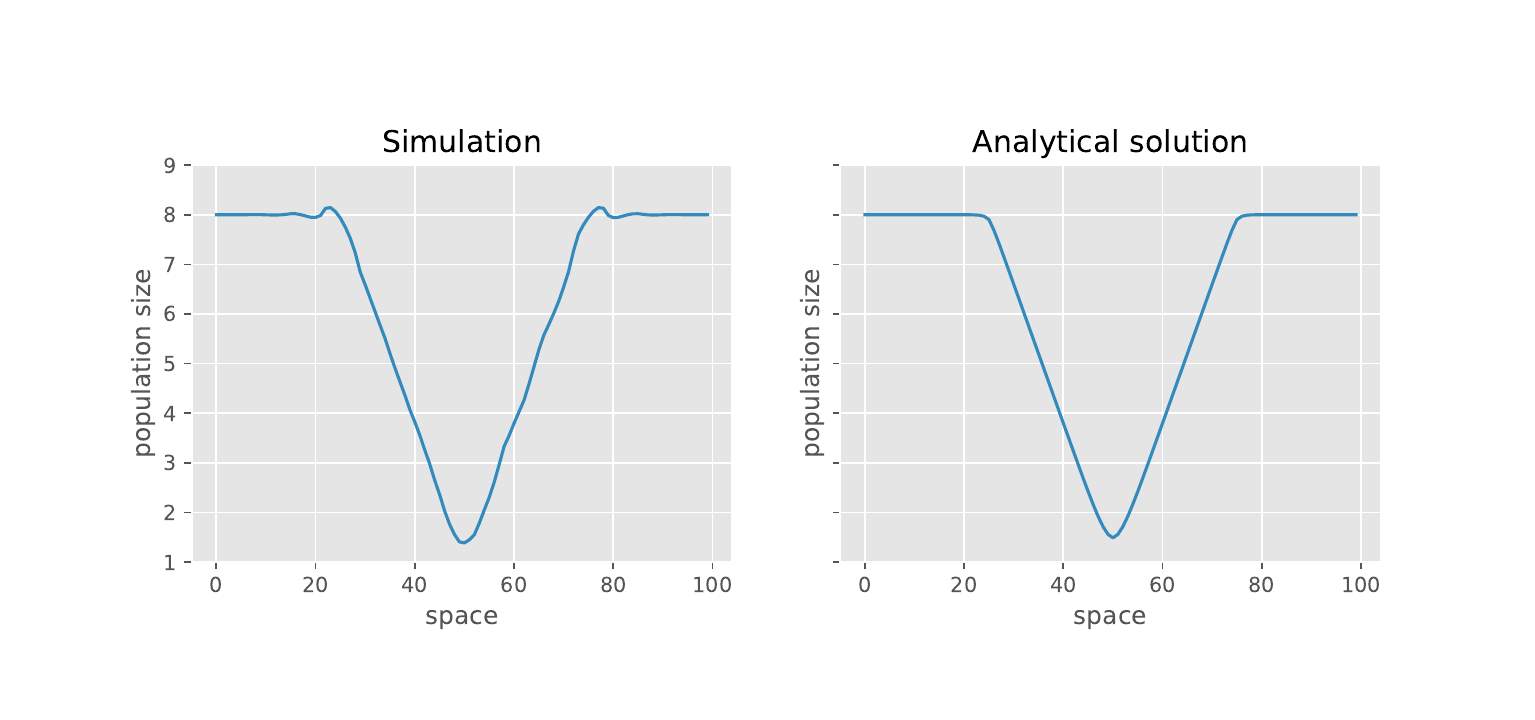}
    \caption{Comparison between the population size of the process run for a time of $125$ on the left and the analytical solution on the right, both corresponding to \eqref{eq:growthfunctionexample}. As the radius of events is positive, in events which fall right next to the edge of the valley, the growth function is positive leading to folds on the rim of the population density valley. The analytical solution depicts the solution to a discretized version of \eqref{eq:reactiondiffusion}.}
    \label{fig:pop_size_comparison}
\end{figure}

\begin{figure}[!htb]
    \centering
    \includegraphics[width=0.9\linewidth, clip=true, trim = 0mm 0mm 0mm 0mm]{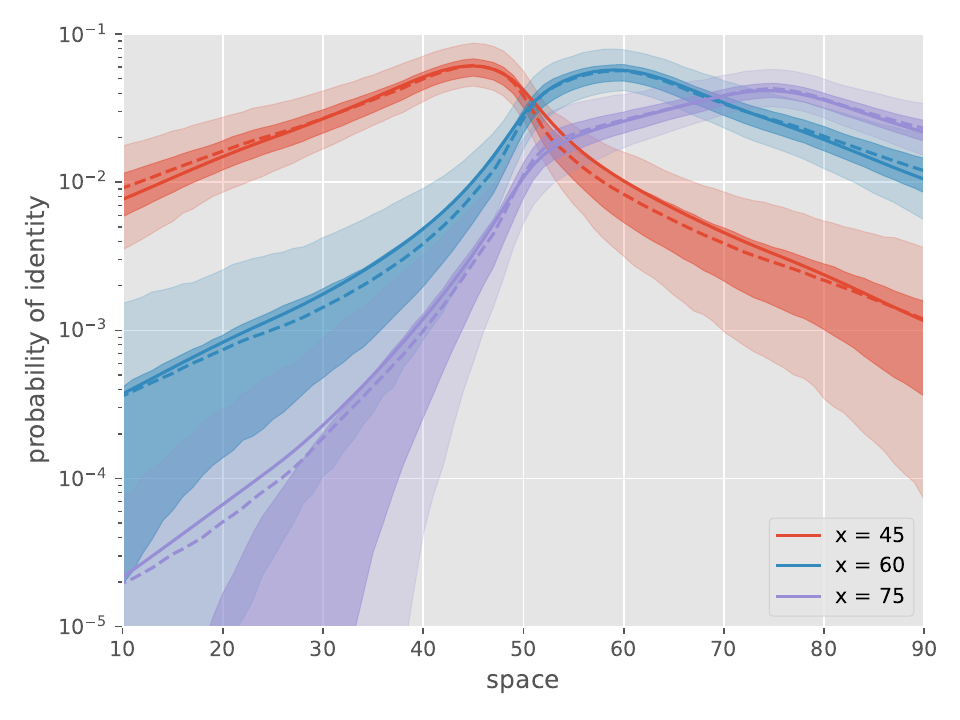}
    \caption{Comparison of the probabilities of identity (on a logarithmic scale) from three reference positions $45, 60, 75$ relative to all other locations. For example, red captures the probability of identity sampling one individual from $45$ and one from the population at the corresponding position on the x-axis. The continuous line shows the prediction $\Theta$ from the analytical calculation based on a calculation of the transition probabilities of ancestral lineages. The dotted lines illustrate the average probability of identity from the simulations of the process, with the lighter regions representing the $50$-th and $90$-th percentiles after $2000$ simulations. }
    \label{fig:1driftcomparison}
\end{figure}

\paragraph{Aim} In this section, we visualize the accuracy of the analytical result \Cref{theo:varying_size_wmf} by comparing it to simulated probabilities of identity in \Cref{fig:1driftcomparison}. We consider a scenario in which the population density is much lower at a central location (see \Cref{fig:pop_size_comparison}). In nature, this might be caused by, for example, a mountain range and will present in our model a barrier to gene flow. Ancestral lineages are attracted by higher population density, as we sample the parental type uniformly among the individuals in the affected region, so more likely from a densily populated area.

We will need to discretize space in both cases and will consider the interval $\{ 0,1,2, \dots, 99 , 100 \}$. The analytical result \Cref{theo:varying_size_wmf} was telling us the motion and coalescence behaviour of ancestral lineages. Discretizing this behaviour, we can calculate the analytical prediction $\Theta$ of the Wright-Mal\'ecot-Formula given by \eqref{eq:wmfformula}. We will compare this prediction with the probability of identity $\mathcal{P}$ resulting from forwards-in-time simulations of the typed process from \Cref{def:mSLFV}. To compare the analytical (limiting) result $\Theta$ with the simulated (prelimiting) probabilities, we `undo' the scaling. Our analyis predicts
\begin{equation*}
N \delta_N \mathcal{P}^N (x) \simeq \Theta \Big( \frac{N}{\delta_N^2} \mu_N , N u_N, \delta_N x \Big).
\end{equation*}
The full details of prediction and simulation will be given in \Cref{appendix:simulation}. 

\paragraph{Setting}
As illustrated by \Cref{fig:pop_size_comparison}, the population size profile will be at stationarity for the following choice of growth function
\begin{equation} \label{eq:growthfunctionexample}
r_x (n) = \max \Bigg\{ \min \Big\{ 14 \frac{\vert x - 50\vert}{50}, 7 \Big\} + 1 - n , - 1 \Bigg\} .
\end{equation}

To simulate the probability of identity, we consider a discretization of the dynamics given in \Cref{def:mSLFV} and follow the proportions of up to $2000$ types. For more details, please see \Cref{appendix:simulation}.

The analytical prediction $\Theta$ of \Cref{fig:pop_size_comparison} is based on calculating the transition probabilities of the lineages of sampled individuals backwards-in-time and incorporating the coalescence at a rate inversely proportional to the population size plus one. The lineages are diffusive with a drift according to the population size profile given in \Cref{fig:pop_size_comparison}.

In \Cref{fig:1driftcomparison} we consider three reference positions $45, 60$ and $75$. For each location $x_0$ on the x-axis, representing space, we visualize probability of identity of a sample from $x_0$ and a sample from the reference position. The prediction from the analytical solution is shown by the continuous line. The dashed line represents the mean of the probability of identity in $2000$ simulations of the process from \Cref{def:mSLFV}, whereas the darker and lighter area of the same colour show the $50$-th and $90$-th percentiles.

\paragraph{Interpretation}
The analytical calculation predicts the simulated probabilities of identities well.

The valley in the population density leads to ancestral lineages on the left of the central position $50$ being attracted to the left and lineages on the right being attracted to the right. Consequently, we would expect a rapid decorrelation of the probabilities of identity by descent around position $50$, which is confirmed by both a simulation of the process and the analytical prediction (see \Cref{fig:1driftcomparison}).

The mean of the simulated probabilities of identity, represented by the dashed line, is in the centre of the $50$-th and $90$-th percentiles close to reference positions $45, 60$ and $75$. This illustrates that in fact even without the expectation,
\[
\frac{\langle Z_t \otimes Z_t , \mathds{1}_{\diagdown} \psi_1 \otimes \psi_2 \rangle}{\langle n_t, \psi_1 \rangle \langle n_t, \psi_2 \rangle} \simeq \mathbb{E} \Bigg[ \frac{\langle Z_t \otimes Z_t , \mathds{1}_{\diagdown} \psi_1 \otimes \psi_2 \rangle}{\langle n_t, \psi_1 \rangle \langle n_t, \psi_2 \rangle} \Bigg]
\]
for sampling positions described by $\psi_1$ and $\psi_2$ being not too far apart. For locations further away from the reference positions a divergence of the mean from the median of the probabilities of identity occurs (particularly visible in the purple part on the left-hand side).

\section{Convergence of the population size process} \label{sec:convergenceofpopulationsize}
In this section, we will explain how the law of large numbers/\Cref{prop:llnpop} for the population size can be proven by adapting the strategy of \parencite[Section 3.2-4.3]{forien_central_2017}.
We consider test functions $\phi : \mathbb{T}^d \rightarrow \mathbb{R}$ depending only on space, so that the process evaluated against them describes the evolution of the population size. In other words, taking $\phi (x,k) = \phi (x)$, we define $\blangle n_t^N , \phi \brangle := \blangle \rho_t^N, \phi \brangle$ under the scaling \eqref{eq:scaling1} and \eqref{eq:scaling2}. We start with the characterisation via a martingale problem.

\begin{proof}[Proof of \Cref{lem:mppop}] 
If we apply the generator \eqref{eq:generator} to test functions only acting on space, we get
\begin{equation*}
\mathcal{G} F_{1, \phi} (n) = \int_{\mathbb{T}^d} (\overline{n} +1) \left( \Big\langle n + \frac{u}{\overline{n} +1} \mathds{1}_{B(x, R)} \Big([1 + r_x(\overline{n}) ]\overline{n} -n \Big) , \phi \Big\rangle - \blangle n, \phi \brangle \right)  dx.
\end{equation*}
The scaling \eqref{eq:scaling1} and \eqref{eq:scaling2} amounts to considering
\begin{equation*}
\blangle \rho_t^N , \phi \brangle = \blangle \mathfrak{p}_{Nt / \delta_N^2}^N (\delta_N \cdot) , \phi^N \brangle
\end{equation*}
for $\phi^N = \delta_N^d \phi (\delta_N \cdot)$. Recall that $\mathfrak{p}^N$ denotes the mSLFV with $u^N = u / N$ and $r^N_\cdot = \delta_N^2 r_{\delta_N \cdot} $ on $\delta_N^{-1} \mathbb{T}^d$. In the following calculation, we will apply the formula for multidimensional substitution 
\[ \int_{\mathbb{T}^d} \phi (x) dx = \int_{\delta_N^{-1} \mathbb{T}^d} \phi (\delta_N x ) \delta_N^d dx. \]
Indeed, we arrive at the desired expression:
\begingroup
\allowdisplaybreaks
\begin{align}
&\mathcal{G}^N F_{1, \phi} (n(\delta_N \cdot)) \nonumber \\
&= \frac{N}{\delta_N^2} \mathcal{G} F_{1, \phi^N} (n(\delta_N \cdot))\nonumber\\
&= \frac{N}{\delta_N^2}  \int_{\delta_N^{-1} \mathbb{T}^d} (\overline{n (\delta_N \cdot)}(x,R) +1) \nonumber \\
& \hspace{1cm} \Bigg[ \Big\langle n (\delta_N \cdot) + \frac{u^N}{\overline{n(\delta_N \cdot) }(x, R) +1} \mathds{1}_{B(x, R)} \Big([1 + r^N_x(\overline{n(\delta_N \cdot)}(x, R)) ]\overline{n(\delta_N \cdot)} (x, R) -n(\delta_N \cdot) \Big) , \phi^N \Big\rangle \nonumber \\
&\hspace{11.5cm} - \blangle n(\delta_N \cdot), \phi^N \brangle \Bigg]  dx\nonumber \\
&= \frac{u}{\delta_N^{d + 2}} \int_{\mathbb{T}^d}
 \Big\langle \mathds{1}_{B(x /\delta_N , R)} \Big([1 + r^N_{x /\delta_N} (\overline{n(\delta_N \cdot)}(x /\delta_N , R)) ]\overline{n(\delta_N \cdot)} (x /\delta_N , R) -n(\delta_N \cdot) \Big) , \phi^N \Big\rangle   dx\nonumber \\
&= \frac{u}{\delta_N^{d + 2}} \int_{\mathbb{T}^d}
 \int_{\delta_N^{-1} \mathbb{T}^d} \mathds{1}_{B(x , \delta_N R)} (\delta_N z) \Big([1 + \delta_N^2 r_{x} (\overline{n(\delta_N \cdot)}(x /\delta_N , R)) ]\overline{n(\delta_N \cdot)} (x /\delta_N , R) -n(\delta_N z) \Big)\nonumber \\
 &\hspace{12.5cm} \delta_N^d \phi( \delta_N z) dz dx\nonumber \\
&= \frac{u}{\delta_N^{d + 2}} \int_{\mathbb{T}^d}
 \int_{\mathbb{T}^d} \mathds{1}_{B(x , \delta_N R)} ( z) \Big([1 + \delta_N^2 r_{x} (\overline{n(\delta_N \cdot)}(x /\delta_N , R)) ]\overline{n(\delta_N \cdot)} (x /\delta_N , R) -n( z) \Big) \phi( z) dz dx\nonumber \\
 &= \frac{u}{\delta_N^{d + 2}} \int_{\mathbb{T}^d}
 \int_{\mathbb{T}^d} \mathds{1}_{B(x , \delta_N R)} ( z) \Big([1 + \delta_N^2 r_{x} (\overline{n}(x ,\delta_N R)) ]\overline{n} (x ,\delta_N R) -n( z) \Big) \phi( z) dz dx\nonumber \\
 &= \frac{u V_{\delta_N R} }{\delta_N^{d + 2}} \int_{\mathbb{T}^d}
\Bigg( \frac{1}{V_{\delta_N R}} \int_{B(z,\delta_N R)} [1 + \delta_N^2 r_{x} (\overline{n}(x ,\delta_N R)) ]\overline{n} (x ,\delta_N R) dx -n( z) \Bigg) \phi( z) dz \nonumber \\ 
&= u V_R \delta_N^{-2} \blangle  \overline{[1 + \delta_N^2 r(\overline{n}) ]\overline{n}} -n ,  \phi  \brangle\nonumber \\
&= u V_R \blangle ( \mathcal{L}^N + \mathcal{R}_n^N ) n , \phi \brangle. \label{eq:self_adjoint}
\end{align}
Note that for $g_1, g_2 : \mathbb{T}^d \rightarrow \mathbb{R}$ (and even for averages of any radius)
\begin{align*}
\langle \overline{g_1} , g_2 \rangle 
&= \int_{\mathbb{T}^d} \frac{1}{V_{\delta_N R}} \int_{\mathbb{T}^d} \mathds{1}_{\{ \vert z -x \vert \leq \delta_N R \}} g_1 (x) dx g_2 (z) dz \\
&= \int_{\mathbb{T}^d} \frac{1}{V_{\delta_N R }} \int_{\mathbb{T}^d} \mathds{1}_{\{ \vert z -x \vert \leq \delta_N R \}} g_2 (z) dz  g_1 (x) dx = \langle g_1 , \overline{g_2} \rangle. 
\end{align*}
This implies
\begin{align*}
\blangle \overline{[1 + \delta_N^2 r(\overline{n}) ]\overline{n}} ,  \phi  \brangle - \blangle n , \phi \brangle
&= \blangle [1 + \delta_N^2 r(\overline{n}) ]\overline{n} ,  \overline{\phi}  \brangle - \blangle n , \phi \brangle\\
&= \blangle \overline{n} , [1 + \delta_N^2 r(\overline{n}) ] \overline{\phi}  \brangle - \blangle n , \phi \brangle \\
&= \blangle n , \overline{[1 + \delta_N^2 r(\overline{n}) ] \overline{\phi}}  \brangle - \blangle n , \phi \brangle,
\end{align*}
which shows that $\mathcal{L}^N$ and $\mathcal{R}_n^N$ are self-adjoint and could also act on $\phi$ in \eqref{eq:self_adjoint}.

Similarly, for the quadratic variation
\begin{align*}
& \frac{N}{\delta_N^{d+2}} \int_{\mathbb{T}^d} (\overline{n} (x, \delta_N R) +1) \\
&\hspace{1cm} \left( \Big\langle n + \frac{u}{(\overline{n} (x, \delta_N R) +1)N} \mathds{1}_{B(x,\delta_N R)} \Big([1 + \delta_N^2 r_x(\overline{n} (x, \delta_N) ]\overline{n} (x,\delta_N) -n \Big) , \phi \Big\rangle - \blangle n, \phi \brangle \right)^2  dx\\
&= \frac{u^2}{\delta_N^{d+2}} \int_{\mathbb{T}^d} \frac{1}{(\overline{n} (x, \delta_N R) +1)N} \\
&\hspace{2cm} \Bigg( \int_{\mathbb{T}^d} \mathds{1}_{B(x,\delta_N R)} (z) \Big([1 + \delta_N^2 r_x(\overline{n} (x, \delta_N ) ]\overline{n} (x, \delta_N R) -n (z) \Big) \phi (z) dz \Bigg)^2  dx. \qedhere
\end{align*}
\endgroup
\end{proof}
The main idea of the proof of \Cref{prop:llnpop} is to separate the convergence of the population size into steps: 
\begin{enumerate}
\item 
The deterministic approximation $m^N$ of \eqref{eq:middleterm} given by
\begin{equation*}
\blangle m_t^N, \phi \brangle := \blangle m_0^N , \phi \brangle + \int_0^t \blangle m_s^N , (\mathcal{L}^N + \mathcal{R}_{m_s^N}^N ) \phi \brangle ds,
\end{equation*}
can be seen to converge to the solution of the reaction-diffusion equation.
\item 
Then, we will prove $\vert n^N - m^N \vert \to 0$ by defining fluctuations $F_t^N = \sqrt{N \delta_N^{2 - d}} (n_t^N - m_t^N)$ and showing that they are sufficiently bounded. Note that 
\begin{equation*}
F_t^N = \sqrt{N \delta_N^{2 - d}} (n_t^N - m_t^N) = \int_0^1 \sqrt{N \delta_N^{2-d}} (\rho_t^N - m_t^N \lambda) (dk) = \mathcal{H} Z_t^N,
\end{equation*}
which will become crucial in later sections.
\end{enumerate} 
The first step is established by the next lemma. As its proof is similar to several results in \Cref{appendix:propertiesofthefunctions}, which we will require for the proof of the main \Cref{theo:clt}, we only provide a sketch.
\begin{lemma}[cf. {\parencite[Proposition 4.6]{forien_central_2017}}] \label{lem:convergencemn}
There exist constants $C, C'> 0$, such that for any multi-index $ 0\leq \vert \kappa \vert \leq 4$
\begin{equation*}
\sup_{0 \leq t \leq T} \Vert \partial_\kappa m_t^N \Vert_\infty \leq C,
\end{equation*}
and uniform convergence holds
\begin{equation*}
\sup_{0 \leq t \leq T} \Vert m_t^N - n_t \Vert_\infty \leq \delta_N C'.
\end{equation*}
\end{lemma}
\begin{proof}
Let $G_t$ be the Gaussian kernel with variance $\frac{u V_R R^2}{d + 2} t$. Then the population size can be represented as
\begin{equation*}
    n_t (x) = G_t * n_0 (x) + \int_0^t (G_{t-s} * r(n_s) n_s) (x) ds.
\end{equation*}
Similarly, if $G_t^N$ denotes the semigroup associated to the generator $\mathcal{L}^N$, we get
\begin{equation*}
    m_t^N = G_t^N * m_0^N + \int_0^t( G_{t-s}^N * \overline{r(\overline{m_s^N}) \overline{m_s^N}} ) (x) ds.
\end{equation*}
The difference between the two terms can now be bounded using Gr\"onwall's inequality. Using \Cref{assumption} the argument of \parencite[Section B.1, Proof of Proposition 4.6]{forien_central_2017} carries over.
\end{proof}
To perform the second step first note that if we assume the same initial conditions $n_0^N = m_0^N a.s.$ the fluctuations satisfy the following semimartingale decomposition
\begin{align*}
\blangle F_t^N, \phi \brangle
&= u V_R \sqrt{N \delta_N^{2 - d}} \int_0^t \blangle (\mathcal{L}^N + \mathcal{R}_{n_s^N}) (n_s^N) - (\mathcal{L}^N + \mathcal{R}_{m_s^N} ) (m_s^N), \phi \brangle ds + \sqrt{N \delta_N^{2 - d}} Y_t^N (\phi)\\
&=  u V_R \int_0^t \blangle \mathcal{L}^N F_s^N + \sqrt{N \delta_N^{2 - d}} \overline{ \big( r(\overline{n^N_s}) \overline{n^N_s} - r(\overline{m^N_s}) \overline{m^N_s} \big)} (\delta_N R), \phi \brangle ds  + \sqrt{N \delta_N^{2 - d}} Y_t^N (\phi).
\end{align*}
The martingale $Y_t^N(\phi)$ had been defined in \Cref{lem:mppop}. This equation is of the same form as \parencite[Equation (45)]{forien_central_2017}. The next lemma controls the martingale part.
\begin{lemma}[cf. {\parencite[Lemma 4.4]{forien_central_2017}}] \label{lem:martingalepopgood}
There exists $C > 0$ such that for all $\phi : \mathbb{R}_+ \rightarrow L^2 (\mathbb{T}^d)$ we have
\begin{equation}
\mathbb{E} \Bigg[ \Big( \sqrt{N \delta_N^{2-d}} \int_{[0,t] \times \mathbb{T}^d} \phi_s (x) Y^N (dx\; ds) \Big)^2  \Bigg] \leq C \int_0^t \Vert \phi_s \Vert_2^2 \; ds.
\end{equation}
\end{lemma}
\begin{proof}
From \Cref{lem:mppop} we know that the left-hand side is equal to
\begin{align*}
&\mathbb{E} \Bigg[ \int_0^t \frac{u^2}{\delta_N^{2d}} \int_{\mathbb{T}^d} \frac{1}{(\overline{n} (x, \delta_N R) +1)} \\
&\hspace{1cm} \Bigg( \int_{\mathbb{T}^d} \mathds{1}_{B(x,\delta_N R)} (z) \Big([1 + \delta_N^2 r_x(\overline{n^N} (x,\delta_N R)) ]\overline{n^N} (x, \delta_NR) -n^N (z) \Big) \phi_s (z) dz \Bigg)^2  dx ds  \Bigg]\\
&\leq \mathbb{E} \Bigg[ \int_0^t \frac{u^2}{\delta_N^{2d}} \int_{\mathbb{T}^d} \Bigg( \int_{\mathbb{T}^d} \mathds{1}_{B(x,\delta_N R)} (z) \; 2 n_{\max} \; \vert \phi_s (z) \vert dz \Bigg)^2  dx ds  \Bigg] \\
&\leq 4 n_{\max}^2 u^2 V_R^2 \mathbb{E} \Bigg[ \int_0^t \int_{\mathbb{T}^d} \Bigg( \frac{1}{V_{\delta_N R}} \int_{B(x, \delta_N R)} \vert \phi_s (z) \vert dz \Bigg)^2  dx ds  \Bigg] \\
&\leq 4 n_{\max}^2 u^2 V_R^2   \int_0^t \Vert \phi \Vert_2^2 \; ds,
\end{align*}
where we have used $ \Vert \overline{\phi} \Vert_2 \leq \Vert \phi \Vert_2$ in the last line.
\end{proof}
The crucial idea to control the fluctuations $F^N$ is that the square of the spatial average is uniformly controlled. To prove the convergence of $Z^N$, we will use \Cref{lem:boundonZsquaredwithtypes}, which can be seen as a generalization of the next result to involve types. 
\begin{lemma}[cf. {\parencite[Lemma 4.5]{forien_central_2017}}] \label{lem:boundonFsquared}
There exists a constant $C> 0$ such that
\begin{equation*}
\sup_{0 \leq t \leq T} \sup_{z \in \mathbb{T}^d} \mathbb{E} \left[ \overline{F_t^N} (z, \delta_N R)^2 \right] \leq \frac{C}{V_{\delta_N R}}.
\end{equation*}
\end{lemma}
\begin{proof}
The remainder of Taylor's expansion to the $k$-th order of a function $\mathfrak{r} \in C^{k}(\mathbb{R})$ about the point $x$ is given by
\begin{equation} \label{eq:taylor_remainder}
\mathfrak{T}_{k} (x,y, \mathfrak{r}) = \int_0^1 \frac{(1- t)^{k-1}}{(k-1)!} \mathfrak{r}^{(k)} (x + t (y - x)) dt.
\end{equation}
Applying Taylor's expansion to $\mathfrak{r} = r(x) x$ yields
\begin{equation} \label{eq:representation_fluc_pop}
\blangle F_t^N, \phi \brangle
=  u V_R \int_0^t \blangle \mathcal{L}^N F_s^N + \overline{ F_s^N \mathfrak{T}_1(\overline{n_s^N}, \overline{m_s^N}, \mathfrak{r} )}, \phi \brangle ds  + \sqrt{N \delta_N^{2 - d}} Y_t^N (\phi).
\end{equation}
We wish to remove $\mathcal{L}^N$ by using an appropriate test function - a technique which we will use frequently. Let $G_t^N$ be again the semigroup associated to the operator $\mathcal{L}^N$. Substituting the time-dependent test function $G_{t-s}^N * \phi$ in \eqref{eq:representation_fluc_pop} yields
\begin{equation*}
\blangle F_t^N, \phi \brangle =  u V_R \int_0^t \Blangle G_{t-s}^N * \Big( \overline{ F_s^N \mathfrak{T}_1(\overline{n_s^N}, \overline{m_s^N}, \mathfrak{r} )} \Big), \phi \Brangle ds + \int_0^t \int_{\mathbb{T}^d} (G_{t-s}^N * \phi)(z) \sqrt{N \delta_N^{2 - d}} Y^N (dz\; ds).
\end{equation*}
Choosing $\phi = \frac{1}{V_R} \mathds{1}_{\vert \cdot - z \vert < \delta _N R}$, the left-hand side $\langle F_t^N , \phi \rangle$ equals $\overline{F_t^N} (z, \delta_N R)$. Squaring the equation, and using Jensen's inequality as well as \Cref{lem:martingalepopgood} to bound the martingale part, one can set up an application of Gr\"onwall's inequality. We provide details of these steps in the more general setting of \Cref{lem:boundonZsquaredwithtypes} (in which we include genetic types) and so we omit the details here.
\end{proof}
\begin{lemma}[cf. {\parencite[Lemma 4.7]{forien_central_2017}}] \label{lem:boundflucpop}
There exists a constant $C > 0$ such that for any twice continuously differentiable function $\phi$ with
\begin{equation}
\max_{p \in \{1,2\} } \max \Big\{ \Vert \phi \Vert_p , \max_{0 \leq \vert \kappa \vert \leq 2} \Vert \partial_\kappa \phi \Vert_p \Big\} \leq 1,
\end{equation}
we have
\begin{equation} \label{eq:lemboundflucdesiredbound}
\mathbb{E} \Bigg[ \sup_{0 \leq t \leq T} \Big\vert \blangle F_t^N , \phi \brangle \Big\vert \Bigg] \leq C.
\end{equation}
\end{lemma}
\begin{proof}
Similar to \parencite[Section 3.2]{forien_central_2017}, using Taylor's expansion of $\mathfrak{r} = r(x) x$, this time to the second order
\begin{equation} \label{eq:firstrepresentationF}
\begin{aligned}
\blangle F_t^N, \phi \brangle
&=  u V_R \int_0^t \blangle F_s^N, \mathcal{L}^N \phi + \overline{ \mathfrak{r}' (\overline{m_s^N} ) \overline{\phi} } (\delta_N R) \brangle ds \\
& \hspace{1cm} + \frac{1}{\sqrt{N \delta_N^{2 - d}}} \int_0^t \blangle (\overline{F_s^N})^2,  \mathfrak{T}_2 (\overline{n_s^N}, \overline{m_s^N}, \mathfrak{r}) \overline{\phi} \brangle ds + \sqrt{N \delta_N^{2 - d}} Y_t^N (\phi)
\end{aligned}
\end{equation}
where the last line is of the same form as \parencite[Equation (48)]{forien_central_2017}. To obtain \eqref{eq:lemboundflucdesiredbound}, we would like to take absolute values, supremum over time and expectation on both sides of \eqref{eq:firstrepresentationF}. Inescapably, this will lead to $\mathbb{E} \big[ \big\vert \blangle F_s^N, \mathcal{L}^N \phi + \overline{ \mathfrak{r}' (\overline{m_s^N} ) \overline{\phi} } (\delta_N R) \brangle \big\vert \big] $ on the right-hand side. We therefore require first a bound on $\mathbb{E} [ \vert \langle F_s^N , \phi \rangle \vert ] $. Indeed, choosing $\phi^N$ to be the time-dependent solution to 
\begin{equation*}
    \begin{dcases}
    \partial_s \phi^N (x,s,t) + \mathcal{L}^N \phi^N (x,s,t) - \overline{\mathfrak{r}'(\overline{m_s^N}) \overline{\phi^N (s,t)}} (x, \delta_N R) = 0 \\
    \phi^N (x,t,t) = \phi (x)
    \end{dcases},
\end{equation*}
one obtains instead the equation (see \parencite[Section 4.1]{forien_central_2017})
\begin{equation} \label{eq:flucpopsecondform}
\begin{aligned}
\blangle F_t^N, \phi \brangle &= \frac{1}{\sqrt{N \delta_N^{2 - d}}} \int_0^t \blangle (\overline{F_s^N})^2, \mathfrak{T}_2 (\overline{n_s^N}, \overline{m_s^N}, \mathfrak{r}) \overline{\phi^N(s,t)} \brangle ds\\
&\hspace{1cm} + \int_0^t \int_{\mathbb{T}^d} \phi^N(x,s,t) \sqrt{N \delta_N^{2 - d}} Y^N(dx\; ds),
\end{aligned}
\end{equation}
which is of the same form as \parencite[Equation (58)]{forien_central_2017}. 
The first term is of the shape necessary to apply \Cref{lem:boundonFsquared}, whereas the second can again be controlled by \Cref{lem:martingalepopgood} for $t \in [0,T]$, yielding
\begin{align*}
\mathbb{E} \big[ \big\vert \langle F_t^N, \phi \rangle \big\vert \big] &\leq \frac{1}{\sqrt{N \delta_N^{2-d}}} T \frac{C}{V_{\delta_N R}} C' \big\Vert \phi^N (s,t) \big\Vert_1 + \Bigg( C'' \int_0^t \big\Vert \phi^N (s,t) \big\Vert_2^2 ds \Bigg)^{1/2} .
\end{align*}
We choose the constant $C'$ to be a uniform bound on all derivatives of the growth function, whose existence was postulated in point 3. of \Cref{assumption}. The norms $\big\Vert \phi^N (s,t) \Vert_1$ and $\Vert \phi^N (s,t) \Vert_2$ can be bounded by $\Vert \phi \Vert_1$ and $\Vert \phi \Vert_2$ respectively. This can be proven similarly to \Cref{lem:semigroupsbounded}, and as this proof is more complex, we omit the details. By the assumption \eqref{eq:assumptionondelta} on $\delta_N$, $(N^{\frac{1}{2}} \delta_N^{1 + \frac{d}{2}})^{-1}$ is bounded for $N \geq 1$, and in total we arrive at
\[
\mathbb{E} \big[ \big\vert \langle F_t^N, \phi \rangle \big\vert \big] \leq \frac{\mathfrak{C}}{N^{\frac{1}{2}} \delta_N^{1 + \frac{d}{2}}} \Vert \phi \Vert_1 + \mathfrak{C}' \Vert \phi \Vert_2,
\]
for two constants $\mathfrak{C}, \mathfrak{C}' \geq 0$.

Returning to \eqref{eq:firstrepresentationF} above, one can then apply this bound to the function $\mathcal{L}^N \phi + \overline{ \mathfrak{r}' (\overline{m_s^N} ) \overline{\phi} } (\delta_N R) $:
\begin{equation*}
\begin{aligned}
\mathbb{E} \Bigg[ \sup_{0 \leq t \leq T} \Big\vert \blangle F_t^N , \phi \brangle \Big\vert \Bigg] &\leq T \frac{\mathfrak{C}}{N^{\frac{1}{2}} \delta_N^{1 + \frac{d}{2}}} \big( \Vert \mathcal{L}^N \phi \Vert_1 + C' \Vert \phi \Vert_1 \big) + T \mathfrak{C}' \big( \Vert \mathcal{L}^N \phi \Vert_2 + C' \Vert \phi \Vert_2 \big) \\
& \hspace{1cm} + \frac{\mathfrak{C}}{N^{\frac{1}{2}} \delta_N^{1 + \frac{d}{2}}} \Vert \phi \Vert_1 + \mathbb{E} \Bigg[ \sup_{t \in [0,T]} \sqrt{N \delta_N^{2-d}} \big\vert Y_t^N (\phi) \big\vert^2 \Bigg]^{1/2}.
\end{aligned}
\end{equation*}
The last summand is bounded by $4 \mathfrak{C}' \Vert \phi \Vert_2$ by Doob's inequality and another application of \Cref{lem:martingalepopgood}. By \cite[Proposition A.1]{forien_central_2017}, $\Vert \mathcal{L}^N \phi \Vert $ can be bounded in terms of the norms of the second derivatives of $\phi$ and the result follows. \qedhere
\end{proof}

We can now conclude the law of large numbers for the population size.
\begin{proof}[Proof of \Cref{prop:llnpop}]
We proceed with the same argument as \parencite[Section 4.3]{forien_central_2017} by combining \Cref{lem:convergencemn} and \Cref{lem:boundflucpop}. Without loss of generality, we assume the existence of a constant $C$ s.t. the sequence of test functions satisfies $ \max_{0 \leq \vert \kappa \vert \leq 2} \Vert \partial_\kappa \phi_n \Vert_p < C$ uniformly for all $n\geq 1$ and $p \in \{ 1, 2\}$ (see \eqref{eq:assumption_on_test_sequence}). We can now bound
\begin{align*}
&\mathbb{E} \Bigg[ \sup_{0 \leq t \leq T} d(n_t^N , n_t ) \Bigg]\\
&\leq \sum_{i = 1}^\infty \frac{1}{2^n} \Bigg( \mathbb{E} \Big[ \sup_{0 \leq t \leq T} \Big\vert \blangle n_t^N , \phi_i \brangle - \blangle m_t^N, \phi_i \brangle \Big\vert \Big] + \sup_{0 \leq t \leq T} \Big\vert \blangle m_t^N , \phi_i \brangle - \blangle n_t, \phi_i \rangle \Big\vert \Bigg) \\
&\leq \sum_{i = 1}^\infty \frac{1}{2^n} \Bigg( (N \delta_N^{2 - d})^{-1/2} \mathbb{E} \Big[ \sup_{0 \leq t \leq T} \Big\vert \blangle F_t^N , \phi_i \brangle \Big\vert \Big] + \sup_{0 \leq t \leq T} \big\Vert m_t^N - n_t \big\Vert_\infty \Vert \phi_n \Vert_1 \Bigg),
\end{align*}
where the right-hand side converges to zero as $N \to \infty$ due to \Cref{lem:convergencemn} and \Cref{lem:boundflucpop}.
\end{proof}
\section{Semimartingale decomposition of the fluctuations} \label{sec:martingaleproblems}
In this section, we will derive a representation of the fluctuations and prove some auxiliary results, in preparation for the convergence proofs in the next section. 
\subsection{Rescaled martingale problem for the process}
We start by rescaling the martingale problem for the mSLFV process.
\begin{lemma} \label{lem:original_mp}
Let $(\rho_t, t \geq 0)$ be an mSLFV process and $\mathcal{G}$ be the generator given by \eqref{eq:generator}. Then, for any test function $\phi \in L^1 (\mathbb{T}^d \times [0,1])$,
\begin{equation} \label{eq:mp_process_martingale}
\mathcal{M}_t (F_{f, \phi}) :=f( \langle \rho_t , \phi \rangle ) - f \langle \rho_0 , \phi \rangle ) - \int_0^t \mathcal{G} F_{f, \phi} (\rho_s ) ds
\end{equation}
is a square integrable martingale with predictable quadratic variation
\begin{equation*}
\begin{aligned}
& \int_{\mathbb{T}^d}  \frac{1}{\int_{B(x,R)\times [0,1]} \rho(y,dk) dy} \int_{B(x,R) \times [0,1]}   \\
& \hspace{1cm} (\overline{n} +1) \Big[ f(\langle \rho + \frac{u}{\overline{n}+1} \mathds{1}_{B(x,R)} (\delta_{k_0} [1 + r(\overline{n} ) ] \overline{n} - \rho) , \phi \rangle)  - f(\langle \rho, \phi \rangle) \Big]^2\\
&\hspace{10cm} \rho (y, dk_0) dy dx.
\end{aligned}
\end{equation*}
\end{lemma}
In very much the same way as for the population size, the dynamics can be expressed with the operators $\mathcal{L}^N$ and $\mathcal{R}_n^N$. Again, we will abbreviate notation by shortening $\overline{n} (x, \delta_N R)$ to $\overline{n}$.
\begin{lemma} \label{lem:mp}
Let $(\rho_t^N, t \geq 0)$ be an mSLFV process under the scaling \eqref{eq:scaling1} and \eqref{eq:scaling2}. Then, for $\phi \in L^1 (\mathbb{T}^d \times [0,1])$,
\begin{equation*}
M_t^N(\phi) := \blangle \rho_t^N, \phi \brangle - \blangle \rho_0^N, \phi \brangle - \int_0^t \Big[ \mu \blangle n_s^N \lambda - \rho_s^N , \phi \brangle +  \blangle \mathcal{L}^N \big( \rho_s^N \big) + \mathcal{R}_{n_s^N}^N \big( \rho_s^N \big) , \phi \brangle \Big] ds
\end{equation*}
is a martingale with predictable variation process
\begin{align*}
&\int_0^t \frac{N}{\delta_N^{2+d}} \int_{\mathbb{T}^d} \frac{1}{\int_{B(x,\delta_N R)} n_s^N(y) dy} \int_{B(x,\delta_N R) \times [0,1]}  (\overline{n_s^N} (x,\delta_N R) +1)\\
& \hspace{1cm} \Bigg[ \Big\langle \rho_s^N +  \frac{u}{(\overline{n_s^N} (x, \delta_NR) +1)} \mathds{1}_{B(x, \delta_N R)} (\cdot) ( \delta_{k_0}[1 + \delta_N^2 r_x(\overline{n_s^N} (x, \delta_N)) ] \overline{n_s^N} (x, \delta_N R) - \rho_s^N) , \phi \Big\rangle\\
&\hspace{9cm} - \Big\langle \rho_s^N, \phi \Big\rangle \Bigg]^2 \rho_s^N (y, dk_0) dy dx ds .
\end{align*}
\end{lemma}
\begin{proof}[Proof of Lemmata \ref{lem:original_mp} and \ref{lem:mp}]
We have
\begin{align*}
&\lim_{\delta t \downarrow 0} \frac{1}{\delta t} \mathbb{E} \Big[ \blangle \rho_{t + \delta t} , \phi \brangle  - \blangle \rho_t , \phi \brangle \Big\vert \rho_t = \rho \Big]\\
& = \int_{\mathbb{T}^d} \Bigg( \frac{1}{\int_{B(x,R)} n(y) dy} \int_{B(x,R) \times [0,1]} (\overline{n} +1) \Big\langle \rho +  \frac{u}{(\overline{n} +1)} \mathds{1}_{B(x,R)} ( \delta_{k_0}[1 + r(\overline{n}) ] \overline{n} - \rho) , \phi \Big\rangle \\
& \hspace{2cm} - (\overline{n} +1) \Big\langle \rho, \phi \Big\rangle \rho (y, dk_0) dy \Bigg) dx + \mu \Big\langle \lambda \cdot n - \rho , \phi \Big\rangle\\
& = u \int_{\mathbb{T}^d} \Bigg( \frac{1}{\int_{B(x,R)} n(y) dy} \int_{B(x,R) \times [0,1]} \Big\langle \mathds{1}_{B(x,R)} ( \delta_{k_0}[1 + r(\overline{n}) ] \overline{n} - \rho) , \phi \Big\rangle \rho (y, dk_0) dy \Bigg) dx \\
& \hspace{1cm} + \mu \Big\langle \lambda \cdot n - \rho , \phi \Big\rangle.
\end{align*}
This implies by e.g. \cite[Proposition 4.1.7]{ethier_markov_1986} that the expression \eqref{eq:mp_process_martingale} is a martingale. The quadratic variation can be obtained by considering a similar small-time limit of the expectation of $\big( \blangle \rho_{t + \delta t} , \phi \brangle  - \blangle \rho_t , \phi \brangle \big)^2$ conditioned on $\rho_t = \rho$.

Next, we again wish to apply \eqref{eq:scaling1} and \eqref{eq:scaling2} to derive the expressions of \Cref{lem:mp} from \Cref{lem:original_mp}. Defining the rescaled test functions $
\phi_N (x) = \delta_N^d \phi (\delta_N x)$,
the rescaled generator will be
\begin{equation*}
\mathcal{G}^N F_{f, \phi}(\rho) = \frac{N}{\delta_N^2} \mathcal{G} F_{f, \phi_N} (\rho (\delta_N \cdot)),
\end{equation*}
where the factor $N /\delta_N^2$ originates from the rescaling of time. Recall that we scale impact by $u^N = u/N$ and mutation rate by $\mu_N = \delta_N^2 \mu /N$. The rescaled growth function was defined as $r^N_\cdot = \delta_N^2 r_{\delta_N \cdot}$. For the mutation part, we obtain
\begin{align*}
&\frac{N}{\delta_N^2} \frac{\delta_N^2 \mu}{N} \int_{\delta_N^{-1} \mathbb{T}^d} \big( \lambda \cdot n (\delta_N z) - \rho (\delta_N z) \big) \phi (\delta_N z ) \delta_N^d dz\\
&= \mu \int_{\mathbb{T}^d} \big(\lambda \cdot n (z)- \rho(z) \big) \phi(z) dz.
\end{align*}
In a similar way to our calculations for the population size, we can derive, abbreviating $\overline{n}(x, \delta_NR)$ by $\overline{n}$,
\begingroup
\allowdisplaybreaks
\begin{align*}
&\frac{u}{N} \frac{N}{\delta_N^{d+2}} \int_{\mathbb{T}^d} \frac{1}{\int_{B(x,\delta_N R)} n ( y) dy } \int_{B(x,\delta_N R) \times [0,1]} \int_{B(x,\delta_N R)} \Bigg[ [1 + \delta_N^2 r_x (\overline{n})] \overline{n} \phi (z,k_0) \\
& \hspace{1cm} - \int_{[0,1]} \phi ( z,k) \rho (z, dk) \Bigg] dz \rho ( y, dk_0 ) dy dx\\
&= \frac{u}{\delta_N^{d+2}} \int_{\mathbb{T}^d} \frac{1}{\int_{B(x,\delta_N R)} n (y) dy } \int_{B(x, \delta_N R) \times [0,1]} \int_{B(x,\delta_N R)} \Big[ [1 + \delta_N^2 r_x (\overline{n})] \overline{n} \phi (z,k_0) \Big] dz \rho (y, dk_0 ) dy dx\\
& \hspace{1cm} - \int_{\mathbb{T}^d} \int_{B(x,\delta_N R)} \int_{[0,1]} \phi (z,k) \rho (z, dk) dz dx\\
&= \frac{u}{\delta_N^{d+2}} \int_{\mathbb{T}^d \times [0,1]}  \int_{B(y,\delta_N R)} \frac{1}{\int_{B(x,\delta_N R)} n (y) dy } \int_{B(x,\delta_N R)} [1 +\delta_N^2 r_x (\overline{n})] \overline{n} \phi (z,k_0) dz dx \rho (y, dk_0 ) dy \\
& \hspace{1cm} - V_{\delta_N R} \int_{\mathbb{T}^d \times [0,1]} \phi (z,k) \rho (z, dk) dz\\
&= \frac{u}{\delta_N^{d+2}} \int_{\mathbb{T}^d \times [0,1]} \Bigg[  \int_{B(z,\delta_N R)} \frac{1}{V_{\delta_N R}} \int_{B(x, \delta_N R)} [1 + \delta_N^2 r_x (\overline{n})] \phi (w,k)  dw dx - V_{\delta_N R} \phi (z,k) \Bigg] \rho (z, dk) dz\\
& = \Big\langle \rho (\cdot, dk), (\mathcal{L}^N + \mathcal{R}_{n}^N) (\phi (\cdot, k )) \Big\rangle.
\end{align*}
\endgroup
The predictable variation formula of \Cref{lem:mp} can be derived by a similar strategy.
\end{proof}
\subsection{Representations of the fluctuations} \label{subsec:representation_fluctuations}
Finally, let us return to the fluctuations defined by
\begin{equation*}
Z_t^N = \sqrt{N \delta_N^{2-d}} \Big(\rho_t^N - m_t^N \lambda \Big).
\end{equation*}
The process $\rho_t^N$ and the deterministic approximation to the population size $m_t^N \lambda$ have the following representations
\begin{align*}
\langle \rho_t^N , \phi \rangle &= \langle \rho_0^N, \phi \rangle + \int_0^t \Big[ \mu \langle n_s^N - \rho_s^N, \phi \rangle + \langle \rho_s^N , \mathcal{L}^N (\phi) + \mathcal{R}_{n_s^N} (\phi) \rangle \Big] ds + M_t^N (\phi), \\
\langle m_t^N \lambda , \phi \rangle &= \langle m_0^N \lambda, \phi \rangle + \int_0^t \langle m_s^N \lambda , \mathcal{L}^N (\phi) + \mathcal{R}_{m_s^N} \phi \rangle ds.
\end{align*}
The are now different options for how to summarize the resulting terms for the fluctuations. Let us define
\[\mathfrak{M}_t^N (\phi) := \sqrt{N \delta_N^{2-d}} M_t^N (\phi).\]
We will use
\begingroup
\allowdisplaybreaks
\begin{align*}
&\langle Z_t^N , \phi \rangle \\
&= \sqrt{N \delta_N^{2-d}} \Bigg[ \langle \rho_0^N, \phi \rangle + \int_0^t \Big[ \mu \langle n_s^N \lambda - \rho_s^N, \phi \rangle + \langle \rho_s^N , \mathcal{L}^N (\phi) + \mathcal{R}_{n_s^N}^N (\phi) \rangle \Big] ds + M_t^N (\phi)\\
& \hspace{1cm} - \langle m_0^N \lambda, \phi \rangle - \int_0^t \langle m_s^N \lambda, \mathcal{L}^N (\phi) + \mathcal{R}_{m_s^N}^N \phi \rangle ds \Bigg]\\
&= \langle Z_0^N , \phi \rangle + \mathfrak{M}_t^N (\phi) + \int_0^t \langle Z_s^N , \mathcal{L}^N (\phi ) \rangle ds\\
&\hspace{1cm} +  \int_0^t \Big[ - \mu \langle Z_s^N, \phi \rangle + \sqrt{N \delta_N^{2-d}} \mu \big\langle (n_s^N - m_s^N) \lambda , \phi \big\rangle \Big] ds \\
& \hspace{1cm} + \int_0^t \Big[ \langle Z_s^N ,\mathcal{R}_{n_s^N} (\phi) \rangle + \sqrt{N \delta_N^{2-d}} \Big\langle m_s^N \lambda , \big( \mathcal{R}_{n_s^N}^N - \mathcal{R}_{m_s^N}^N \big) \phi \Big\rangle \Big] ds\\
&= \langle Z_0^N , \phi \rangle + \mathfrak{M}_t^N (\phi) + \int_0^t \big\langle Z_s^N , \mathcal{L}^N (\phi ) + \mathcal{R}_{n_s^N}^N (\phi) - \mu \phi \big\rangle ds\\
&\hspace{1cm} + \int_0^t \mu \big\langle Z_t^N , \mathcal{H} \phi \big\rangle ds + \sqrt{N \delta_N^{2-d}} u V_R \int_0^t \Big\langle m_s^N \lambda , \overline{r (\overline{n_s^N}) \overline{\phi}} (\cdot, \delta_N R) - \overline{r (\overline{m_s^N}) \overline{\phi}} (\cdot, \delta_N R) \Big\rangle ds\\
&= \langle Z_0^N , \phi \rangle + \mathfrak{M}_t^N (\phi) + \int_0^t \big\langle Z_s^N , \mathcal{L}^N (\phi ) + \mathcal{R}_{n_s^N}^N (\phi) + \mu ( \mathcal{H} \phi - \phi) \big\rangle ds\\
&\hspace{1cm} + \sqrt{N \delta_N^{2-d}} u V_R \int_0^t \Big\langle m_s^N \lambda , \overline{ [r (\overline{n_s^N}) - r(\overline{m_s^N}) ] \overline{\phi}} (\cdot, \delta_N R) \Big\rangle ds,
\end{align*}
where again $\mathcal{H} \phi (x) = \int_0^1 \phi (x,k) dk$. The last summand will make a non-trivial contribution to the central limit theorem. To see this, we apply Taylor's expansion to obtain
\begin{align*}
&\sqrt{N \delta_N^{2-d}} u V_R \int_0^t \Big\langle \overline{m_s^N} \lambda , \big[r'(\overline{m_s^N}) (\overline{n_s^N} - \overline{m_s^N}) + \mathfrak{T}_2 (\overline{n_s^N}, \overline{m_s^N}, r) (\overline{n_s^N} - \overline{m_s^N})^2 \big] \overline{\phi} \Big\rangle ds\\
&= u V_R \int_0^t \Big\langle \overline{m_s^N} r'(\overline{m_s^N}) \overline{Z}_s^N, \overline{\mathcal{H} \phi} \Big\rangle ds + u V_R \int_0^t \Big\langle \overline{m_s^N} \mathfrak{T}_2 (\overline{n_s^N}, \overline{m_s^N}, r) (\overline{n_s^N} - \overline{m_s^N}) \overline{Z}_s^N, \overline{\mathcal{H} \phi} \Big\rangle ds.
\end{align*}
\endgroup
Returning to the complete expression, this yields
\begin{align*}
\langle Z_t^N , \phi \rangle &= \langle Z_0^N , \phi \rangle + \int_0^t \big\langle Z_s^N , \mathcal{L}^N (\phi ) + \mathcal{R}_{n_s^N}^N (\phi) + \mu ( \mathcal{H} \phi - \phi) + uV_R \overline{\overline{m_s^N} r'(\overline{m_s^N}) \overline{\mathcal{H} \phi}} \big\rangle ds\\
&\hspace{1cm} + \mathfrak{M}_t^N (\phi) +  u V_R \int_0^t \Big\langle \overline{m_s^N} \mathfrak{T}_2 (\overline{n_s^N}, \overline{m_s^N}, r) (\overline{n_s^N} - \overline{m_s^N}) \overline{Z}_s^N, \overline{\mathcal{H} \phi} \Big\rangle ds.
\end{align*}
From now on, we will choose the initial condition $\rho_0 = \lambda$ and $m_0 = n_0$ leading to $Z_0 = 0$ almost surely. At this stage, we are already able to read off the familiar terms of the SPDE in \Cref{theo:clt}. Again, for later analysis, we would like to eliminate the term involving $\mathcal{L}^N$ by choosing an appropriate test function.
For a time-dependent test function \[g : \{ (s,t) : 0 \leq s \leq t \leq T \} \rightarrow L^1 (\mathbb{T}^d \times [0,1] )\]
the representation becomes
\begin{align*}
\langle Z_t^N , g_{t,t} \rangle
&= \int_0^t \big\langle Z_s^N , \partial_s g_{s,t} + \mathcal{L}^N (g ) + \mathcal{R}_{n_s^N}^N (g_{s,t}) + \mu ( \mathcal{H} g_{s,t} - g_{s,t}) + uV_R \overline{\overline{m_s^N} r'(\overline{m_s^N}) \overline{\mathcal{H} g_{s,t}}} \big\rangle ds\\
&\hspace{1cm} + \int_0^t \blangle g_{s,t}, d \mathfrak{M}_s^N \brangle +  u V_R \int_0^t \Big\langle \overline{m_s^N} \mathfrak{T}_2 (\overline{n_s^N}, \overline{m_s^N}, r) (\overline{n_s^N} - \overline{m_s^N}) \overline{Z}_s^N, \overline{\mathcal{H} g_{s,t}} \Big\rangle ds.
\end{align*}
We now take $g^N$ to be the solution to
\begin{equation*}
\left\{
\begin{aligned}
- \partial_s g_{s,t}^N &= \mathcal{L}^N (g_{s,t}^N) + \mathcal{R}_{m_s^N}^N (g_{s,t}^N) + \mu \mathcal{H} g_{s,t}^N\\
& \hspace{0.5cm} - \mu g_{s,t}^N  + uV_R \overline{\overline{m_s^N} r'(\overline{m_s^N}) \overline{\mathcal{H} g_{s,t}^N}}, \\
 g_{t,t}^N &= \phi,
\end{aligned}
\right.
\end{equation*}
where we have replaced $\mathcal{R}_{n_s^N}^N$ by $\mathcal{R}_{m_s^N}^N$. Then we get
\begin{equation} \label{eq:fluctuationsrepresentation}
\begin{aligned}
\langle Z_t^N , g_{t,t}^N \rangle &= \int_0^t \blangle g_{s,t}^N, d \mathfrak{M}_s^N \brangle +\int_0^t \big\langle Z_s^N , \mathcal{R}_{n_s^N} (g_{s,t}^N) - \mathcal{R}_{m_s^N} (g_{s,t}^N) \big\rangle ds\\
&\hspace{1cm}  +  u V_R \int_0^t \Big\langle \overline{m_s^N}  \mathfrak{T}_2 (\overline{n_s^N}, \overline{m_s^N}, r) (\overline{n_s^N} - \overline{m_s^N}) \overline{Z}_s^N, \overline{\mathcal{H} g_{s,t}^N} \Big\rangle ds.
\end{aligned}
\end{equation}
We need to use $\mathcal{R}_{n_s^N}^N$ instead of $\mathcal{R}_{m_s^N}^N$ in the definition of $g^N$, so that the stochastic integral with respect to $\mathfrak{M}^N$ is well-defined.

We can regard the  martingales $\mathfrak{M}_t^N(\phi)$ as evaluations of a single martingale measure $\mathfrak{M}_t^N$. In \parencite{walsh_introduction_1986} a theory of stochastic integration with respect to martingale measures has been developed, which we will use to prove the central limit theorem.
\begin{definition}[{\parencite[Chapter 2]{walsh_introduction_1986}, \parencite[Definition 1.1]{Cho1995}}]
A martingale measure is a random set function $\mathfrak{M} : \mathbb{R}_+ \times \mathbb{T}^d \rightarrow \mathbb{R}$ such that
\begin{enumerate}
    \item 
    $\mathfrak{M}_t (S_1 \cup S_2 ) = \mathfrak{M}_t(S_1) + \mathfrak{M}_t(S_2) $ almost surely for any two disjoint Borel sets $S_1, S_2 \in \mathcal{B}( \mathbb{T}^d)$,
    \item
    $\mathfrak{M}_t(S)$ is a martingale for any Borel set $ S \in \mathcal{B} ( \mathbb{T}^d)$,
    \item
    there exists $(E_n, n \geq 0)$ such that $E_n$ is increasing with $\cup_{n \geq 1} E_n = \mathbb{T}^d$ and if $\mathcal{E}_n = \mathcal{B} (E_n)$ denotes the Borel-algebra, $\forall n \geq 0$ the supremum $ \sup_{S \in \mathcal{E}_n} \mathbb{E} (\mathfrak{M}_t(S)^2) $ is finite,
    \item for sequences $(S_n)_{n \geq 1} \subset \mathcal{E}_n$ with $S_n \downarrow \emptyset$, holds $\mathbb{E} (\mathfrak{M}_t (S_n)^2) \to 0$ as $n \to \infty$.
\end{enumerate}
\end{definition}
\begin{definition}[{\parencite[p. 291, Definition]{walsh_introduction_1986}}] \label{def:dominatingmeasure}
Let $\mathfrak{M}$ be a martingale measure. Suppose there exists a random measure $D$ on $\mathcal{B} ( \mathbb{R}_+ \times \mathbb{T}^d \times \mathbb{T}^d )$ such that
\begin{enumerate}
\item
$D$ is symmetric and positive definite in the second and third argument,
\item 
there exists a sequence $(E_n , n \geq 0 )$ exhausting $\mathbb{T}^d$ such that for all $n \geq 0$ the expectation $\mathbb{E} \big[ D([0,t] \times E_n \times E_n) \big]$ is finite,
\item 
for any $S_1 , S_2 \in \mathcal{B}(\mathbb{T}^d)$ the function $t \mapsto D([0,t] \times S_1 \times S_2)$ is predictable, 
\item
for any $S_1, S_2 \in \mathcal{B} (\mathbb{T}^d)$ and times $0 \leq s < t$
\begin{equation*}
    \Big\vert \blangle \mathfrak{M}_t(S_1), \mathfrak{M}_t(S_2) \brangle - \blangle \mathfrak{M}_s (S_1) , \mathfrak{M}_s(S_2) \brangle \Big\vert \leq D ((s, t] \times S_1 \times S_2 ) \enspace a.s.
\end{equation*}
\end{enumerate}
Then we call $\mathfrak{M}$ worthy and $D$ a dominating measure for $\mathfrak{M}$.
\end{definition}

\begin{lemma} \label{lem:mpfluctuations}
The martingale $\mathfrak{M}_t^N (\phi) = \sqrt{N \delta_N^{2-d}} M_t^N (\phi)$ is square-integrable with predictable quadratic variation
\begin{equation*}
\begin{aligned}
\Blangle \mathfrak{M}^N (\phi) \Brangle_t
&= \frac{u^2}{ \delta_N^{2 d}} \int_0^t \int_{\mathbb{T}^d} \frac{1}{\int_{B(x,\delta_N R)} n_s^N (y) dy} \int_{B(x, \delta_N R) \times [0,1]} \\
&\hspace{1cm} \frac{1}{\overline{n_s^N} (x, \delta_N R) +1} \Bigg[ \int_{B(x,\delta_N R)} \Big[ [1 + \delta_N^2 r_x(\overline{n_s^N}(x, \delta_N R)) ] \overline{n_s^N} (x, \delta_N R) \phi (z, k_0) \\
&\hspace{4cm} - \int_{[0,1]} \phi (z,k) \rho_s^N (z, dk) \Big] dz \Bigg]^2 \rho_s^N (y, dk_0) dy  dx ds
\end{aligned}
\end{equation*}
and $\mathfrak{M}^N$ is a worthy martingale measure.
\end{lemma}
\begin{proof}
The quadratic variation is given by
\begin{align*}
&\begin{aligned}
& \frac{N \delta_N^{2 - d} N}{\delta_N^{d + 2}} \int_{\mathbb{T}^d} \frac{1}{\int_{B(x,\delta_N R)} n_s^N(y) dy} \int_{B(x, \delta_N R) \times [0,1]} (\overline{n_s^N} (x, \delta_N R) +1)\\
&\hspace{0.5cm} \Bigg[ \Big\langle \rho_s^N +  \frac{u}{(\overline{n_s^N} (x, \delta_N R) +1)N} \mathds{1}_{B(x,\delta_N R)} (\cdot) ( \delta_{k_0}[1 + \delta_N^2 r_x(\overline{n_s^N} (x, \delta_N R)) ] \overline{n_s^N} (x, \delta_N R) - \rho) , \phi \Big\rangle\\
&\hspace{9cm} - \Big\langle \rho_s^N, \phi \Big\rangle \Bigg]^2  \rho_s^N (y, dk_0) dy dx \\
\end{aligned}\\
&\begin{aligned}
&= \frac{u^2}{ \delta_N^{2 d}} \int_0^t \int_{\mathbb{T}^d} \frac{1}{\int_{B(x,\delta_N R)} n_s^N(y) dy} \int_{B(x, \delta_N R) \times [0,1]} \frac{1}{\overline{n_s^N} (x, \delta_N R) +1} \\
&\hspace{0.5cm} \Bigg[ \int_{B(x,\delta_N R)} \Big[ [1 + \delta_N^2 r_x(\overline{n_s^N} (x, \delta_N R)) ] \overline{n_s^N} (x, \delta_N R) \phi (z, k_0) - \int_{[0,1]} \phi (z,k) \rho_s^N (z, dk) \Big] dz \Bigg]^2\\
& \hspace{10cm} \rho_s^N (y, dk_0) dy  dx ds.\\
\end{aligned}
\end{align*}
Using \Cref{assumption} 
a dominating measure $D^N$ satisfying \Cref{def:dominatingmeasure} is given by
\begin{align*}
\Big\langle D^N , \phi \otimes \phi \Big\rangle
&= \frac{u^2 n_{\max}^2}{\delta_N^{2d}} \int_0^t \int_{\mathbb{T}^d} \frac{1}{\int_{B(x,\delta_N R)} n_{s-}^N (y) dy} \int_{B(x, \delta_N R) \times [0,1]} \\
&\hspace{1cm} \Bigg[ \int_{B(x,\delta_N R)} \Big( \phi (z,k_0) + \int_{[0,1]} \phi (z,k) \rho_{s-}^N (z, dk) \Big) dz \Bigg]^2  \rho_{s-}^N (y, dk_0) dy dx ds.
\end{align*}
Here, $n_{\max}$ is the upper bound on the population size guaranteed by \Cref{assumption}.
\end{proof}
\begin{lemma}
\label{lem:boundondominatingmeasure}
The sequence of dominating measures  $D^N$ of \Cref{lem:mpfluctuations} 
are uniformly bounded in the sense that there exists a constant $C$ such that for any $0 \leq s \leq t$ and $\phi \in \mathcal{D} ( \mathbb{T}^d \times [0,1])$
\begin{equation}
\big\langle D^N , \mathds{1}_{[s,t]} \phi \otimes \phi \big\rangle \leq C \vert t - s \vert  \Vert \phi \Vert_2^2\enspace a.s.
\end{equation}
\end{lemma}
\begin{proof}
Indeed, we can bound the dominating measure by
\begin{align*}
\Big\langle D^N , \phi \otimes \phi \Big\rangle
&= \frac{u^2 n_{\max}^2}{\delta_N^{2d}} \int_0^t \int_{\mathbb{T}^d} \frac{1}{\int_{B(x,\delta_N R)} n_{s-}^N (y) dy} \int_{B(x, \delta_N R) \times [0,1]} \\
&\hspace{1cm} \Bigg[ \int_{B(x,\delta_N R)} \Big( \phi (z,k_0) + \int_{[0,1]} \phi (z,k) \rho_{s-}^N (z, dk) \Big) dz \Bigg]^2 \rho_{s-} (y, dk_0) dy dx ds \\
&\leq \frac{4 u^2 n_{\max}^2}{\delta_N^{2d}} \int_0^t \int_{\mathbb{T}^d} \Bigg[ \int_{B(x,\delta_N R)} \sup_{k \in [0,1]} \vert \phi (z, k ) \vert dz \Bigg]^2 dx ds\\
&\leq 4 u^2 n_{\max}^2 t \Vert \phi \Vert_2^2.
\end{align*}
The last step uses the identity (see \parencite[Proof of Lemma 3.4]{forien})
\begin{equation} \label{eq:identity_with_intersection}
\int_{(\mathbb{T}^d)^2} V_{\delta_N R} (z_1, z_2) \sup_{k \in [0,1]} \vert \phi (z_1, k ) \vert \sup_{k \in [0,1]} \vert \phi (z_2, k ) \vert dz_1 dz_2 \leq V_{\delta_N R}^2 \Vert \phi \Vert_2^2,
\end{equation}
where $V_R (z_1, z_2)$ denotes the volume of the intersection of two balls of radius $R$ centred at positions $z_1$ and $z_2$. In turn, the identity can be shown using the Cauchy-Schwarz inequality and $\int_{\mathbb{R}^d} V_r (x_1, x_2) dx_1 = V_r^2$.
\end{proof}


\subsection{Bound on the square of local averages of the fluctuations}
We end this section with a result, which will be the crucial ingredient for proving tightness of the sequence $Z^N$.
\begin{lemma} \label{lem:boundonZsquaredwithtypes}
Let $h : [0,1] \rightarrow \mathbb{R}$ be a measurable test function acting only on types and assume \[\Vert h \Vert_\infty = \sup_{k \in [0,1]} \vert h(k) \vert  < \infty.\]Then there exists a constant $K > 0$, independent of $h$, such that
\begin{equation*}
\sup_{0 \leq t \leq T} \sup_{x \in \mathbb{T}^d} \mathbb{E} \Bigg[ \Bigg( \frac{1}{V_{\delta_N R}} \int_{B(x,\delta_NR) \times [0,1]} h (k) Z_t^N(y,dk) dy \Bigg)^2 \Bigg] \leq K \frac{\Vert h \Vert_\infty^2}{V_{\delta_N R}}.
\end{equation*}
\end{lemma}
\begin{corollary} \label{cor:boundonZsquaredwithtypes}
Suppose $h : \mathbb{T}^d \times [0,1] \rightarrow \mathbb{R}$ and assume 
\begin{equation*}
\Vert h \Vert_2 = \Bigg( \int_{\mathbb{T}^d} \sup_{k \in [0,1]} \vert h (x,k) \vert^2 dx \Bigg)^{\frac{1}{2}} < \infty.
\end{equation*}
Then there exists a constant $K > 0$, independent of $h$, such that
\begin{equation*}
\sup_{0 \leq t \leq T} \int_{\mathbb{T}^d} \mathbb{E} \Bigg[ \Bigg( \frac{1}{V_{\delta_N R}} \int_{B(x,\delta_NR) \times [0,1]} h (x,k) Z_t^N(y,dk) dy \Bigg)^2 \Bigg] dx \leq K \frac{\Vert h \Vert_2^2}{V_{\delta_N R}}.
\end{equation*}
\end{corollary}
\begin{proof}[Proof of \Cref{lem:boundonZsquaredwithtypes}]
This proof can be seen as a generalization of \cite[Lemma 4.5]{forien_central_2017} to include types and to accomodate our setting. We will be using the representation
\begin{equation} \label{eq:lemZsquaredrepresentation}
\begin{aligned}
\langle Z_t^N , \phi \rangle &= \mathfrak{M}_t^N (\phi) + \int_0^t \big\langle Z_s^N , \mathcal{L}^N (\phi ) + \mathcal{R}_{n_s^N}^N (\phi) + \mu ( \mathcal{H} \phi - \phi) \big\rangle ds\\
&\hspace{1cm} + \sqrt{N \delta_N^{2-d}} u V_R \int_0^t \Big\langle m_s^N \lambda , \overline{ [r (\overline{n_s^N}) - r(\overline{m_s^N}) ] \overline{\phi}} (\cdot, \delta_N R) \Big\rangle ds,
\end{aligned}
\end{equation}
which we deduced in \Cref{subsec:representation_fluctuations}. Ultimately, we wish to set up an application of Gr\"onwall's inequality. We would like to substitute the solution to
\begin{equation*}
\left\{
\begin{aligned}
- \partial_s \zeta_{s,t}^N &= \mathcal{L}^N (\zeta_{s,t}^N) - \mu \zeta_{s,t}^N + \mu \mathcal{H} \zeta_{s,t}^N, \\
 \zeta^N_{t,t} &= \phi,
\end{aligned}
\right.
\end{equation*}
to remove the corresponding operators from the first integral. If we let $G_t^N$ be the semigroup of the operator $\mathcal{L}^N$, one can see that the differential equation is solved by
\begin{equation*}
\begin{aligned}
\zeta_{s,t}^N &= G_{t-s}^N * \Big( e^{- \mu (t-s)} \phi + \big(1 - e^{- \mu (t-s)} \big) \mathcal{H} \phi  \Big).
\end{aligned}
\end{equation*}
We intend to choose the test function $\phi^N (z,k) = h(k) \frac{1}{V_{\delta_N R}} \mathds{1}_{\vert x- z \vert < \delta_N R} $, so that the evaluation $\langle Z_t^N , \phi \rangle $ will lead to the average in the statement of the lemma. This leads to 
\begin{align*}
&G_{t-s}^N * \Big( e^{- \mu (t-s)} \phi^N + \big(1 - e^{- \mu (t-s} \big) \mathcal{H} \phi^N  \Big) (z,k)\\
&= \Bigg( G_{t-s}^N * \frac{1}{V_{\delta_N R}} \mathds{1}_{\vert x - \cdot \vert < \delta_N R} \Bigg) (z) \Big( e^{- \mu (t-s)} h(k) + \big(1 - e^{- \mu (t-s)} \big) \mathcal{H} h \Big).
\end{align*}
Let us define 
\begin{equation*}
h_{t-s} (k) := \Big( e^{- \mu (t-s)} h(k) + \big(1 - e^{- \mu (t-s)} \big) \mathcal{H} h \Big)
\end{equation*}
and for $n, m \geq 0$ the quotient
\begin{equation} \label{eq:quotient}
q_z (n, m ) := 
\begin{dcases}
\frac{r_z(n ) - r_z(m)}{n - m} & n \neq m, \\
r_z'(n) & n = m.
\end{dcases}
\end{equation}
Inserting $\zeta_{s,t}^N$ into \eqref{eq:lemZsquaredrepresentation} and using $\mathcal{R}_{m_s^N}^N (\zeta_{s,t}^N) = u V_R \overline{r(\overline{m_s^N}) \overline{\zeta_{s,t}^N}} (\cdot, \delta_N R)$ yields
\begin{align*}
\langle Z_t^N , \phi \rangle 
&= \int_0^t \int_{[0,t] \times \mathbb{T}^d \times [0,1]} \zeta_{s,t}^N (z,k) \mathfrak{M}^N (ds\;dz\;dk) + u V_R\int_0^t \big\langle \overline{Z_s^N} , r(\overline{n_s^N}) \overline{\zeta_{s,t}^N} \rangle ds\\
&\hspace{1cm} + u V_R \int_0^t \Big\langle \overline{m_s^N \lambda} , \overline{F_s^N} q(\overline{n_s^N} , \overline{m_s^N} ) \overline{\zeta_{s,t}^N} \Big\rangle ds.
\end{align*}
If we now choose $\phi^N (z,k) = h(k) \frac{1}{V_{\delta_N R}} \mathds{1}_{\vert x- z \vert < \delta_N R} $, we obtain
\begin{align*}
\langle Z_t^N , \phi^N \rangle
&= \int_{[0,t] \times \mathbb{T}^d \times [0,1]} \overline{G_{t-s}^N} (x - z) h_{t-s} (k) \mathfrak{M}^N (ds\;dz\;dk)\\
&\hspace{1cm} + u V_R\int_0^t \big\langle \overline{Z_s^N} , r(\overline{n_s^N}) \overline{\overline{G_{t-s}^N}} (x - \cdot) h_{t-s} \rangle ds \\
&\hspace{1cm} + u V_R \int_0^t \Big\langle \overline{m_s^N \lambda} , \overline{F_s^N} q(\overline{n_s^N} , \overline{m_s^N} ) \overline{\overline{G_{t-s}^N}} (x- \cdot ) h_{t-s} \Big\rangle ds.
\end{align*}
Applying the inequality
\begin{equation*}
\Big( \sum_{i = 1}^n a_i \Big)^2 \leq n \Big( \sum_{i = 1}^n a_i^2\Big),
\end{equation*}
we arrive at
\begin{equation} \label{eq:fluctuationsaftersquaretrick}
\begin{aligned}
\langle Z_t^N , \phi^N \rangle^2
&\leq 3 \Bigg( \int_{[0,t] \times \mathbb{T}^d \times [0,1]} \overline{G_{t-s}^N} (x - z) h_{t-s} (k) \mathfrak{M}^N (ds\;dz\;dk) \Bigg)^2\\
&\hspace{0.5cm} + 3 u^2 V_R^2 \Bigg( \int_0^t \int_{\mathbb{T}^d \times [0,1]} \overline{Z_s^N}(z, dk) r(\overline{n_s^N} (z)) \overline{\overline{G_{t-s}^N}} (x - z) h_{t-s} (k) dz ds \Bigg)^2 \\
&\hspace{0.5cm} + 3u^2 V_R^2 \Bigg( \int_0^t \int_{\mathbb{T}^d} \overline{m_s^N} (z) \overline{F_s^N} (z) q(\overline{n_s^N}(z) , \overline{m_s^N}(z) ) \overline{\overline{G_{t-s}^N}} (x- z) \mathcal{H} h \; dz ds \Bigg)^2.
\end{aligned}
\end{equation}
One is able to apply Jensen's inequality as $\int_{\mathbb{T}^d} G_{t-s}^N (z) dz = 1$, which we will use to treat the second and third term. For the second term we get the bound
\begin{align*}
3 u^2 V_R^2 \Vert r \Vert_\infty^2 t \int_0^t \int_{\mathbb{T}^d} \Bigg( \int_{[0,1]} h_{t-s} (k) \overline{Z_s^N}(z, dk) \Bigg)^2 \overline{\overline{G_{t-s}^N}} (x - z) dz ds.
\end{align*}
In turn, using Jensen's inequality
\begin{align*}
\Bigg( \int_{[0,1]} h_{t-s} (k) \overline{Z_s^N}(z, dk) \Bigg)^2
&= \Bigg(e^{- \mu (t-s)} \int_{[0,1]} h(k) \overline{Z_s^N} (z,dk) + \big(1 - e^{- \mu (t-s)} \big) \mathcal{H}h \cdot \overline{F_s^N} (z) \Bigg)^2\\
&\leq \Bigg( \int_{[0,1]} h(k) \overline{Z_s^N} (z,dk) \Bigg)^2 + \Vert h \Vert_\infty^2 \overline{F_s^N} (z)^2,
\end{align*}
where we used $\Vert h \Vert_1 \leq \Vert h \Vert_\infty$.
The third term of \eqref{eq:fluctuationsaftersquaretrick} results in
\begin{align*}
    3u^2 V_R^2 \Vert h \Vert_\infty^2 n_{\max}^2 \Vert r' \Vert_\infty^2 t \int_0^t \int_{\mathbb{T}^d} \overline{F_s^N} (z)^2 \overline{\overline{G_{t-s}^N}} (x- z) dz ds,
\end{align*}
where the square of the quotient \eqref{eq:quotient} was bounded by $\sup_{z \in \mathbb{T}^d} \Vert r_z' \Vert_\infty^2$.
Taking expectations on both sides of \eqref{eq:fluctuationsaftersquaretrick} and exchanging integral and expectation, we arrive at
\begin{align*}
&\mathbb{E} \Bigg[ \Bigg( \int_{[0,1]} h(k) \overline{Z_t^N} (x,dk) \Bigg)^2 \Bigg]\\
&\leq 3 \mathbb{E} \Bigg[ \Bigg( \int_{[0,t] \times \mathbb{T}^d \times [0,1]} \overline{G_{t-s}^N} (x - z) h_{t-s} (k) \mathfrak{M}^N (ds\;dz\;dk) \Bigg)^2 \Bigg]\\
&\hspace{1cm} + 3 u^2 V_R^2 \Vert r \Vert_\infty^2 t \int_0^t \int_{\mathbb{T}^d} \mathbb{E} \Bigg[ \Bigg( \int_{[0,1]} h(k) \overline{Z_s^N} (z,dk) \Bigg)^2 \Bigg] \overline{\overline{G_{t-s}^N}} (x - z) dz ds\\
&\hspace{1cm} + 3 u^2 V_R^2 \Vert h \Vert_\infty^2 \Vert r \Vert_\infty^2  t \int_0^t \int_{\mathbb{T}^d}  \mathbb{E} \big[ 
\overline{F_s^N} (z)^2 \big] \overline{\overline{G_{t-s}^N}} (x - z) dz ds \\
&\hspace{1cm} + 3u^2 V_R^2 \Vert h \Vert_\infty^2 n_{\max}^2 \Vert r' \Vert_\infty^2 t \int_0^t \int_{\mathbb{T}^d} \mathbb{E} \big[ \overline{F_s^N} (z)^2 \big] \overline{\overline{G_{t-s}^N}} (x- z) dz ds.
\end{align*}
By \Cref{lem:boundondominatingmeasure}, the term involving the martingale can be bounded by 
\begin{align*}
&\mathbb{E} \Bigg[ \Bigg( \int_0^t \int_{\mathbb{T}^d \times [0,1]} \overline{G_{t-s}^N} (x - z) h_{t-s} (k) \mathfrak{M}^N (dz,dk) \Bigg)^2 \Bigg]\\
&\leq C \int_0^t \int_{\mathbb{T}^d} \sup_{k \in [0,1]} \Big\vert \overline{G_{t-s}^N} (x - z) h_{t-s} (k) \Big\vert^2 dz ds \\
&\leq C \Vert h \Vert_\infty^2 \int_0^t \Big\Vert \overline{G_{t-s}^N} (x - z) \Big\Vert^2_2 ds.
\end{align*}
If we let $(X_{t-s}^N)_{0 \leq s \leq t}$ be a random walk started in the origin with generator $\mathcal{L}^N$, we can write 
\begin{align*}
\Big\Vert \overline{G_{t-s}^N} (x - z) \Big\Vert^2_2 &= \int_0^t \int_{\mathbb{T}^d} \mathbb{E}_0 \Bigg[ \frac{1}{V_{\delta_N R}} \mathds{1}_{\vert X_{t-s}^N - (x-z) \vert < \delta_N R} \Bigg]^2 dz ds\\
&\leq \int_0^t  \mathbb{E}_0 \Bigg[ \int_{\mathbb{T}^d} \Bigg( \frac{1}{V_{\delta_N R}} \mathds{1}_{\vert X_{t-s}^N - (x-z) \vert < \delta_N R} \Bigg)^2 dz \Bigg] ds\\
&= \frac{t}{V_{\delta_N R}}.
\end{align*}
Summarizing, we obtain
\begin{equation} \label{eq:last_annoying_step}
\begin{aligned}
&\mathbb{E} \Bigg[ \Bigg( \int_{[0,1]} h(k) \overline{Z_t^N} (x,dk) \Bigg)^2 \Bigg]\\
&\leq 3 C \Vert h \Vert_\infty^2 \frac{t^2}{V_{\delta_N R}} + 3 u^2 V_R^2 \Vert r \Vert_\infty^2 t \int_0^t \sup_{z \in \mathbb{T}^d}  \mathbb{E} \Bigg[ \Bigg( \int_{[0,1]} h(k) \overline{Z_s^N} (z,dk) \Bigg)^2 \Bigg] ds\\
&\hspace{1cm} + 3 u^2 V_R^2 t^2 \Vert h \Vert_\infty^2 \big( \Vert r \Vert_\infty^2 + n_{\max}^2 \Vert r' \Vert_\infty^2 \big) \sup_{s \in [0,t]} \sup_{z \in \mathbb{T}^d} \mathbb{E} \big[ \overline{F_s^N} (z)^2 \big].
\end{aligned}
\end{equation}
In turn, by \Cref{lem:boundonFsquared}, there exists a constant $C' > 0$ such that
\[ \sup_{s \in [0,t]} \sup_{z \in \mathbb{T}^d} \mathbb{E} \big[ \overline{F_s^N} (z)^2 \big] \leq \frac{C'}{V_{\delta_N R}}. \]
As the right-hand side of \eqref{eq:last_annoying_step} does not depend on $x$, we can take a supremum on the left-hand side
\begin{align*}
\sup_{x \in \mathbb{T}^d} \mathbb{E} \Bigg[ \Bigg( \int_{[0,1]} h(k) \overline{Z_t^N} (x,dk) \Bigg)^2 \Bigg]
&\leq \frac{\Vert h \Vert_\infty^2}{V_{\delta_N R}} 3 t^2 \Big( C + u^2 V_R^2 \big( \Vert r \Vert_\infty^2 + n_{\max}^2 \Vert r' \Vert_\infty^2 \big) C' \Big) \\
&\hspace{1cm} + 3 u^2 V_R^2 \Vert r \Vert_\infty^2 t \int_0^t \sup_{x \in \mathbb{T}^d}  \mathbb{E} \Bigg[ \Bigg( \int_{[0,1]} h(k) \overline{Z_s^N} (x,dk) \Bigg)^2 \Bigg] ds.
\end{align*}
The function
\[ t \mapsto \sup_{x \in \mathbb{T}^d} \mathbb{E} \Bigg[ \Bigg( \int_{[0,1]} h(k) \overline{Z_t^N} (x,dk) \Bigg)^2 \Bigg] \]
is bounded as 
\[ \Bigg\vert \int_{[0,1]} h(k) \overline{Z_t^N} (x,dk) \Bigg\vert \leq 2 \Vert h \Vert_\infty \sqrt{N \delta_N^{2-d}} n_{\max}. \]
Hence, by Gr\"onwall's inequality
\begin{align*}
\sup_{x \in \mathbb{T}^d} \mathbb{E} \Bigg[ \Bigg( \int_{[0,1]} h(k) \overline{Z_t^N} (x,dk) \Bigg)^2 \Bigg]
&\leq \frac{\Vert h \Vert_\infty^2}{V_{\delta_N R}} \Big[ 3 T^2 \Big( C + u^2 V_R^2 \big( \Vert r \Vert_\infty^2 + n_{\max}^2 \Vert r' \Vert_\infty^2 \big) C' \Big) \Big] e^{3 u^2 V_R^2 \Vert r \Vert_\infty^2 T t}\\
&\leq C'' \frac{\Vert h \Vert_\infty^2}{V_{\delta_NR }}. \qedhere
\end{align*}
\end{proof}

\section{Proof of the central limit theorem} \label{sec:proofclt}
Recall the representation \eqref{eq:fluctuationsrepresentation}
\begin{align*}
\langle Z_t^N , \phi \rangle &= \int_0^t \blangle g_{s,t}^N, d \mathfrak{M}_s^N \brangle +\int_0^t \big\langle Z_s^N , \mathcal{R}_{n_s^N} (g_{s,t}^N) - \mathcal{R}_{m_s^N} (g_{s,t}^N) \big\rangle ds\\
&\hspace{1cm}  +  u V_R \int_0^t \Big\langle \overline{m_s^N}  \mathfrak{T}_2 (\overline{n_s^N}, \overline{m_s^N}, r) (\overline{n_s^N} - \overline{m_s^N}) \overline{Z}_s^N, \overline{\mathcal{H} g_{s,t}^N} \Big\rangle ds.
\end{align*}
The first part will be the essential one, the second and third part are error terms which will vanish. The proof classically separates into: 
\begin{enumerate}[leftmargin=4\parindent]
\item[Tightness]
We show compact containment and control of the modulus of continuity of the fluctuations to apply the Aldous-Rebolledo criterion.
\item[Uniqueness]
We will show convergence of the finite-dimensional distributions by proving the convergence of $g^N$, as well as the convergence of the martingale measures $\mathfrak{M}^N$ to a Gaussian random field.
\end{enumerate}
Our strategy will be an adaptation and synthesis of the ideas of \cite{forien_central_2017} and \cite{forien}.
\subsection{Properties of the functions \texorpdfstring{$g^N$}{}}
An important prerequisite is the control of the behaviour of the functions $g^N$. Recall that $g^N$ is the solution to
\begin{equation*}
\left\{
\begin{aligned}
- \partial_s g_{s,t}^N &= \mathcal{L}^N (g_{s,t}^N) + \mathcal{R}_{m_s^N} (g_{s,t}^N) - \mu g_{s,t}^N\\
& \hspace{0.5cm} + \mu \mathcal{H} g_{s,t}^N + uV_R \overline{\overline{m_s^N} r'(\overline{m_s^N}) \overline{\mathcal{H} g_{s,t}^N}}, \\
g_{t,t}^N &= \phi,
\end{aligned}
\right.
\end{equation*}
where $\mathcal{H} g^N = \int_0^1 g^N (\cdot, k) dk$. Unsurprisingly, the limit $g$ is the solution to 
\begin{equation*}
\left\{
\begin{aligned}
- \partial_s g_{s,t} &= u V_R \Big( \frac{R^2}{d+2} \Delta + r(n_s) \Big) g_{s,t} - \mu g_{s,t}\\
& \hspace{0.5cm} + \mu \mathcal{H} g_{s,t} + uV_R n_s r'(n_s) \mathcal{H} g_{s,t}, \\
 g_{t,t} &= \phi.
\end{aligned}
\right.
\end{equation*}
Let us remind ourselves that in \eqref{eq:pnormdefinition} we defined the $L^p$-norm to involve the supremum over types. The proofs of the following lemmas can be found in \Cref{appendix:propertiesofthefunctions}. The first result provides uniform bounds on the $L^p$-norms of the functions and derivatives.

\begin{lemma}
\label{lem:semigroupsbounded}

For $p \in [1,2] $ and $\beta \in \mathbb{N}^d, \vert \beta \vert_1 \leq 4$, for any $T > 0$ there exists constants $C, C' > 0$ with $C$ being independent of $\phi$ such that
\begin{align*}
\sup_{0 \leq s \leq t \leq T} \big\Vert g_{s,t}^N \big\Vert_p &\leq C \Vert \phi \Vert_p, \\
\sup_{0 \leq s \leq t \leq T} \big\Vert \partial_\beta g_{s,t}^N \big\Vert_p &\leq C'.
\end{align*}
\end{lemma}
The next lemma confirms our intuition that $g$ should be the limit of $g^N$ as $N \to \infty$.
\begin{lemma} \label{lem:semigroupconvergence}
For any $T > 0$, there exists a constant $C > 0$, depending on $\phi$, such that for $p \in [ 1,2]$
\begin{equation*}
\sup_{0 \leq s \leq t \leq T} \Big\Vert g_{s,t}^N - g_{s,t} \Big\Vert_p \leq C \delta_N.
\end{equation*}
\end{lemma}
We end this interlude with a characterisation of the continuity of $g^N$ in the second time component.
\begin{lemma} \label{lem:semigroupcontinuity}
For any $T > 0$, there exists $C > 0$ depending on $\phi$ such that for $N \geq 1, 0 \leq s , t,t ' \leq T$ and $p \in [1,2]$
\begin{equation*}
\Big\Vert  g_{s,t'}^N - g_{s,t}^N \Big\Vert_p \leq C \vert t' - t \vert.
\end{equation*}
\end{lemma}

\subsection{Tightness} \label{subsec:tightness}
By \Cref{theo:convergenceofstochasticprocesses}, the tightness of the processes $(Z^N)_{N \geq 1}$ can be reduced to the tightness of $((\langle Z_t^N , \phi \rangle)_{t \geq 0}, N \geq 1)$ for any $\phi \in \mathcal{D} (\mathbb{T}^d \times [0,1])$. We will conclude by Aldous' criterion \cite{aldouscriterion}.
The bounds on the error terms we develop in the compact containment proof will also show that these terms vanish in $L^1$ uniformly in time as $N \to \infty$.

\subsubsection{Compact containment and vanishing error terms}

\begin{lemma} \label{lem:compactcontainment}
 There exists a constant $C > 0$, independent of $N \geq 1$, such that for any $\phi : \mathbb{T}^d \times [0,1] \rightarrow \mathbb{R}$ with $\sup_{p \in \{1, 2 \} } \Vert \phi \Vert_p < 1$
\begin{equation*}
\mathbb{E} \Bigg[ \sup_{0 \leq t \leq T} \big\vert \blangle Z_t^N, \phi \brangle \big\vert \Bigg] \leq C.
\end{equation*}
\end{lemma}

\begin{proof}[Proof of \Cref{lem:compactcontainment}]
Again, for convenience, let us recall the representation from \eqref{eq:fluctuationsrepresentation} given by 
\begin{align*}
\langle Z_t^N , \phi \rangle &= \int_0^t \blangle g_{s,t}^N, d \mathfrak{M}_s^N \brangle +\int_0^t \big\langle Z_s^N , \mathcal{R}_{n_s^N} (g_{s,t}^N) - \mathcal{R}_{m_s^N} (g_{s,t}^N) \big\rangle ds\\
&\hspace{1cm}  +  u V_R \int_0^t \Big\langle \overline{m_s^N}  \mathfrak{T}_2 (\overline{n_s^N}, \overline{m_s^N}, r) (\overline{n_s^N} - \overline{m_s^N}) \overline{Z}_s^N, \overline{\mathcal{H} g_{s,t}^N} \Big\rangle ds.
\end{align*}
For the integral with respect to the martingale measure we can use \Cref{theo:tightnessstochasticintegrals}, whose prerequisites have been verified in \Cref{lem:boundondominatingmeasure}, \Cref{lem:semigroupsbounded} and \Cref{lem:semigroupcontinuity}. To treat the second integral, recall that the remainder in Taylor expansion was defined in \Cref{eq:taylor_remainder}.
We can write for the integrand of the second integral for $\psi : \mathbb{T}^d \times [0,1] \rightarrow \mathbb{R}$
\begin{align*}
&\big\langle Z_s^N , \overline{\big[ r (\overline{n_s^N}) - r (\overline{m_s^N} ) \big] \overline{\psi}}\big\rangle \\
&= \big\langle \overline{Z_s^N} , \big[ r (\overline{n_s^N}) - r (\overline{m_s^N} ) \big] \overline{\psi} \big\rangle\\
&= \int_{\mathbb{R}^d \times [0,1]} \overline{Z_s^N} (x, dk ) (\overline{n_s^N} (x) - \overline{m_s^N} (x)) \mathfrak{T}_1(\overline{n_s^N} (x), \overline{m_s^N} (x), r_x) \overline{\psi} (x,k) dx\\
&= \frac{1}{\sqrt{N \delta_N^{2-d}}} \int_{\mathbb{T}^d} \int_{[0,1]} \overline{\psi} (x,k) \overline{Z_s^N} (x, dk ) \overline{F_s^N} (x) \mathfrak{T}_1(\overline{n_s^N}(x), \overline{m_s^N} (x), r_x)  dx\\
&\leq \frac{1}{\sqrt{N \delta_N^{2-d}}} \Bigg( \int_{\mathbb{T}^d} \Bigg( \int_{[0,1]} \overline{\psi} (x,k) \overline{Z_s^N} (x, dk ) \Bigg)^2 dx \Bigg)^{1/2} \Bigg( \int_{\mathbb{T}^d} \overline{F_s^N} (x)^2 \mathfrak{T}_1(\overline{n_s^N} (x), \overline{m_s^N} (x), r_x)^2  dx \Bigg)^{1/2}.
\end{align*}
Taking the expectation on both sides and using the Cauchy-Schwarz inequality again
\begin{align*}
&\mathbb{E} \Big[ \big\vert \big\langle Z_s^N , \overline{\big[ r (\overline{n_s^N}) - r (\overline{m_s^N} ) \big] \overline{\psi}} \big\rangle \big\vert \Big] \\
&\leq \frac{1}{\sqrt{N \delta_N^{2-d}}} \mathbb{E} \Bigg[ \int_{\mathbb{T}^d} \Bigg( \int_{[0,1]} \overline{\psi} (x,k) \overline{Z_s^N} (x, dk ) \Bigg)^2 dx \Bigg]^{1/2} \mathbb{E} \Bigg[ \int_{\mathbb{T}^d} \overline{F_s^N} (x)^2 \mathfrak{T}_1(\overline{n_s^N} (x), \overline{m_s^N} (x), r_x)^2  dx \Bigg]^{1/2} \\
&\leq \frac{1}{\sqrt{N \delta_N^{2-d}}}  \frac{\sqrt{K} \Vert \psi \Vert_2}{\sqrt{V_{\delta_N R}}} \cdot \sup_{z \in \mathbb{T}^d} \Vert r_z' \Vert_\infty Vol (\mathbb{T}^d)^{1/2} \Bigg( \sup_{x \in \mathbb{T}^d} \mathbb{E} \Big[ \overline{F_s^N} (x)^2 \Big] \Bigg)^{1/2} \\
&\leq \Big( N^{1/2} \delta_N^{1+ \frac{d}{2}} \Big)^{-1} \frac{\sqrt{KC'}}{V_R} \sup_{z \in \mathbb{T}^d} \Vert r_z' \Vert_\infty Vol (\mathbb{T}^d)^{1/2} \Vert \psi \Vert_2 .
\end{align*}
Here, we used \Cref{cor:boundonZsquaredwithtypes} and \Cref{lem:boundonFsquared} and denote the constants in the corresponding bounds by $K$ and $C'$ respectively. Returning to the integral over time
\begin{align*}
&\mathbb{E} \Bigg[ \sup_{0 \leq t \leq T} \Big\vert \int_0^t \blangle Z_s^N, \overline{\big[ r (\overline{n_s^N}) - r (\overline{m_s^N} ) \big] \overline{g_{s,t}^N}} (\cdot, \delta_N R) \brangle ds \Big\vert \Bigg]\\
&\leq \mathbb{E} \Bigg[  \int_0^T \Big\vert \blangle Z_s^N, \overline{\big[ r (\overline{n_s^N}) - r (\overline{m_s^N} ) \big] \overline{g_{s,t}^N}} (\cdot, \delta_N R) \brangle \Big\vert ds \Bigg]\\
&\leq T \Big( N^{1/2} \delta_N^{1+ \frac{d}{2}} \Big)^{-1} \frac{\sqrt{KC'}}{V_R}  \sup_{z \in \mathbb{T}^d} \Vert r_z' \Vert_\infty Vol (\mathbb{T}^d)^{1/2} \sup_{0 \leq s \leq T} \Big\Vert \overline{ g_{s,t}^N } \Big\Vert_2.
\end{align*}
We assumed in \eqref{eq:assumptionondelta} that $(\delta_N)_{N \geq 1}$ is chosen in such a way that $N^{-1/2} \delta_N^{-(1+ \frac{d}{2})}$ vanishes. The norm of the average of $g^N$ can be bounded noting $\Vert \overline{\psi} \Vert_p \leq \Vert \psi \Vert_p$ and
\begin{align*}
\int_{\mathbb{T}^d} \sup_{k \in [0,1]} \big\vert \overline{g_{s,t}^N} \big\vert^2 dx
&= \int_{\mathbb{T}^d} \Big( \overline{ \sup_{k \in [0,1]} \big\vert g_{s,t}^N \big\vert} \Big)^2 dx \leq \big\Vert g_{s,t}^N \big\Vert_2^2,
\end{align*}
which is bounded due to \Cref{lem:semigroupsbounded}.

For the integrand of the third integral we can use 
\begin{align*}
&\Big\langle \overline{m_s^N}  \mathfrak{T}_2 (\overline{n_s^N}, \overline{m_s^N}, r) (\overline{n_s^N} - \overline{m_s^N}) \overline{Z}_s^N, \overline{\mathcal{H} g_{s,t}^N} \Big\rangle\\
&= \frac{1}{\sqrt{N \delta_N^{2 - d}}} \int_{\mathbb{T}^d} \overline{m_s^N}  \mathfrak{T}_2 (\overline{n_s^N} (x), \overline{m_s^N} (x), r_x) \Big[ \overline{F}_s^N (x) \Big]^2 \overline{\mathcal{H} g_{s,t}^N} (x) dx,
\end{align*}
which can be bounded similarly to the previous term. As $\mathcal{H} g^N = \int_0^1 g^N dk$ we have no dependence on the types, which simplifies the calculation. We are thus able to bound the expectation of the supremum over time of the three parts of the representation \eqref{eq:fluctuationsrepresentation}, which concludes the proof.
\end{proof}

\subsubsection{Bound on the modulus of continuity}

\begin{lemma} \label{lem:modulusofcontinuity}
For any function $\phi \in \mathcal{D} (\mathbb{T}^d \times [0,1])$, sequence of stopping times $(T_N)_{N \geq 1}$ in $[0,T]$ and sequence of positive reals $(\gamma_N)_{N \geq 1}$ converging to zero,
\begin{equation*}
\langle Z_{T_N + \gamma_N}^N , \phi \rangle - \langle Z_{T_N}^N, \phi \rangle \xrightarrow[N \to \infty]{} 0
\end{equation*}
in probability.
\end{lemma}
\begin{proof}
We will use a strategy from \cite[Theorem 7.13]{walsh_introduction_1986} as in \cite[Proposition 4.8]{forien_central_2017}. We extend $g^N$ from $0 \leq s \leq t \leq T$ to arbitrary $s,t \in [0,T]$ by defining
\begin{equation*}
g_{s,t}^N (x,k) := g_{s \land t , t}^N (x,k),
\end{equation*}
and define
\begin{equation*}
V_t^N = \int_0^T \langle g_{s,t}^N, d \mathfrak{M}_s^N \rangle.
\end{equation*}
We can then write
\begin{align*}
&\int_0^{T^N + \gamma^N} \blangle g_{s,T^N + \gamma^N}^N, d \mathfrak{M}_s^N \brangle - \int_0^{T^N} \blangle g_{s,T^N}^N, d \mathfrak{M}_s^N \brangle\\
&= ( V_{T^N + \gamma^N}^N - V_{T^N}^N ) + \int_{T^N}^{T^N + \gamma^N} \blangle \phi , d \mathfrak{M}_s^N \brangle,
\end{align*}
where we used $g^N (x,k,s,t) = \phi$ for $s \geq t$. Thus, the representation \eqref{eq:fluctuationsrepresentation} results in
\begingroup
\allowdisplaybreaks
\begin{align*}
&\langle Z_{T^N + \gamma^N}^N , \phi \rangle - \langle Z_{T^N}^N , \phi \rangle \\
&= (V_{T^N + \gamma^N}^N - V_{T^N}^N ) + \int_{T^N}^{T^N + \gamma^N} \blangle \phi, d \mathfrak{M}_s^N \brangle\\
&\hspace{1cm} + u V_R\int_0^{T^N + \gamma^N} \big\langle Z_s^N , \overline{\big[ r (\overline{n_s^N}) - r (\overline{m_s^N} ) \big] \overline{g_{s,T^N + \gamma^N}^N}} \big\rangle ds\\
&\hspace{1cm}  +  u V_R \int_0^{T^N + \gamma^N} \Big\langle \overline{m_s^N}  \mathfrak{T}_2 (\overline{n_s^N}, \overline{m_s^N}, r) (\overline{n_s^N} - \overline{m_s^N}) \overline{Z}_s^N, \overline{\mathcal{H} g_{s,T^N + \gamma^N}^N} \Big\rangle ds\\
&\hspace{1cm} - u V_R\int_0^{T^N} \big\langle Z_s^N , \overline{\big[ r (\overline{n_s^N}) - r (\overline{m_s^N} ) \big] \overline{g_{s,T^N}^N}} \big\rangle ds\\
&\hspace{1cm} -  u V_R \int_0^{T^N} \Big\langle \overline{m_s^N}  \mathfrak{T}_2 (\overline{n_s^N}, \overline{m_s^N}, r) (\overline{n_s^N} - \overline{m_s^N}) \overline{Z}_s^N, \overline{\mathcal{H}g_{s,T^N}^N} \Big\rangle ds.
\end{align*}
\endgroup
Note that the four error terms converge to zero in $L^1$ using the same bounds as in the compact containment proof. By \Cref{lem:boundondominatingmeasure} and \Cref{lem:semigroupcontinuity} we obtain for any $t' \geq t$
\begin{align*}
\mathbb{E} \Big[ \big\vert V_{t'}^N - V_{t}^N \big\vert^2 \Big] &= \mathbb{E} \Bigg[ \Bigg( \int_0^T \int_{\mathbb{R}^d} (g_{s,t'}^N (x,k) - g_{s,t}^N (x,k)) M^N (dx dk ds) \Bigg)^2 \Bigg] \\
&\leq C \int_0^T \Big\Vert g_{s,t'}^N - g_{s,t}^N  \Big\Vert_2^2 ds \\
&\leq C' T \vert t' - t \vert^2.
\end{align*}
Hence, by Kolmogorov's continuity criterion, there exists a random variable $B^N$ such that for any $0 \leq t, t' \leq T, 0 < \alpha < 1/2$,
\begin{equation*}
\vert V_{t'}^N - V_{t}^N \vert \leq B^N \vert t' - t \vert^{\alpha} 
\end{equation*}
and $\mathbb{E} \big[ (B^N)^2 \big] < C''$, where the constant $C'' \geq 0$ is independent of $N \geq 1$. Therefore, we get for $\alpha = 1/4$
\begin{equation*}
\mathbb{E} \Big[ \big\vert V_{T^N + \gamma^N}^N - V_{T^N}^N \big\vert^2 \Big] \leq \mathbb{E} \big[ (B^N)^2 \big] \left( \gamma^N \right)^{1/2} \leq C'' \left( \gamma^N \right)^{1/2}.
\end{equation*}
Additionally, again by \Cref{lem:boundondominatingmeasure}
\begin{equation*}
\begin{aligned}
\mathbb{E} \Bigg[ \Bigg( \int_{T^N}^{T^N + \gamma^N} \blangle \phi, d \mathfrak{M}_s^N \brangle \Bigg)^2 \Bigg] &\leq C \mathbb{E} \Bigg[ \int_{T^N}^{T^N + \gamma^N} \Vert \phi\Vert_2^2 ds \Bigg] \\
&\leq C \gamma^N \Vert \phi \Vert_2^2.
\end{aligned}
\end{equation*}
Both terms converge to zero in $L^2$ and thus in probability.
\end{proof}
\Cref{lem:compactcontainment} gives, by the Markov inequality, the tightness of fixed time marginals. Together with \Cref{lem:modulusofcontinuity} this amounts precisely to the requisites of Aldous' criterion \cite{aldouscriterion} and tightness follows.
\subsection{Convergence of the martingale measures}
The convergence of $(\mathfrak{M}^N)_{N \geq 1}$ is the last missing ingredient for the convergence of the finite-dimensional distributions, as we have already convergence of the integrands $g^N$ of the stochastic integrals with respect to $\mathfrak{M}^N$.
\begin{theorem}[{\parencite[VIII Theorem 3.11]{JS03}}] \label{theo:convergenceofmartingales}
Suppose $((X_t^N)_{t \in [0,T]}, N \geq 1)$ is a sequence of local martingales in $\mathbb{D}(\mathbb{R}_+, \mathbb{R}^d)$ such that 
\begin{enumerate}
    \item 
there exists $C > 0$, such that for any $N \geq 1$, $  \big\vert X_t^N - X_{t-}^N \big\vert \leq C$,
\item
the jumps vanish, i.e.
$   \sup_{t \in [0,T]} \big\vert X_t^N - X_{t-}^N \big\vert \xrightarrow[N \to \infty]{P} 0,$
\item 
there exists a $d$-dimensional, continuous Gaussian martingale $(X_t)_{t \geq 0}$ such that for any rational time $t$ the quadratic covaration
$\big\langle X^{i, N}, X^{j, N} \big\rangle_t \xrightarrow[N \to \infty]{} \big\langle X^i , X^j \big\rangle_t$
in probability. Here, $i$ and $j$ indicate the $i$-th and $j$-th component of the martingales $X^N$ and $X$.
\end{enumerate}
Then $X^N$ converges to $X$ in distribution with limit $X$ in $C(\mathbb{R}_+, \mathbb{R}^d)$.
\end{theorem}

\begin{lemma} 
\label{lem:martignaleconvergence}
In $\mathbb{D}(\mathbb{R}_+, \mathcal{D}' ( \mathbb{T}^d \times [0,1] ) $, the sequence of martingale measures $\mathfrak{M}^N$ converges in distribution to a continuous Gaussian martingale limit $\mathfrak{M}$, such that for any $\phi \in \mathcal{D} (\mathbb{T}^d \times [0,1])$
\begin{equation} \label{eq:limitingcovariance}
    \big\langle \mathfrak{M}_t(\phi) \big\rangle = t \big\langle \mathcal{Q}, \phi \otimes \phi \big\rangle,
\end{equation}
where $\mathcal{Q}$ is given in \eqref{eq:limitingcovariation}.
\end{lemma}
\begin{proof}
The uniform upper and lower bounds $n_{\min}, n_{\max}$ on the population size allow one to check the first two assumptions of \cref{theo:convergenceofmartingales} for $M_t (\phi)$ in a similar way to \parencite[Lemma 3.10, fixed radius case]{forien}, which also inspired the strategy of the rest of this proof.
Let us recall the expression for the covariation of $\mathfrak{M}^N$ from \Cref{lem:mpfluctuations}:
\begin{equation*}
\begin{aligned}
\Blangle \mathfrak{M}^N (\phi) \Brangle_t
&= \frac{u^2}{ \delta_N^{2 d}} \int_0^t \int_{\mathbb{T}^d} \frac{1}{\int_{B(x,\delta_N R)} n_s^N (y) dy} \int_{B(x, \delta_N R) \times [0,1]} \\
&\hspace{1cm} \frac{1}{\overline{n_s^N} (x, \delta_N R) +1} \Bigg[ \int_{B(x,\delta_N R)} \Big[ [1 + \delta_N^2 r_x(\overline{n_s^N}(x, \delta_N R)) ] \overline{n_s^N} (x, \delta_N R) \phi (z, k_0) \\
&\hspace{4cm} - \int_{[0,1]} \phi (z,k) \rho_s^N (z, dk) \Big] dz \Bigg]^2 \rho_s^N (y, dk_0) dy  dx ds
\end{aligned}
\end{equation*}
We start by proving that the term involving $\delta_N^2 r(\overline{n_s^N})$ does not contribute in the limit. Bounding the term in question, we can see 
\begin{align*}
&\begin{aligned}
 &\Bigg\vert \frac{u^2}{ \delta_N^{2 d}} \int_{\mathbb{T}^d} \frac{1}{\int_{B(x,\delta_N R)} n_s^N (y) dy} \int_{B(x, \delta_N R) \times [0,1]} \frac{1}{\overline{n_s^N} +1} \Bigg[ \int_{B(x,\delta_N R)} \delta_N^2 r(\overline{n_s^N}) \overline{n_s^N} \phi (z, k_0) dz \Bigg]\\
&\hspace{1cm} \times \Bigg[ \int_{B(x,\delta_N R)} \Big[ [1 + \delta_N^2 r(\overline{n_s^N}) ] \overline{n_s^N} \phi (z, k_0) - \int_{[0,1]} \phi (z,k) \rho_s^N (z, dk) \Big] dz \Bigg] \rho_s^N (y, dk_0) dy dx \Bigg\vert
\end{aligned}\\
&\begin{aligned}
 &\leq \frac{u^2}{ \delta_N^{2 d}} \int_{\mathbb{T}^d} \frac{1}{\int_{B(x,\delta_N R)} n_s^N (y) dy} \int_{B(x, \delta_N R) \times [0,1]} \frac{1}{\overline{n_s^N} +1} \delta_N^2 \Vert r \Vert_\infty n_{\max} \int_{B(x,\delta_N R)} \sup_{k \in [0,1]} \big\vert \phi (z, k) \big\vert dz \\
&\hspace{1cm} \times (2 n_{\max} + \delta_N^2 \Vert r \Vert_\infty) \int_{B(x,\delta_N R)} \sup_{k \in [0,1]} \big\vert \phi (z, k) \big\vert dz \; \rho_s^N (y, dk_0) dy dx.
\end{aligned}
\end{align*}
The integrand no longer depends on $k_0$ anymore and the integral with respect to $\rho_s^N (y, dk_0) dy$ cancels with the average of the population size over the ball $B(x, \delta_N R)$. As before, let $V_{\delta_N R} (z_1, z_2)$ denote the volume of the intersection of two balls of radius $\delta_N R$ centered at $z_1$ and $z_2$. Then, using the identity \eqref{eq:identity_with_intersection}, we obtain the upper bound
\begin{align*}
&\begin{aligned}
&\delta_N^2 \Vert r \Vert_\infty n_{\max} (2 n_{\max} + \delta_N^2 \Vert r \Vert_\infty ) \frac{u^2}{ \delta_N^{2 d}} \int_{\mathbb{T}^d} \int_{B(x,\delta_N R)} \sup_{k \in [0,1]} \big\vert \phi (z, k) \big\vert dz \int_{B(x,\delta_N R)} \sup_{k \in [0,1]} \big\vert \phi (z, k) \big\vert dz \; dx 
\end{aligned}\\
&\begin{aligned}
&\leq \delta_N^2 \Vert r \Vert_\infty n_{\max} (2 n_{\max} + \Vert r \Vert_\infty ) \frac{u^2}{ \delta_N^{2 d}} \int_{(\mathbb{T}^d)^2} V_{\delta_N R} (z_1, z_2) \sup_{k \in [0,1]} \big\vert \phi (z_1, k) \big\vert \sup_{k \in [0,1]} \big\vert \phi (z_2, k) \big\vert dz_1 dz_2 
\end{aligned}\\
&\leq \delta_N^2 \Vert r \Vert_\infty n_{\max} (2 n_{\max} + \Vert r \Vert_\infty ) u^2 V_R^2 \Vert \phi \Vert_2^2.
\end{align*}
We can therefore focus our analysis on
\begin{equation} \label{eq:martingaleconvergenceremainingterms}
\begin{aligned}
&\frac{u^2}{ \delta_N^{2 d}} \int_0^t \int_{\mathbb{T}^d} \frac{1}{\int_{B(x,\delta_N R)} n_s^N(y) dy} \int_{B(x, \delta_N R) \times [0,1]} \\
&\hspace{1cm} \frac{1}{\overline{n_s^N} +1} \int_{(\mathbb{T}^d)^2} \mathds{1}_{B(x, \delta_N R)} (z_1) \mathds{1}_{B(x, \delta_N R)} (z_2) \Big( \overline{n_s^N} \phi (z_1, k_0) - \int_{[0,1]} \phi (z_1,k) \rho_s^N (z_1, dk) \Big)\\
& \hspace{2cm} \times \Big( \overline{n_s^N} \phi (z_2, k_0) - \int_{[0,1]} \phi (z_2,k) \rho_s^N (z_2, dk) \Big)dz_1 dz_2 \rho_s^N (y, dk_0) dy  dx ds.
\end{aligned}
\end{equation}
The product of the two brackets in \eqref{eq:martingaleconvergenceremainingterms} will result in four terms, which we will treat separately. Introducing the average
\begin{equation*}
\overline{\rho} (x, dk, \delta_N R) = \frac{1}{\int_{B(x,\delta_N R)} n(y) dy} \int_{B(x, \delta_N R)} \rho (y, dk) dy,
\end{equation*}
the first term can be written as
\begin{align*}
&\begin{aligned}
 &\frac{u^2}{ \delta_N^{2 d}} \int_{\mathbb{T}^d} \frac{1}{\int_{B(x,\delta_N R)} n_s^N (y) dy} \int_{B(x, \delta_N R) \times [0,1]} \\
&\hspace{1cm} \frac{1}{\overline{n_s^N} +1} \int_{(\mathbb{T}^d)^2} \mathds{1}_{B(x, \delta_N R)} (z_1) \mathds{1}_{B(x, \delta_N R)} (z_2) \overline{n_s^N} (x) \phi (z_1,k_0) \overline{n_s^N} (x) \phi (z_2,k_0) \\
&\hspace{5cm} dz_1 dz_2 \rho_s^N (y, dk_0) dy dx
\end{aligned}\\
&= u^2 V_{R}^2 \Bigg\langle \overline{\rho_s^N} (\cdot, \delta_N R ) , \frac{\overline{n_s^N}^2}{\overline{n_s^N} + 1} \overline{\phi} (\cdot, \delta_N R)^2 \Bigg\rangle.
\end{align*}
The second term is given by
\begin{equation} \label{eq:second_covariation_term}
\begin{aligned}
 &\frac{u^2}{ \delta_N^{2 d}} \int_{\mathbb{T}^d} \frac{1}{\overline{n_s^N} +1} \int_{(\mathbb{T}^d)^2} \mathds{1}_{B(x, \delta_N R)} (z_1) \mathds{1}_{B(x, \delta_N R)} (z_2) \\
&\hspace{1cm} \times \int_{[0,1]} \overline{n_s^N} (x) \phi (z_1,k_0) \overline{\rho_s^N} (x,d k_0, \delta_N R) \int_{[0,1]} \phi (z_2,k) \rho_s^N (z_2, dk) dz_1 dz_2  dx.
\end{aligned}
\end{equation}
Replacing $\phi (z_i, k)$ by $\phi (x,k)$ in the second term results in
\begin{equation} \label{eq:second_covariation_term_modified}
\begin{aligned}
& u^2 V_R^2 \int_{\mathbb{T}^d \times [0,1]^2} \frac{1}{\overline{n_s^N} (x) +1} \overline{n_s^N} (x)^2 \phi (x,k_0)  \phi (x,k) \overline{\rho_s^N} (x,d k_0, \delta_N R) \overline{\rho_s^N} (x, dk, \delta_N R)  dx\\
&= u^2 V_R^2 \Bigg\langle \overline{\rho_s^N} (\cdot , \delta_N R ) \otimes \overline{\rho_s^N} (\cdot, \delta_N R) ,\frac{\overline{n_s^N}^2}{\overline{n_s^N} +1} \phi \otimes \phi \Bigg\rangle.
\end{aligned}
\end{equation}
We can control the error made by the replacement as
\begin{equation*}
\Bigg\vert \int_{B(x, \delta_N R) \times [0,1]} \big(\phi (z_2, k) - \phi (x,k ) \big) \rho_s^N (z_2, dk) dz_2 \Bigg\vert \leq V_{\delta_N R}n_{\max} \cdot \delta_N R \max_{\vert \kappa \vert = 1} \Vert \partial_\kappa \phi \Vert_\infty.
\end{equation*}
To bound the distance between the original second term \eqref{eq:second_covariation_term} and the modified one \eqref{eq:second_covariation_term_modified}, we can add and subtract an intermediate term with only one of the replacements and apply the triangle inequality. This allows us to bound the distance by
\begin{equation*}
2 \frac{u^2}{\delta_N^{2d}} \cdot n_{\max} V_{\delta_N R} \Vert \phi \Vert_1 \cdot V_{\delta_N R}n_{\max} \delta_N R \max_{\vert \kappa \vert = 1} \Vert \partial_\kappa \phi \Vert_\infty =\delta_N 2 u^2 n_{\max}^2 R V_R^2 \max_{\vert \kappa \vert = 1} \Vert \partial_\kappa \phi \Vert_\infty \Vert \phi \Vert_1.
\end{equation*}
The third term is of the same form as the second, and the fourth term can be seen to equal
\begin{equation*}
u^2 V_R^2 \Bigg\langle \overline{\rho_s^N} (\cdot, \delta_N R) \otimes \overline{\rho_s^N} (\cdot, \delta_N R) ,\frac{\overline{n_s^N}^2}{\overline{n_s^N} +1} \phi \otimes \phi \Bigg\rangle.
\end{equation*}
We conclude
\begin{align*}
&\begin{aligned}
&\Bigg\vert \blangle \mathfrak{M}_t^N (\phi ) \brangle - u^2 V_R^2 \int_0^t \Bigg[ \Bigg\langle \overline{\rho_s^N} (\cdot, \delta_N R ) , \frac{\overline{n_s^N}^2}{\overline{n_s^N} + 1} \overline{\phi} (\cdot, \delta_N R)^2 \Bigg\rangle \\
& \hspace{1cm} - \Bigg\langle \overline{\rho_s^N} (\cdot , \delta_N R ) \otimes \overline{\rho_s^N} (\cdot, \delta_N R) ,\frac{\overline{n_s^N}^2}{\overline{n_s^N} +1} \phi \otimes \phi \Bigg\rangle \Bigg] ds \Bigg\vert\\
&\leq \delta_N 4 t u^2 n_{\max}^2 R V_{R}^2 \max_{\vert \kappa \vert = 1} \Vert \partial_\kappa \phi \Vert_\infty \Vert \phi \Vert_1.
\end{aligned}
\end{align*}
Observe that
\begin{align*}
&\Bigg\vert \frac{\overline{n_s^N}^2}{\overline{n_s^N} + 1} - \frac{\overline{m_s^N}^2}{\overline{m_s^N} + 1} \Bigg\vert \\
&\leq \Big\vert \overline{n_s^N}^2 (\overline{m_s^N} +1) - \overline{m_s^N}^2 ( \overline{n_s^N} +1 ) \pm \overline{n_s^N}^2 ( \overline{n_s^N} +1 ) \Big\vert\\
&= \Big\vert \overline{n_s^N}^2 (\overline{m_s^N} - \overline{n_s^N}) + (\overline{n_s^N} - \overline{m_s^N} ) ( \overline{n_s^N} + \overline{m_s^N} ) ( \overline{n_s^N} +1 ) \Big\vert\\
&\leq n_{\max}^2 \big\vert \overline{n_s^N} - \overline{m_s^N} \big\vert + \big\vert \overline{n_s^N} - \overline{m_s^N} \big\vert 2 n_{\max} (n_{\max} + 1)\\
&\leq 3 n_{\max} ( n_{\max} +1) \big\vert \overline{n_s^N} - \overline{m_s^N} \big\vert.
\end{align*}
Setting $C = 3 n_{\max} ( n_{\max} +1)$, we obtain the bound
\begin{align*}
& \mathbb{E} \Bigg\vert \Bigg\langle \overline{\rho_s^N} (\cdot, \delta_N R ) , \Bigg( \frac{\overline{n_s^N}^2}{\overline{n_s^N} + 1} - \frac{\overline{m_s^N}^2}{\overline{m_s^N} + 1} \Bigg) \overline{\phi} (\cdot, \delta_N R)^2 \Bigg\rangle \Bigg\vert \\
&\leq \frac{C}{\sqrt{N \delta_N^{2-d}}} \mathbb{E} \Bigg[ \int_{\mathbb{T}^d} \overline{F_s^N} (x, \delta_N R)^2 dx \Bigg]^{1/2} \mathbb{E} \Bigg[ \int_{\mathbb{T}^d} \Bigg( \int_{[0,1]} \big\vert \overline{\phi} (x, \delta_N R) \big\vert^2 \overline{\rho_s^N} (x, dk, \delta_N R) \Bigg)^2 dx \Bigg]^{1/2}\\
&\leq \Big( N \delta_N^2 \Big)^{- \frac{1}{2}} C \sqrt{\frac{C'}{V_R}} \Vert \phi \Vert_4^2,
\end{align*}
which converges to zero by assumption \eqref{eq:assumptionondelta}. The same bound also holds for 
\begin{equation*}
\mathbb{E} \Bigg\vert \Bigg\langle \overline{\rho_s^N} (\cdot , \delta_N R ) \otimes \overline{\rho_s^N} (\cdot, \delta_N R) , \Bigg( \frac{\overline{n_s^N}^2}{\overline{n_s^N} + 1} - \frac{\overline{m_s^N}^2}{\overline{m_s^N} + 1} \Bigg) \phi \otimes \phi \Bigg\rangle \Bigg\vert.
\end{equation*}
As $N \to \infty$, the convergence $\rho_s^N \to n_s \lambda$ in probability is implied by the central limit theorem. The same convergence follows for the averages $\overline{\rho_s^N} (\cdot, \delta_N R) $ by \Cref{prop:llnpop} for the averaged population size $\overline{n_s^N} (\cdot, \delta_N R)$. As we are evaluating against deterministic test functions, the complete expression
\begin{align*}
&u^2 V_R^2 \int_0^t \Bigg[ \Bigg\langle \overline{\rho_s^N} (\cdot, \delta_N R ) , \frac{\overline{m_s^N}^2}{\overline{m_s^N} + 1} \overline{\phi} (\cdot, \delta_N R)^2 \Bigg\rangle - \Bigg\langle \overline{\rho_s^N} (\cdot , \delta_N R ) \otimes \overline{\rho_s^N} (\cdot, \delta_N R) ,\frac{\overline{m_s^N}^2}{\overline{m_s^N} +1} \phi \otimes \phi \Bigg\rangle \Bigg] ds
\end{align*}
and thus $\langle \mathfrak{M}^N (\phi) \rangle_t$ converges to \eqref{eq:limitingcovariance}. 
We can now proceed as in \parencite[Section 4.5]{forien_central_2017}. By \Cref{theo:convergenceofmartingales}, the sequence $(\mathfrak{M}_t^N (\phi ))_{t \geq 0}$ converges in distribution to $(\mathfrak{M}_t(\phi))_{t \geq 0}$ in $\mathbb{D}(\mathbb{R}_+, \mathcal{D}' (\mathbb{T}^d \times [0,1]))$ and is in particular tight. Additionally, \Cref{theo:convergenceofmartingales} yields the convergence of finite-dimensional distributions of $\mathfrak{M}^N$, as $\langle \mathfrak{M}^N (\phi_i ) , \mathfrak{M}^N(\phi_j) \rangle_t$ is uniquely described by $\langle \mathfrak{M}^N (\phi_i + \phi_j ) \rangle_t$ and $\langle \mathfrak{M}^N (\phi_i - \phi_j ) \rangle_t$ using the polarization identity. An application of \Cref{theo:convergenceofstochasticprocesses} finishes this proof.
\end{proof}
\subsection{Conclusion of the proof of the central limit theorem}
In the previous subsections we have collected all the ingredients for a proof of \Cref{theo:clt}. Tightness of the sequences $((\langle Z_t^N , \phi \rangle)_{t \geq 0}, N \geq 1)$ has been proven in \Cref{subsec:tightness}. The convergence of finite-dimensional distributions follows from an application of \Cref{theo:convergencefinitedistributions}. Its prerequisites are the convergence properties of the integrands $g^N$ (\Cref{lem:semigroupsbounded}, \Cref{lem:semigroupconvergence}) and the convergence of the martingale measures \Cref{lem:convergencemn}. \Cref{theo:clt} now follows from an application of the following theorem.
\begin{theorem}[{\parencite[Theorem 6.15]{walsh_introduction_1986}}] \label{theo:convergenceofstochasticprocesses}
Suppose for a sequence $((X_t^N)_{t \geq 0}, N \geq 0)$ of distribution-valued stochastic processes in $\mathbb{D} ([0,T] , \mathcal{D}' (\mathbb{T}^d ))$
\begin{enumerate}
    \item 
    for every test function $\phi \in \mathcal{D} (\mathbb{T}^d)$ the sequence $((X^N_t (\phi))_{t \geq 0}, N \geq 0) $ is tight,
    \item
    for all test functions $\phi_1,\ldots,\phi_k \in \mathcal{D} (\mathbb{T}^d)$ and times $t_1,\ldots,t_k \in [0,T]$, $(X_{t_1}^N (\phi_1) ,\ldots, X_{t_k}^N (\phi_k) ) $ converges weakly in distribution.
\end{enumerate}
Then there exists limit $X = (X_t)_{t \geq 0}$ in distribution of $X^N$ in $\mathbb{D}([0,T], \mathcal{D}' (\mathbb{T}^d))$.
\end{theorem}

\section{Proof of the Wright-Mal\'ecot formula} \label{sec:proofofwmf}
We can describe the probability of identity by descent for two probability densities functions $\psi_1, \psi_2 \in \mathcal{D} (\mathbb{T}^d)$ by
\begin{equation*}
N \eta_N P_t^N (\phi, \psi) = \mathbb{E} \Bigg[ \frac{\langle Z_t^N \otimes Z_t^N , \mathds{1}_{\diagdown} \psi_1 \otimes \psi_2 \rangle}{\langle n_t^N, \psi_1 \rangle \langle n_t^N, \psi_2 \rangle} \Bigg].
\end{equation*}
A potential exchange of the limit $N \to \infty$ and the expectation was discussed in \Cref{remark:technicality}. Here, we will compute the limiting Wright-Mal\'ecot formula
\begin{equation} \label{eq:wmfquantity}
\mathbb{E} \Bigg[ \frac{\langle Z_t \otimes Z_t , \mathds{1}_{\diagdown} \psi_1 \otimes \psi_2 \rangle}{\langle n_t, \psi_1 \rangle \langle n_t, \psi_2 \rangle} \Bigg]
\end{equation}
given in \Cref{theo:varying_size_wmf}.

\begin{proof}[Proof of \Cref{theo:varying_size_wmf}]
For $\phi \in \mathcal{D} (\mathbb{T}^d \times[0,1] )$, the mild solution to the SPDE \eqref{eq:spde} can be written as
\begin{equation*}
\langle Z_t , \phi \rangle = \int_{[0,t] \times \mathbb{T}^d \times [0,1]} g_{s,t} (z,k) \mathfrak{M} (ds\; dz\; dk),
\end{equation*}
where $g$ is the solution to
\begin{equation*}
\left\{
\begin{aligned}
- \partial_s g_{s,t} &= u V_R \Big( \frac{R^2}{d+2} \Delta + r(n_s) \Big) g_{s,t} - \mu g_{s,t}\\
& \hspace{0.5cm} + \mu \mathcal{H} g_{s,t} + uV_R n_s r'(n_s) \mathcal{H} g_{s,t}, \\
 g_{t,t} &= \phi.
\end{aligned}
\right.
\end{equation*}
As $\mathfrak{M}$ is Gaussian, $\langle Z_t, \phi \rangle$ is also Gaussian with variance 
\begin{equation} \label{eq:expression_for_variance}
\mathbb{E} \big[ \blangle Z_t \otimes Z_t , \phi_1 \otimes \phi_2  \brangle \big] = \int_0^t \int_{(\mathbb{T}^d)^2} \int_{[0,1]^2} f_{s,t} (x_1, k_1, x_2, k_2) \mathcal{Q}_s (dx_1 dk_1 dx_2 dk_2 ) ds,
\end{equation}
for the covariation measure $\mathcal{Q}_s$ given in \eqref{eq:limitingcovariation}. Here, 
\[ (g_{s,t} \otimes g_{s,t} ) \Big( \phi_1 \otimes \phi_2 \Big) (x_1, k_1, x_2, k_2) := f_{s,t} (x_1, k_1, x_2, k_2) \] 
is defined as the solution $f_{s,t}$ to
\begin{equation*}
\left\{
\begin{aligned}
- \partial_s f_{s,t} &= u V_R \Big( \frac{R^2}{d+2} (\Delta \otimes Id + Id \otimes \Delta) + r_{x_1}(n_s (x_1)) + r_{x_2}(n_s(x_2)) \Big) f_{s,t}\\
& \hspace{0.5cm} + \mu \big( ( \mathcal{H} - Id) \otimes Id\big) f_{s,t} + \mu \big( Id \otimes ( \mathcal{H} - Id)\big) f_{s,t} \\
& \hspace{0.5cm} + uV_R n_s (x_1) r_{x_1}'(n_s (x_1)) \big(\mathcal{H} \otimes Id \big) f_{s,t} + uV_R n_s (x_2) r_{x_2}'(n_s (x_2)) \big(Id \otimes \mathcal{H} \big) f_{s,t}, \\
 f_{t,t} &= \phi_1 \otimes \phi_2.
\end{aligned}
\right.
\end{equation*}
The operators $\Delta \otimes Id$ and $\mathcal{H} \otimes Id$ act on the first space and type parameter - correspondingly $Id \otimes \Delta$ and $Id \otimes \mathcal{H}$ act on the second. Using linearity and approximating sequences, we can generalize the expression \eqref{eq:expression_for_variance} to functions $\Phi \in \mathcal{D} \big( (\mathbb{T}^d \times [0,1])^2 \big)$. We obtain
\begin{equation*}
\mathbb{E} \big[ \blangle Z_t \otimes Z_t , \Phi  \brangle \big] = \int_0^t \int_{(\mathbb{T}^d)^2} \int_{[0,1]^2} f_{s,t} (x_1, k_1, x_2, k_2) \mathcal{Q}_s (dx_1 dk_1 dx_2 dk_2 ) ds,
\end{equation*}
where $\mathcal{Q}$ is now defined as
\begin{equation} \label{eq:modified_covariation}
\blangle \mathcal{Q}_s , \Phi \brangle = u^2 V_R^2 \int_{\mathbb{T}^d} \frac{n_s (z)^2}{n_s(z) + 1} \Bigg( \int_{[0,1]} \Phi (z, k, z, k) dk - \int_{[0,1]^2} \Phi (z, k_1, z, k_2) dk_1 dk_2 \Bigg) dz.
\end{equation}

Let us now return to the quantity \eqref{eq:wmfquantity}. Again we define $\vartheta_i = \frac{n_t \psi_i}{\langle n_t , \psi_i \rangle}$, where $\psi_i$ only depends on space, and set $\Phi = \mathds{1}_{\diagdown} \psi_1 \otimes \psi_2$.
We get
\begin{align*}
\mathbb{E} \Bigg[ \frac{\langle Z_t \otimes Z_t , \mathds{1}_{\diagdown} \psi_1 \otimes \psi_2 \rangle}{\langle n_t, \psi_1 \rangle \langle n_t, \psi_2 \rangle} \Bigg] 
&= \mathbb{E} \Big[ \Big\langle Z_t \otimes Z_t , \mathds{1}_{\diagdown} \frac{\vartheta_1}{n_t} \otimes \frac{\vartheta_2}{n_t} \Big\rangle \Big] \\
&= \int_0^t \int_{(\mathbb{T}^d)^2} \int_{[0,1]^2} (g_{s,t} \otimes g_{s,t} ) \Big( \mathds{1}_{\diagdown} \frac{\vartheta_1}{n_t} \otimes \frac{\vartheta_2}{n_t} \Big) \mathcal{Q}_s (dx_1 dk_1 dx_2 dk_2 ) ds.
\end{align*}
We can show with a calculation, which we provide in \Cref{lem:termsdisappear} below, that this is equal to 
\begin{equation} \label{eq:wmfcalculation}
\int_0^t \int_{(\mathbb{T}^d)^2} \int_{[0,1]^2} \mathds{1}_{\diagdown} (k_1, k_2) P_{s,t} \left( \frac{\vartheta_1}{n_t} \right) (x_1) P_{s,t} \left(\frac{\vartheta_2}{n_t} \right) (x_2) Q_s (dx_1 dk_1 dx_2 dk_2 ) ds,
\end{equation}
where $P_{s,t} \vartheta$ solves, for $\vartheta \in \mathcal{D}(\mathbb{T}^d)$,
\begin{equation}
\left\{
\begin{aligned}
- \frac{d}{ds} P_{s,t} \vartheta &= \frac{u V_R R^2}{d + 2} \Delta P_{s,t} \vartheta + u V_R r(n_s) P_{s,t} \vartheta - \mu P_{s,t} \vartheta,\\
P_{t,t} \vartheta &= \vartheta .
\end{aligned}
\right.
\end{equation}
The indicator $\mathds{1}_{\diagdown}$ was defined as $\mathds{1}_{k_1 = k_2}$ and hence
\[ \int_{[0,1]} \mathds{1}_{\diagdown} (k, k) dk = 1, \hspace{1cm} \int_{[0,1]^2} \mathds{1}_{\diagdown} (k_1, k_2) dk_1 dk_2 = 0.\]
Therefore, under $\Phi = \mathds{1}_{\diagdown} \psi_1 \otimes \psi_2$ the negative term of \eqref{eq:modified_covariation} vanishes and \eqref{eq:wmfcalculation} is equal to
\begin{equation} \label{eq:varying_size_wmf_penultimate}
   u^2 V_R^2 \int_0^t \int_{\mathbb{T}^d} \frac{n_s(z)^2}{n_s(z) +1} P_{s,t} \left( \frac{\vartheta_1}{n_t} \right) (z) P_{s,t} \left(\frac{\vartheta_2}{n_t} \right) (z) dz ds.
\end{equation}
Similar to Doob's h-transform we define
\begin{equation*}
\tilde{P}_{s,t} \vartheta := n_s P_{s,t} \left( \frac{\vartheta}{n_t} \right).
\end{equation*}
We will see in \Cref{lem:doobstransform} below that $\tilde{P}_{s,t} \vartheta$ solves the equation
\begin{equation} \label{eq:wmfcalculationalmostG}
\left\{
\begin{aligned}
- \frac{d}{ds} \tilde{P}_{s,t} \vartheta &= A_s^* \tilde{P}_{s,t} \vartheta - \mu \tilde{P}_{s,t} \vartheta, \\
\tilde{P}_{t,t} &= \vartheta,
\end{aligned}
\right.
\end{equation}
where the operator $A_s$ is defined as
\begin{equation*}
A_s \vartheta = \frac{u V_R R^2}{d + 2} \Big( \Delta \vartheta + 2\frac{\nabla n_s}{n_s} \cdot \nabla \vartheta \Big).
\end{equation*}
Finally, the solution $\tilde{P}_{s,t} \vartheta$ to \eqref{eq:wmfcalculationalmostG} can be written as $e^{- \mu (t-s)} G_{s,t} \vartheta$, where $G_{s,t} \vartheta$ solves
\begin{equation*}
\left\{
\begin{aligned}
- \frac{d}{ds} G_{s,t} \vartheta &= A_s^* G_{s,t} \vartheta, \\
G_{t,t} \vartheta &= \vartheta. 
\end{aligned}
\right.
\end{equation*}
Substituting $e^{- \mu (t-s)} G_{s,t} \vartheta$ in \eqref{eq:varying_size_wmf_penultimate} results in \Cref{theo:varying_size_wmf}.
\end{proof}
We end this section with the two lemmas we used in the proof of \Cref{theo:varying_size_wmf}.
\begin{lemma} \label{lem:termsdisappear}
For $\vartheta_1 , \vartheta_2 \in \mathcal{D} (\mathbb{T}^d)$, the following equality holds
\begin{equation*}
(g_{s,t} \otimes g_{s,t} ) \big( \mathds{1}_{\diagdown} \vartheta_1 \otimes \vartheta_2 \big) = \mathds{1}_{\diagdown} P_{s,t} \big( \vartheta_1\big) P_{s,t} \big( \vartheta_2 \big),
\end{equation*}
where $P_{s,t} \vartheta$ is defined as the solution to
\begin{equation}
\left\{
\begin{aligned}
- \frac{d}{ds} P_{s,t} \vartheta &= \frac{u V_R R^2}{d + 2} \Delta P_{s,t} \vartheta + u V_R r(n_s) P_{s,t} \vartheta - \mu P_{s,t} \phi,\\
P_{t,t} \vartheta &= \vartheta .
\end{aligned}
\right.
\end{equation}
\end{lemma}
\begin{proof}
Recall that
\begin{equation*}
(g_{s,t} \otimes g_{s,t} ) \big( \mathds{1}_{\diagdown} \vartheta_1 \otimes \vartheta_2 \big) := f_{s,t} (x_1, k_1 , x_2 , k_2 ),
\end{equation*}
where $f_{s,t}$ satisfies the differential equation
\begin{equation} \label{eq:doublesemigroup}
\left\{
\begin{aligned}
- \frac{d}{ds} f_{s,t} &= u V_R \Big( \frac{R^2}{d+2} (\Delta \otimes Id + Id \otimes \Delta) + r_{x_1}(n_s (x_1)) + r_{x_2}(n_s(x_2)) \Big) f_{s,t} \\
& \hspace{0.5cm} + \mu \big( ( \mathcal{H} - Id) \otimes Id\big) f_{s,t} + \mu \big( Id \otimes ( \mathcal{H} - Id)\big) f_{s,t} \\
& \hspace{0.5cm} + uV_R n_s (x_1) r_{x_1}'(n_s (x_1)) \big(\mathcal{H} \otimes Id \big) f_{s,t} + uV_R n_s (x_2) r_{x_2}'(n_s (x_2)) \big(Id \otimes \mathcal{H} \big) f_{s,t}, \\
f_{t,t} &= \mathds{1}_{\diagdown} \vartheta_1 \otimes \vartheta_2.
\end{aligned}
\right.
\end{equation}
We aim to show
\begin{equation*}
    ( \mathcal{H} \otimes Id) f_{s,t} = (Id \otimes \mathcal{H} ) f_{s,t} = 0.
\end{equation*}
Abbreviating $( \mathcal{H} \otimes Id) f_{s,t}$ by $\tilde{f}_{s,t}$, its derivative and initial condition are
\begin{equation*}
\left\{
\begin{aligned}
- \frac{d}{ds} \tilde{f}_{s,t} &= u V_R \Big( \frac{R^2}{d+2} (\Delta \otimes Id + Id \otimes \Delta) + r_{x_1} (n_s (x_1)) + r_{x_2} (n_s(x_2)) \Big) \tilde{f} (s,t) \\
& \hspace{0.5cm} + \mu \big( Id \otimes ( \mathcal{H} - Id)\big) \tilde{f}_{s,t} \\
& \hspace{0.5cm} + uV_R n_s (x_1) r_{x_1}'(n_s (x_1)) \tilde{f}_{s,t} + uV_R n_s (x_2) r_{x_2}'(n_s (x_2)) \big(Id \otimes \mathcal{H} \big) \tilde{f}_{s,t}, \\
\tilde{f}_{t,t} &= (\mathcal{H} \otimes Id) \big(\mathds{1}_{\diagdown} \vartheta_1 \otimes \vartheta_2 \big) = 0.
\end{aligned}
\right.
\end{equation*}
For the initial condition we used $\mathcal{H} (\mathds{1}_{\diagdown}) = \int_0^1 \mathds{1}(k_1, k_2) dk_1 = 0$. Indeed, as the initial condition is zero, we obtain that $\tilde{f}_{s,t} = 0$. Returning to \eqref{eq:doublesemigroup}, this shows
\begin{equation*}
\left\{
\begin{aligned}
- \frac{d}{ds} f_{s,t} &= u V_R \Big( \frac{R^2}{d+2} (\Delta \otimes Id + Id \otimes \Delta) + r_{x_1} (n_s (x_1)) + r_{x_2} (n_s(x_2)) \Big) f_{s,t} - 2 \mu f_{s,t} ,\\
f_{t,t} &= \mathds{1}_{\diagdown} \vartheta_1 \otimes \vartheta_2.
\end{aligned}
\right.
\end{equation*}
As the operator only depends on the spatial components, we can now write
\begin{align*}
    (g_{s,t} \otimes g_{s,t} ) \big( \mathds{1}_{\diagdown} \vartheta_1 \otimes \vartheta_2 \big) &=  (P_{s,t} \otimes P_{s,t} ) \big( \mathds{1}_{\diagdown} \vartheta_1 \otimes \vartheta_2 \big)\\
    &= \mathds{1}_{\diagdown}  (P_{s,t} \otimes P_{s,t} ) \big( \vartheta_1 \otimes \vartheta_2 \big)\\
    &= \mathds{1}_{\diagdown} P_{s,t} \big( \vartheta_1 \Big) P_{s,t} \Big( \vartheta_2 \big). \qedhere
\end{align*}
\end{proof}
\begin{lemma} \label{lem:doobstransform}
For any $\vartheta \in \mathcal{D}(\mathbb{T}^d)$, define
\begin{equation*}
\tilde{P}_{s,t} \vartheta := n_s P_{s,t} \left( \frac{\vartheta}{n_t} \right),
\end{equation*}
where $(n_t)_{t \geq 0}$ solves the differential equation in \Cref{prop:llnpop}. Then $\tilde{P}_{s,t} \vartheta$ solves the differential equation
\begin{equation*}
\left\{
\begin{aligned}
- \frac{d}{ds} \tilde{P}_{s,t} \vartheta &= A_s^* \tilde{P}_{s,t} \vartheta - \mu \tilde{P}_{s,t} \vartheta, \\
\tilde{P}_{t,t} &= \vartheta,
\end{aligned}
\right.
\end{equation*}
where the operator $A_s$ is defined as
\begin{equation*}
A_s \phi = \frac{u V_R R^2}{d + 2} \Big( \Delta \vartheta + 2\frac{\nabla n_s}{n_s} \cdot \nabla \vartheta \Big).
\end{equation*}
\end{lemma}
\begin{proof}
Integrating by parts, we can see that
\begin{align*}
\int \frac{\partial_x n }{n} \partial_x \vartheta \cdot \phi &= - \int \vartheta \cdot \partial_x \Big( \frac{\partial_x n }{n } \cdot \phi \Big) \\
&= - \int \vartheta \Big( \frac{\partial_{xx} n}{n} - \frac{(\partial_x n)^2}{n^2} \Big) \phi - \int \vartheta \frac{\partial_x n}{n} \partial_x \phi.
\end{align*}
Thus, the adjoint of $A$ will be of form
\begin{equation*}
A^* \vartheta = \frac{u V_R R^2}{d + 2} \Big( \Delta \vartheta - 2\frac{\Delta n}{n} \vartheta + 2\frac{\Vert \nabla n \Vert^2}{n^2} \vartheta - 2\frac{\nabla n}{n} \cdot \nabla \vartheta \Big).
\end{equation*}
Recall that $\partial_t n_t = \frac{u V_R R^2}{d + 2} \Delta n_t + u V_R r (n_t ) n_t$ and 
\begin{align*}
\frac{d}{ds} P_{s,t} \vartheta &= - \Big( \frac{u V_R R^2}{d + 2} \Delta P_{s,t} \vartheta + u V_R r(n_s) P_{s,t} \vartheta - \mu P_{s,t} \vartheta \Big).
\end{align*}
Then we can calculate
\begin{align*}
\frac{d}{ds} \tilde{P}_{s,t} \vartheta 
&= \Big( \frac{d}{ds} n_s \Big) P_{s,t} \left( \frac{\vartheta}{n_t} \right) + n_s \frac{d}{ds} P_{s,t} \left( \frac{\vartheta}{n_t} \right)\\
&=  \Big( \frac{u V_R R^2}{d + 2} \Delta n_s + u V_R r(n_s) n_s \Big) \frac{1}{n_s} \tilde{P}_{s,t} \vartheta \\
& \hspace{1cm} - n_s \Big( \frac{u V_R R^2}{d + 2} \Delta + u V_R r(n_s) - \mu \Big) \left( \frac{1}{n_s} \tilde{P}_{s,t} \vartheta \right)\\
&= \frac{u V_R R^2}{d + 2} \Bigg[ \frac{\Delta n_s}{n_s} \cdot \tilde{P}_{s,t} \vartheta - n_s \Delta \left( \frac{1}{n_s} \tilde{P}_{s,t} \vartheta \right) \Bigg] + \mu \tilde{P}_{s,t} \vartheta.
\end{align*}
To calculate $\Delta \left( \frac{1}{n_s} \tilde{P}_{s,t} \vartheta \right)$ note that
\begin{align*}
\partial_x \left( \frac{\phi}{n_s} \right) &= \frac{\partial_x \phi}{n_s} - \frac{\partial_x n_s}{n_s^2} \cdot \phi, \\
\partial_{xx} \left( \frac{\phi}{n_s} \right) &= \frac{\partial_{xx} \phi}{n_s} - \frac{\partial_x n_s}{n_s^2} \partial_x \phi - \frac{\partial_{xx} n_s}{n_s^2} \phi - \frac{\partial_x n_s}{n_s^2} \partial_x \phi + 2 \frac{(\partial_x n_s )^2}{n_s^3} \phi\\
&= \frac{\partial_{xx} \phi}{n_s} - 2 \frac{\partial_x n_s}{n_s^2} \partial_x \phi - \frac{\partial_{xx} n_s}{n_s^2} \phi + 2 \frac{(\partial_x n_s)^2}{n_s^3} \phi.
\end{align*}
This implies first
\begin{align*}
n_s \partial_{xx} \Big( \frac{1}{n_s} \tilde{P}_{s,t} \vartheta \Big) 
&= \partial_{xx} \tilde{P}_{s,t} \vartheta - 2 \frac{\partial_x n_s}{n_s} \partial_x \tilde{P}_{s,t} \vartheta - \frac{\partial_{xx} n_s}{n_s} \tilde{P}_{s,t} \vartheta + 2 \frac{(\partial_x n_s )^2}{n_s^2} \tilde{P}_{s,t} \vartheta
\end{align*}
and in turn
\begin{equation*}
n_s \Delta \Big( \frac{1}{n_s} \tilde{P}_{s,t} \vartheta \Big) = \Delta \tilde{P}_{s,t} \vartheta - 2 \frac{\nabla n_s}{n_s} \nabla \tilde{P}_{s,t} \vartheta - \frac{\Delta n_s}{n_s} \tilde{P}_{s,t} \vartheta + 2 \frac{\Vert \nabla n \Vert^2}{n_s^2} \tilde{P}_{s,t} \vartheta.
\end{equation*}
In total, we can conclude
\begin{align*}
\frac{d}{ds} \tilde{P}_{s,t} \vartheta &= \frac{u V_R R^2}{d + 2} \Bigg[ \frac{\Delta n_s}{n_s} \cdot \tilde{P}_{s,t} \vartheta - n_s \Delta \left( \frac{1}{n_s} \tilde{P}_{s,t} \vartheta \right)  \vartheta \Bigg]  + \mu \tilde{P}_{s,t} \vartheta\\
&= \frac{u V_R R^2}{d + 2} \Bigg[ \frac{\Delta n_s}{n_s} \cdot \tilde{P}_{s,t} \vartheta - \Delta \tilde{P}_{s,t} \vartheta + 2 \frac{\nabla n_s}{n_s} \nabla \tilde{P}_{s,t} \vartheta + \frac{\Delta n_s}{n_s} \tilde{P}_{s,t} \vartheta - 2 \frac{\Vert \nabla n \Vert^2}{n_s^2} \tilde{P}_{s,t} \vartheta \Bigg] + \mu \tilde{P}_{s,t} \vartheta\\
&= \frac{u V_R R^2}{d + 2} \Bigg[ - \Delta \tilde{P}_{s,t} \vartheta + 2 \frac{\nabla n_s}{n_s} \nabla \tilde{P}_{s,t} \vartheta + 2 \frac{\Delta n_s}{n_s} \tilde{P}_{s,t} \vartheta - 2 \frac{\Vert \nabla n \Vert^2}{n_s^2} \tilde{P}_{s,t} \vartheta \Bigg] + \mu \tilde{P}_{s,t} \vartheta\\
&= - (A^* - \mu ) \tilde{P}_{s,t} \vartheta. \qedhere
\end{align*}
\end{proof}

\printbibliography

\appendix
\section{Details on the simulation} \label{appendix:simulation}
This section details the principles underlying \Cref{fig:1driftcomparison}.
\subsection{Simulation of the process}
We consider a discretized version of the process defined by \Cref{def:mSLFV} on $\{ 0 ,1,2,..., 99, 100 \}$. We represent the uniform proportion of the population by one type, while allowing for up to $2000$ other types in the population. If the uniform proportion type is chosen as the parental type in an event, we choose instead the next still empty type or replace the type with the smallest total mass in the population. We start the process with a uniform mass of three everywhere and let it evolve for a time of $125$ under the growth function \eqref{eq:growthfunctionexample}. Additionally, we choose
\begin{align*}
\mu = 0.0001, \hspace{1cm} u = 0.04, \hspace{1cm} R = 4.
\end{align*}

\subsection{Analytical prediction}
For the analytical part we consider ancestral lineages backwards in time. Again, we consider the discretized space $\{ 0, ..., 100\}$. We consider the population at stationarity solving 
\begin{equation*}
\frac{u V_R R^2}{3} \Big( n (x + 1) + n (x - 1) - 2 n (x) \Big) + r_x (n (x)) n(x) = 0 \enspace \forall x \in \{1, ..., 99\},
\end{equation*}
with boundary condition $n(0) = 8 = n(100)$ and precisely the same growth function $r_x$ as for the process simulation.
To describe the movement of ancestral lineages, we consider the continuous time random walk with transition rates
\[
Q(i,j) = \begin{dcases}
\frac{n(j)}{n (i)} & j = i \pm 1,\\
- \frac{n (i + 1)}{n(i)} - \frac{n(i-1)}{n (i)} & j = i, \\
0 & \text{else. }
\end{dcases}
\]
This corresponds to a discretization of the lineages evolving according to the operator \eqref{eq:generatorwmflineage}.
We can calculate the transition probabilities of the ancestral lineages as the matrix exponential
\begin{equation*}
\mathfrak{P} = e^{Q}.
\end{equation*}
We can then consider the sequence $K_0 = diag \{1,..., 1 \}$, the $100 \times 100$ matrix with ones only on the diagonal, and $K_t = K_{t-1} \mathfrak{P}$. The $(i,j)$ entry of $K_t$ will encapsulate the probability to transition from $i$ to $j$ over the period $t$. We can therefore describe the probability of identity by descent for two positions $l_1, l_2$ by
\[
\Theta(l_1, l_2) = \sum_{t \geq 1}^{t_{\max}} \sum_{0 \leq z \leq 100} e^{-2 \mu i} \frac{1}{n(z) + 1} K_i(l_1, z) K_i (l_2, z).
\]
Defining the diagonal matrix $D(n) = \{ (n(0) + 1)^{-1} ,..., (n(100) + 1)^{-1} \}$, we can calculate the probability of identity between all possible combinations of positions as the matrix
\[
\Theta = \sum_{t \geq 1}^{t_{\max}} e^{-2 \mu t} K_t D(n) (K_t)^T,
\]
where $(K_t)^T$ is the transposition of $K_t$. For \Cref{fig:1driftcomparison} we chose $t_{\max} = 28$ as due to the finite range the ancestral lineages start to equilibriate, leading to a much shallower decay as on would expect on a large domain.

\subsection{Alignment of simulation and prediction}
As mentioned in the introduction, The analytical result $A$ corresponds to the limiting Wright-Mal\'ecot formula (in particular not necessarily being a probability), whereas we obtain probabilities of identity $\mathcal{P}$ from the simulation of the process. The outcome of simulation $\mathcal{P}^N (x)$ and analytical prediction $\Theta(\mu, u, x)$ for locations at distance $x$ can be connected through
\begin{equation*}
N \delta_N \mathcal{P}^N (x) \simeq \Theta \Big( \frac{N}{\delta_N^2} \mu_N , N u_N, \delta_N x \Big), 
\end{equation*}
corresponding to our scaling of the probability of identity, space and the parameters. As we do not consider a sequence of processes, in our case we can set $\mu_N = \mu = 0.0001$ and $u_N = u = 0.04$, the same parameters as we used for the simulation of the process. To avoid complications about the spatial scaling, we set simply set $\delta_N = 1$. Choosing $N = 2.8$ led to the observed alignment of simulation and prediction.

\section{Properties of \texorpdfstring{$g^N$}{}}\label{appendix:propertiesofthefunctions}
\subsection{Boundedness}

\begin{proof}[Proof of \Cref{lem:semigroupsbounded}]
Similar to \cite[Lemma 4.3]{forien_central_2017}, we can represent the solution as
\begin{equation*}
\begin{aligned}
g_{s,t}^N (x,k) &= G_{t-s}^{N} * \phi (x,k) + \int_s^t G_{u-s}^{N} * \Big( \mathcal{R}_{m_u^N} (g_{u,t}^N) - \mu g_{u,t}^N\\
&\hspace{1cm} + \mu \mathcal{H} g_{u,t}^N + uV_R \overline{\overline{m}_u^N r'(\overline{m}_u^N) \overline{\mathcal{H} g_{u,t}^N}} \Big) (x,k) du,
\end{aligned}
\end{equation*}
again with $G^N$ being the semigroup associated to $\mathcal{L}^N$. The inequality $\Vert G_{t-s}^N * \phi \Vert_p \leq \Vert \phi\Vert_p$ for functions $\phi \in \mathcal{D}(\mathbb{T}^d)$ extends to our functions $\phi \in \mathcal{D} (\mathbb{T}^d \times [0,1])$ with a type component and the extended $p$-norms of \eqref{eq:pnormdefinition}
\begin{align*}
\Vert G_{t-s}^N * \phi \Vert_p 
&= \Bigg( \int_{\mathbb{T}^d} \sup_{k \in [0,1]} \Bigg\vert \int_{\mathbb{T}^d} G_{t-s}^N (z-y) \phi (y, k) dy \Bigg\vert^p dz \Bigg)^{1/p} \\
&\leq \Bigg( \int_{\mathbb{T}^d} \sup_{k \in [0,1]} \int_{\mathbb{T}^d} G_{t-s}^N (z-y) \big\vert \phi(y,k) \big\vert^p dy dz \Bigg)^{1/p}\\
&\leq \Bigg( \int_{\mathbb{T}^d} \int_{\mathbb{T}^d} G_{t-s}^N (z-y) \sup_{k \in [0,1]} \big\vert \phi(y,k) \big\vert^p dy dz \Bigg)^{1/p}\\
&= \Bigg( \int_{\mathbb{T}^d} \sup_{k \in [0,1]} \big\vert \phi(y,k) \big\vert^p dy \Bigg)^{1/p} = \Vert \phi \Vert_p.
\end{align*}
Applying $(a+b)^2 \leq 2a^2 + 2 b^2$ and Jensen's inequality again
\begin{align*}
\big\Vert g_{s,t}^N \big\Vert_p^p &\leq (2 (t-s))^{p-1} \int_s^t \Big\Vert u V_R \overline{r(\overline{m}_t^N) \overline{g_{u,t}^N}}  - \mu g_{u,t}^N \\
& \hspace{1cm} + \mu \mathcal{H} g_{u,t}^N + uV_R \overline{\overline{m}_u^N r'(\overline{m}_u^N) \overline{\mathcal{H} g_{u,t}^N}} \Big\Vert_p^p du+  2^{p - 1} \big\Vert \phi \big\Vert_p^p\\
&\leq (2 (t-s))^{p-1} \int_s^t \Big( u V_R r_{\max} + 2 \mu + u V_R n_{\max} \sup_{z \in \mathbb{T}^d} \Vert r_z' \Vert_\infty \Big) \big\Vert g_{s,t}^N \big\Vert_p^p du \\
& \hspace{1cm} + 2^{p - 1} \big\Vert \phi \big\Vert_p^p.
\end{align*}
We also used that $\Vert \overline{\phi} \Vert_p \leq \Vert \phi \Vert_p$ and that, as $[0,1]$ is finite, \[\big\vert \mathcal{H} g_{u,t}^N \big\vert \leq \sup_{k \in [0,1]} \big\vert g_{u,t}^N \big\vert.\] 
Gr\"onwall's inequality allows us to conclude the statement and an inductive argument shows the generalization to the derivatives.
\end{proof}

\subsection{Convergence}

\begin{restatable}{lemma}{lemconvergencegenerators} \label{lem:convergenceofthegenerators}
Let $m^N$ be the sequence of deterministic approximate population sizes of \eqref{eq:middleterm} converging to the solution $n$ of the reaction-diffusion equation \eqref{eq:reactiondiffusion} by \Cref{lem:convergencemn}.
Let $\phi: \mathbb{T}^d \rightarrow \mathbb{R}$ be a four times continuously-differentiable function and assume $\max_{\vert \kappa \vert \leq 4} \Vert \partial_\kappa \phi \Vert_p < \infty$ for $1 \leq p \leq \infty$. 
Then there exists constants $C_1, C_2 > 0$ such that
\begin{equation} \label{eq:boundongenerator}
\big\Vert \big( \mathcal{L}^N + \mathcal{R}^N_{m_t^N} \big)  \phi \big\Vert_p \leq C_1 \max_{\vert \kappa \vert \leq 2} \Vert \partial_\kappa \phi \Vert_p \; 
\end{equation}
and
\begin{equation} \label{eq:convergenceofgenerators}
\begin{aligned}
    &\bigg\Vert \big( \mathcal{L}^N + \mathcal{R}^N_{m_t^N} \big) \phi - u V_R \Big( \frac{R^2}{d+2} \Delta + r(n_s) \Big) \phi \bigg\Vert_p \leq \delta_N C_2 .
\end{aligned}
\end{equation}
\end{restatable}

\begin{proof} 
From \cite[Proposition A.1]{forien_central_2017} we know
\begin{equation*}
\Bigg\Vert \mathcal{L}^N \phi - u V_R \frac{R^2}{d + 2} \Delta \phi \Bigg\Vert_p \leq u V_R R^4 \frac{d^3}{3} (\delta_N)^2 \max_{\vert \kappa \vert \leq 4} \Vert \partial_{\kappa} \phi \Vert_p
\end{equation*}
and 
\begin{equation} \label{eq:boundphiaverageminusphi}
\Big\Vert \overline{\phi} (\cdot , \delta_N R) - \phi \Big\Vert_p \leq \frac{d}{2} (\delta_N R)^2 \max_{\vert \kappa \vert = 2} \Vert \partial_\kappa \phi \Vert_p.
\end{equation}
The growth operator $\mathcal{R}^N_{m_t^N}$ was defined as
\begin{equation*}
\begin{aligned}
\mathcal{R}_{m_t^N}^N \phi (z) &= u V_R \frac{1}{V_{\delta_N R}} \int_{B(z, \delta_N R)} \frac{1}{V_{\delta_N R}} \int_{B(x,\delta_N R)} r_x(\overline{m_s^N}(x, \delta_NR)) \phi (y) dy dx\\
&= u V_R \overline{r (\overline{m_t^N}) \overline{\phi}} (\cdot, \delta_N R).
\end{aligned}
\end{equation*}
Adding and substracting terms, we obtain
\begin{align*}
&\Big\Vert \overline{r (\overline{m_t^N}) \overline{\phi}} (\cdot, \delta_N R) - r(n_t) \phi \Big\Vert_p\\
& \leq \Big\Vert \overline{r (\overline{m_t^N}) \overline{\phi}} (\cdot, \delta_N R) - \overline{r (\overline{n_t}) \overline{\phi}} (\cdot, \delta_N R) \Big\Vert_p + \Big\Vert \overline{r (\overline{n_t}) \overline{\phi}} (\cdot, \delta_N R) - r(n_t) \phi \Big\Vert_p\\
& \leq \Big\Vert \big[ r (\overline{m_t^N}) - r ( \overline{n_t}) \big] \overline{\phi} (\cdot, \delta_N R) \Big\Vert_p + \Big\Vert \overline{r (\overline{n_t}) \overline{\phi}} (\cdot, \delta_N R) - r (\overline{n_t}) \overline{\phi}\Big\Vert_p \\
& \hspace{1cm} + \Big\Vert r (\overline{n_t}) \overline{\phi} - r (\overline{n_t}) \phi \Big\Vert_p + \Big\Vert r (\overline{n_t} ) \phi - r(n_t) \phi \Big\Vert_p\\
& \leq \sup_{z \in \mathbb{T}^d} \Vert r_z' \Vert_\infty \big\Vert \vert \overline{m_t^N} - \overline{n_t^N} \vert \overline{\phi} \big\Vert_p + \frac{d}{2} (\delta_N R)^2 \max_{\vert \kappa \vert = 2} \Big\Vert \partial_\kappa \big( r (\overline{n_t}) \overline{\phi} \big) \Big\Vert_p \\
& \hspace{1cm} + r_{\max} \frac{d}{2} (\delta_N R)^2 \max_{\vert \kappa \vert = 2} \Vert \partial_\kappa \phi \Vert_p + \sup_{z \in \mathbb{T}^d} \Vert r_z' \Vert_\infty \Big\Vert \big\vert \overline{n_t} - n_t \big\vert \phi \Big\Vert_p\\
& \leq \sup_{z \in \mathbb{T}^d} \Vert r_z' \Vert_\infty \delta_N C \Vert \phi \Vert_p + \frac{d}{2} (\delta_N R)^2 \max_{\vert \kappa \vert = 2} \Big\Vert \partial_\kappa \big( r (\overline{n_t}) \overline{\phi} \big) \Big\Vert_p \\
& \hspace{1cm} + r_{\max} \frac{d}{2} (\delta_N R)^2 \max_{\vert \kappa \vert = 2} \Vert \partial_\kappa \phi \Vert_p + \sup_{z \in \mathbb{T}^d} \Vert r_z' \Vert_\infty \; \delta_N R \max_{\vert \kappa \vert = 1} \Vert \partial_\kappa n_t \Vert_\infty \Vert \phi \Vert_p,
\end{align*}
where we used \Cref{lem:convergencemn}.
\end{proof}
\begin{lemma} \label{lem:convergenceofthegenerators2}
For any twice-continuously differentiable function $\phi : \mathbb{T}^d \rightarrow \mathbb{R}$ there exists a constant $C > 0$ such that
\begin{equation*}
\Big\Vert uV_R \overline{\overline{m}_u^N r'(\overline{m}_u^N) \overline{\phi}} - uV_R n_u r'(n_u) \phi \Big\Vert_p \leq \delta_N C.
\end{equation*}
\end{lemma}
\begin{proof}
We can write
\begin{align*}
&\Big\Vert \overline{\overline{m}_u^N r'(\overline{m}_u^N) \overline{\phi}} - \overline{n_u r'(\overline{m}_u^N) \overline{\phi}} \Big\Vert_p + \Big\Vert \overline{n_u r'(\overline{m}_u^N) \overline{\phi}} - \overline{n_u r'(n_u) \overline{\phi}} \Big\Vert_p \\
& \hspace{1cm} + \Big\Vert \overline{n_u r'(n_u) \overline{\phi}} - \overline{n_u r'(n_u) \phi} \Big\Vert_p + \Big\Vert \overline{n_u r'(n_u) \phi}  - n_u r'(n_u) \phi \Big\Vert_p \\
&\leq \Big\Vert \overline{( \overline{m}_u^N - n_u ) r'(\overline{m}_u^N) \overline{\phi}} \Big\Vert_p + \Big\Vert \overline{n_u (r'(\overline{m}_u^N) - r'(n_u)) \overline{\phi}} \Big\Vert_p\\
&\hspace{1cm} + \Big\Vert \overline{n_u r'(n_u) (\overline{\phi} - \phi)} \Big\Vert_p + \Big\Vert \overline{n_u r'(n_u) \phi}  - n_u r'(n_u) \phi \Big\Vert_p\\
&\leq \delta_N C \sup_{z \in \mathbb{T}^d} \Vert r_z' \Vert_\infty \Vert \phi \Vert_p  + n_{\max} \sup_{z \in \mathbb{T}^d} \Vert r_z'' \Vert_\infty \delta_N C \Vert \phi \Vert_p\\
&\hspace{1cm} + n_{\max} \sup_{z \in \mathbb{T}^d} \Vert r_z' \Vert_\infty \frac{d}{2} (\delta_N R)^2 \max_{\vert \kappa \vert = 2} \Vert \partial_\kappa \phi \Vert_p + \frac{d}{2} (\delta_N R)^2 \max_{\vert \kappa \vert = 2} \big\Vert \partial_\kappa (n_u r' (n_u) \phi ) \big\Vert_p.
\end{align*}
We used \Cref{lem:convergencemn} for the first and second term and \eqref{eq:boundphiaverageminusphi} for the third and fourth term.
\end{proof}

\begin{proof}[Proof of \Cref{lem:semigroupconvergence}]
To show the convergence, we use the representation
\begin{equation*}
\begin{aligned}
g_{s,t} &= G_{t-s} * \phi (x) + \int_s^t G_{u-s} * \Big( u V_R r(n_u) (g_{u,t}) - \mu g_{u,t}\\
&\hspace{1cm} + \mu \mathcal{H}g_{u,t} + uV_R n_u r'(n_u) \mathcal{H}g_{u,t} \Big) du,
\end{aligned}
\end{equation*}
\begin{equation*}
\begin{aligned}
g_{s,t}^N &= G_{t-s} * \phi (x) + \int_s^t G_{u-s} * \Big( \mathcal{L}^N (g_{u,t}^N) - u V_R \frac{R^2}{d+2} \Delta g_{u,t}^N +  \mathcal{R}_{m_s^N} (g_{u,t}^N)\\
&\hspace{1cm} - \mu g_{u,t}^N + \mu \mathcal{H} g_{u,t}^N + uV_R \overline{\overline{m}_u^N r'(\overline{m}_u^N) \overline{\mathcal{H} g_{u,t}^N}} \Big) du.
\end{aligned}
\end{equation*}
Here, $G_t$ denotes the Gaussian kernel associated to the Laplacian $uV_R \frac{R^2}{d+2} \Delta$. With the help of \Cref{lem:convergenceofthegenerators} and \Cref{lem:convergenceofthegenerators2}, we obtain
\begin{align*}
&\Big\Vert g_{s,t}^N - g_{s,t} \Big\Vert_p \\
&\begin{aligned}
&\leq \int_s^t \Bigg\Vert  \mathcal{L}^N (g_{u,t}^N) - u V_R \frac{R^2}{d+2} \Delta g_{u,t}^N +  \mathcal{R}_{m_s^N} (g_{u,t}^N) - u V_R r(n_u) (g_{u,t})\\
&\hspace{1cm} - \mu \big( g_{u,t}^N - g_{u,t} \big) + \mu \big( \mathcal{H} g_{u,t}^N - \mathcal{H}g_{u,t} \big) \\
& \hspace{1cm} + uV_R \overline{\overline{m}_u^N r'(\overline{m}_u^N) \overline{\mathcal{H} g_{u,t}^N}} - uV_R n_u r'(n_u) \mathcal{H}g_{u,t} \Bigg\Vert_p du
\end{aligned}\\
&\begin{aligned}
&\leq (t - s) \delta_N C + \int_s^t \Big\Vert u V_R r(n_u) (g_{u,t}^N) - u V_R r(n_u) (g_{u,t}) - \mu \big( g_{u,t}^N - g_{u,t} \big) \\
& \hspace{1cm} + \mu \big( \mathcal{H} g_{u,t}^N - \mathcal{H}g_{u,t} \big) + uV_R n_u r'(n_u) \mathcal{H}g_{u,t}^N - uV_R n_u r'(n_u) \mathcal{H}g_{u,t} \Big\Vert_p du
\end{aligned}\\
&\begin{aligned}
&\leq (t - s) \delta_N C + \big( u V_R (\sup_{z \in \mathbb{T}^d} \Vert r \Vert_\infty + n_{max} \sup_{z \in \mathbb{T}^d} \Vert r_z' \Vert_\infty) + 2 \mu \big) \int_s^t  \Big\Vert  g_{u,t}^N - g_{u,t} \Big\Vert_p du .
\end{aligned}
\end{align*}
Again Gr\"onwall's inequality yields the desired result.
\end{proof}

\subsection{Continuity in the second time component}

\begin{proof}[Proof of \Cref{lem:semigroupcontinuity}]
This proof is similar to \cite[Lemma 4.9]{forien_central_2017}. As in \Cref{lem:modulusofcontinuity}, we extend $g^N$ to $s,t \in [0,T]$ by setting
\begin{equation*}
g_{s,t}^N (x,k) = g_{s \land t , t}^N (x,k),
\end{equation*}
and assume $t' > t \geq s$. We can then write
\begin{equation*}
\begin{aligned}
g_{s,t}^N (x,k) &= G_{t-s}^{N} * \phi (x,k) \\
&\hspace{1cm} + \int_s^T G_{u-s}^{N} * \Big( \mathcal{R}_{m_u^N}^N (g_{u,t}^N) - \mu g_{u,t}^N + \mu \mathcal{H} g_{u,t}^N + uV_R \overline{\overline{m}_u^N r'(\overline{m}_u^N) \overline{\mathcal{H} g_{u,t}^N}} \Big) (x,k) du\\
&\hspace{1cm} - \int_t^T G_{u-s}^{N} * \Big( \mathcal{R}_{m_u^N}^N (\phi) - \mu \phi + \mu \mathcal{H} \phi + uV_R \overline{\overline{m}_u^N r'(\overline{m}_u^N) \overline{\mathcal{H} \phi}} \Big) (x,k) du.
\end{aligned}
\end{equation*}
This transforms the difference into
\begin{equation*}
\begin{aligned}
&\big( g_{s, t'}^N - g_{s,t}^N \big)(x,k)\\
&= G_{t'-s}^{N} * \phi (x,k) - G_{t-s}^{N} * \phi (x,k) \\
&\hspace{1cm} + \int_s^T G_{u-s}^{N} * \Bigg[ \mathcal{R}_{m_u^N} (g_{u,t'}^N) - \mu g_{u,t'}^N + \mu \mathcal{H}g_{u,t'}^N + uV_R \overline{\overline{m}_u^N r'(\overline{m}_u^N) \overline{\mathcal{H}g_{u,t'}^N}} \\
& \hspace{2cm} - \Big( \mathcal{R}_{m_u^N} (g_{u,t}^N) - \mu g_{u,t}^N + \mu \mathcal{H} g_{u,t}^N + uV_R \overline{\overline{m}_u^N r'(\overline{m}_u^N) \overline{\mathcal{H} g_{u,t}^N}} \Big) \Bigg] (x,k) du\\
&\hspace{1cm} + \int_{t}^{t'} G_{u-s}^{N} * \Big( \mathcal{R}_{m_u^N} (\phi) - \mu \phi + \mu \mathcal{H} \phi + uV_R \overline{\overline{m}_u^N r'(\overline{m}_u^N) \overline{\mathcal{H} \phi}} \Big) (x,k) du.
\end{aligned}
\end{equation*}
We continue similar as in the proof of \Cref{lem:semigroupsbounded}
\begin{equation} \label{eq:bigdisplaycontinuityintime}
\begin{aligned}
&\Big\Vert g_{s,t'}^N - g_{s,t}^N  \Big\Vert_p^p\\
&\leq 3 \Big\Vert G_{t'-s}^{N} * \phi - G_{t-s}^{N} * \phi \Big\Vert_p^p \\
&\hspace{1cm} + 3^{p-1} (T-s)^{p-1} \int_s^T \Big\Vert \mathcal{R}_{m_u^N} (g_{u,t'}^N) - \mu g_{u,t'}^N + \mu \mathcal{H}g_{u,t'}^N + uV_R \overline{\overline{m}_u^N r'(\overline{m}_u^N) \overline{\mathcal{H}g_{u,t'}^N}} \\
& \hspace{2cm} - \Big( \mathcal{R}_{m_u^N} (g_{u,t}^N) - \mu g_{u,t}^N + \mu \mathcal{H} g_{u,t}^N + uV_R \overline{\overline{m}_u^N r'(\overline{m}_u^N) \overline{\mathcal{H} g_{u,t}^N}} \Big) \Big\Vert_p^p du\\
&\hspace{1cm} + 3^{p-1} (t' -t)^{p-1} \int_{t}^{t'} \Big\Vert \mathcal{R}_{m_u^N} (\phi) - \mu \phi + \mu \mathcal{H} \phi + uV_R \overline{\overline{m}_u^N r'(\overline{m}_u^N) \overline{\mathcal{H} \phi}} \Big\Vert_p^p du.
\end{aligned}
\end{equation}
The bound on $\Vert \overline{\phi} - \phi \Vert_\infty$ from \eqref{eq:boundphiaverageminusphi}, yields for the first term
\begin{align*}
\Big\Vert G_{t'-s}^N * \phi - G_{t-s}^N * \phi \Big\Vert_p^p
&= \Bigg\Vert \int_t^{t'} G_{u-s}^{N} * \mathcal{L}^N \phi du \Bigg\Vert_p^p \\
&= ( t' - t )^{p-1} \int_t^{t'} \int_{\mathbb{T}^d} \sup_{k \in [0,1]} \big\vert \mathcal{L}^N \phi (x,k) \big\vert^p dx du\\
&\leq C ( t' -t )^{p} Vol (\mathbb{T}^d) \Big( \sup_{k \in [0,1]} \max_{\vert \beta \vert = 2} \Vert \partial_\beta \phi (\cdot, k) \Vert_\infty \Big)^p.
\end{align*}
Both the second and the third summand of \eqref{eq:bigdisplaycontinuityintime} can be bounded similar to the previous lemma resulting in
\begin{align*}
\begin{aligned}
&\Big\Vert g_{s,t'}^N - g_{s,t}^N \Big\Vert_p^p\\
&\leq 3^{p-1} C ( t' - t )^p + 3^{p-1} ( t' -t )^{p} ( u V_R (r_{\max} + n_{max} \Vert r' \Vert_\infty) + 2 \mu ) \big\Vert \phi \big\Vert_p^p\\
&\hspace{0.5cm} +  3^{p-1} (T-s)^{p-1} \big( u V_R (\sup_{z \in \mathbb{T}^d} \Vert r \Vert_\infty + n_{max} \sup_{z \in \mathbb{T}^d} \Vert r_z' \Vert_\infty) + 2 \mu \big)
 \int_s^T \Big\Vert g_{u,t'} - g_{u, t} \Big\Vert_p^p du.
\end{aligned}
\end{align*}
A final application of Gr\"onwall's inequality completes this section.
\end{proof}

\section{Convergence properties of worthy martingale measures integrals} \label{appendix:theorems}
\subsection{Tightness criterion and uniform bound}
\begin{theorem}[{\parencite[Theorem 3.7]{forien}}] \label{theo:tightnessstochasticintegrals}
Let 
\begin{equation*}
U_t^N = \int_{[0,t] \times \mathbb{T}^d \times [0,1]} \psi_{s,t}^N (x,k) M^N (ds\; dx\; dk)
\end{equation*}
be a sequence of stochastic integrals of functions
\begin{equation*}
    \psi^N : \big\{ (s,t) : 0 \leq s \leq t \big\} \times \mathbb{T}^d \times [0,1] \rightarrow \mathbb{R}
\end{equation*}
with respect to a sequence $(M^N)_{N \geq 1}$ of worthy martingale measures. Then, if two conditions, which we specify below, are satisfied, $(U^N)_{N \geq 1} \subset \mathbb{D}(\mathbb{R}_+, \mathbb{R})$ is a tight sequence and for any $T > 0$ there exists constants $C(T), C_1, C_2, C_3 > 0$ such that uniformly in $N$
\begin{equation*}
    \mathbb{E} \Bigg( \sup_{t \in [0,T]} \big\vert U_t^N \big\vert^2 \Bigg) \leq 2 C(T) C_1 (C_2^2 + C_3^2).
\end{equation*}
Here, $C(T)$ depends only on time and $ C_1, C_2$ and $C_3$ are subject to the conditions:
\begin{enumerate}
    \item
    Uniform bound on the dominating measures: for any $0 \leq s \leq t, N \geq 1$ and $\phi \in \mathcal{D} (\mathbb{T}^d \times [0,1])$ holds
    \begin{equation*}
    \blangle D^N , \mathds{1}_{[s,t]} \phi \otimes \phi \brangle \leq C_1 \vert t -s \vert \Vert \phi \Vert_2^2.
    \end{equation*}
    \item 
    Well-behavedness of the integrand: $\psi_{s,t}^N$ is continuous in both time dimensions and for any $0 \leq s \leq s' \leq t \leq t', N \geq 1$, $\psi_{s,t}^N \in \mathcal{D} (\mathbb{T}^d \times [0,1] )$ and for any $p \in [1, 2]$
    \begin{align*}
        \big\Vert \psi_{s,t}^N \big\Vert_p &\leq C_2,\\
        \big\Vert \psi_{s,t'}^N - \psi_{s,t}^N \big\Vert_p &\leq C_3 \vert t' - t \vert,\\
    \end{align*}
\end{enumerate}
\end{theorem}
\subsection{Convergence of finite-dimensional distributions}

\begin{theorem}[{\parencite[Theorem E.1]{forien}}] \label{theo:convergencefinitedistributions}
Let $((\psi^{i,N})_{0 \leq i \leq n}, N \geq 1)$ be a sequence of functions, where $\psi^{i,N} : [0,T] \times \mathbb{T}^d \times [0,1] \rightarrow \mathbb{R}$ and $M^N$ be a sequence of worthy martingale measures. Suppose that the following conditions are satisfied:
\begin{enumerate}
    \item 
    Boundedness of the integrands: for any $N \geq 1, 1 \leq i \leq n$ and $s \in [0,T]$ we have $\psi^{i,N}_s \in \mathcal{D} (\mathbb{T}^d \times [0,1])$ and there exists a bound on the corresponding $L^p$-norms for $p \in [1,2]$ which is uniform in $N$, $i$ and $s$.
    \item 
    Convergence of the integrands: the sequence $((\psi^{i, N})_{0 \leq i \leq n}, N \geq 1)$ converges to limiting functions $\psi^i : [0,T] \times \mathbb{T}^d \times [0,1]$ in the sense that for $p \in \{ 1,2 \}$
    \begin{equation*}
        \lim_{N \to \infty} \sup_{s \in [0,T]} \big\Vert \psi^{i, N}_s - \psi^i_s \big\Vert_p = 0.
    \end{equation*}
    \item 
    Boundedness and convergence of integrators: the sequence $(M^N, N \geq 1)$ satisfies the analog of \Cref{lem:boundondominatingmeasure} and converges in distribution to a martingale measure limit $M$ in $\mathbb{D}(\mathbb{R}_+, \mathcal{D}' (\mathbb{T}^d \times [0,1])$.
\end{enumerate}
Then, as $N \to \infty$,
\begin{align*}
&\Bigg( \int_0^t \blangle \psi^{1,N} , dM^N \brangle ,\ldots, \int_0^t \blangle \psi^{n,N}, dM^N \brangle \Bigg)_{t \in [0,T]} \\
& \hspace{2cm} \rightarrow \Bigg( \int_0^t \blangle \psi^1 , dM \brangle ,\ldots, \int_0^t \blangle \psi^n, dM \brangle \Bigg)_{t \in [0,T]}
\end{align*}
in distribution in $\mathbb{D}([0,T] , \mathbb{R}^n)$.
\end{theorem}

\end{document}